\documentclass[12pt]{article}
\usepackage[]{amsmath,amssymb}
\usepackage{amscd}
\usepackage{latexsym}
\usepackage{cite}

\newtheorem{definition}{Definition}[section]
\newtheorem{theorem}[definition]{Theorem}
\newtheorem{lemma}[definition]{Lemma}
\newtheorem{corollary}[definition]{Corollary}

\newtheorem{example}[definition]{Example}

\newtheorem{note}[definition]{Note}

\newtheorem{proposition}[definition]{Proposition}

\def\Z{\mathbb Z}
\def\F{\mathbb F}

\begin{document}
\title{\bf  
Finite-dimensional irreducible 
$U_q(\mathfrak{sl}_2)$-modules 
from the equitable point of view
}
\author{
Paul Terwilliger}
\date{}

\maketitle
\begin{abstract}
We consider the quantum algebra
$U_q(\mathfrak{sl}_2)$
with $q$ not a root of unity.
We describe the
finite-dimensional
irreducible 
$U_q(\mathfrak{sl}_2)$-modules
from the point of view of the
equitable presentation.

 \end{abstract}

\section{Introduction}
The  quantum universal enveloping
algebra 
$U_q(\mathfrak{sl}_2)$ appears extensively
in the literature; see for example
\cite{chari,jantzen, kassel}.   
In 
\cite{equit} 
the equitable
presentation for
$U_q(\mathfrak{sl}_2)$ was introduced.
This presentation is linked to
tridiagonal pairs of linear transformations
\cite{uqsl2hat,nonnil},
Leonard pairs of linear transformations
\cite{alnajjar},
the
$q$-tetrahedron
algebra
\cite{neubauer2,qtet, qinv, miki},
bidiagonal pairs of linear transformations
\cite{neubauer},
 $Q$-polynomial distance-regular graphs
\cite{drg,drg2,boyd},
Poisson algebras \cite{jordan},
  and the     universal Askey-Wilson algebra \cite{uawe}.
The equitable presentation 
concept has been applied  to
symmetrizable Kac-Moody algebras 
\cite{tersym}
and the Lie algebra 
$\mathfrak{sl}_2$ \cite{bt}.

\medskip
\noindent
In the representation theory of
$U_q(\mathfrak{sl}_2)$, perhaps  the most fundamental
objects are the
finite-dimensional
irreducible 
$U_q(\mathfrak{sl}_2)$-modules with $q$ not a root of
unity. For these objects one desires
 a comprehensive description from the equitable
point of view. Some of the articles mentioned above
contain results in this direction, but  a comprehensive
treatment is lacking.
The goal of
the present paper is to provide this comprehensive treatment.
Our treatment has a  linear algebraic
and geometric flavor.

\medskip
\noindent 
Our treatment is summarized as follows.
Let $\F$ denote a field
 and consider the algebra $U_q(\mathfrak{sl}_2)$  over $\F$.
Let $x,y^{\pm 1},z$ denote the equitable generators for 
 $U_q(\mathfrak{sl}_2)$ 
and let
$n_x, n_y, n_z$ denote their  nilpotent
relatives (formal definitions begin in Section 2).
We retain the notation 
$x,y^{\pm 1},z$ and $n_x,n_y, n_z$ for the corresponding
elements in 
 $U_{q^{-1}}(\mathfrak{sl}_2)$.
We display an $\F$-algebra antiisomorphism 
$\dagger:
 U_q(\mathfrak{sl}_2) \to
 U_{q^{-1}}(\mathfrak{sl}_2)$ that sends
$\xi \mapsto \xi $ and $n_\xi \mapsto -n_\xi$ for
$\xi \in \lbrace x,y,z\rbrace$. 
Fix an integer $d\geq 0$ and
let $V$ denote
an irreducible 
 $U_q(\mathfrak{sl}_2)$-module of type 1 and dimension $d+1$.
Let $V^*$ denote the dual space for $V$ and note that
$V^*$ has dimension $d+1$. Define
a bilinear form
$(\,,\,): V\times V^* \to \F$ such that
$(u,f)=f(u)$ for all $u \in V$ and $f \in V^*$.
We show that $V^*$
becomes a  
 $U_{q^{-1}}(\mathfrak{sl}_2)$-module such that
$(\zeta u, v) = (u,\zeta^\dagger v)$ for all $ u \in V$,  $v \in V^*$,
$\zeta \in 
 U_q(\mathfrak{sl}_2)$.
The 
 $U_{q^{-1}}(\mathfrak{sl}_2)$-module $V^*$ is irreducible 
 of type 1.
We show that on $V$ and $V^*$, each of $x,y,z$ is diagonalizable
with eigenvalues
$\lbrace q^{d-2i}\rbrace_{i=0}^d$.
For $V$ and $V^*$ we display three flags,
six decompositions, and twelve bases. 
We consider (i) how these objects are related to
each other; (ii) 
how these objects are related via the bilinear form;
(iii) 
how these objects are acted upon by
$x,y,z$ and $n_x,n_y,n_z$.
Among the  objects the easiest to
describe are the decompositions, so we begin with these.

\medskip
\noindent 
Each of the six decompositions is an eigenspace decomposition
for 
 one of $x,y,z$.
The corresponding sequence of eigenvalues is
$\lbrace q^{d-2i}\rbrace_{i=0}^d $ or
$\lbrace q^{2i-d}\rbrace_{i=0}^d $.
For each of the six decompositions, the inverted
decomposition is included among the six.
For each of the six decompositions
of $V$, the dual decomposition 
with respect to
$(\, , \,)$ is included among the
 six decompositions for $V^*$.
For $V$ or $V^*$ and 
$\xi \in \lbrace x,y,z\rbrace $ 
we describe
the actions of $\xi$ and $n_\xi$ on the six decompositions.
For these six decompositions
the action of $\xi $ is diagonal on two,
 quasi-lowering on two, and quasi-raising on two.
The action of $n_\xi$ is tridiagonal on two,
lowering on two, and raising on two.

\medskip
\noindent 
Turning to the three flags, 
we show that for $\xi \in \lbrace x,y,z\rbrace$ the subspace
$n^i_\xi V$ has dimension $d-i+1$ for $0\leq i \leq d$
and
$n^{d+1}_\xi V=0$.
Therefore the nested sequence
$\lbrace n^{d-i}_\xi V\rbrace_{i=0}^d$
is a flag on $V$. This gives three flags on 
$V$, and we similarly obtain three flags on
$V^*$.
We show that for $V$ or $V^*$
the three flags are mutually
opposite. These flags are related to the six decompositions
as follows.
For $V$ or $V^*$ let $\lbrace V_i\rbrace_{i=0}^d$ denote
one of the six decompositions. Define
$U_i = V_0 + \cdots + V_i$ for $0 \leq i \leq d$.
We show that the sequence 
$\lbrace U_i\rbrace_{i=0}^d$ is among the three
flags.
To characterize the three flags on $V$,
we show that
for 
$\xi \in \lbrace x,y,z\rbrace$ and
$0 \leq i \leq d+1$, $n^i_\xi V$ is the unique
$(d-i+1)$-dimensional subspace of $V$ that
is invariant under those elements among 
$ x,y,z$ other than $\xi$.
A similar result applies to $V^*$.
We also show that for
$\xi \in \lbrace x,y,z\rbrace$  and
$0 \leq i \leq d+1$,
the subspaces
$n^i_\xi V$ and $n^{d-i+1}_\xi V^*$ are orthogonal
complements with respect to
$(\, , \,)$.

\medskip
\noindent 
Turning to the twelve bases,
each of these bases induces one of the six decompositions.
For each of the twelve bases, the inverted basis is included
among the twelve. For each of the twelve bases for
$V$, the dual basis with respect to 
$(\, , \,)$ is included among the twelve bases for $V^*$.
For $V$ or $V^*$ and each of the twelve bases, we give the matrices that
represent $x,y,z$. Of the resulting three matrices
one is diagonal, one is lower bidiagonal,
and one is upper bidiagonal. 
In each case the sequence  of diagonal entries
is $\lbrace q^{d-2i}\rbrace_{i=0}^d$ or
 $\lbrace q^{2i-d}\rbrace_{i=0}^d$.
In each bidiagonal case
the matrix has constant row sum or constant column sum.
For $V$ or $V^*$ 
and each of the twelve bases we also give the matrices
that represent $n_x,n_y,n_z$.
For $V$ or $V^*$ and each of the twelve bases,
we give the
transition matrix to three other bases among the
twelve. Of the resulting three matrices one is diagonal,
one is lower triangular, and one is
the identity matrix reflected about a vertical axis.

\medskip
\noindent 
Throughout the paper
we employ 
an element in ${\rm End}(V)$ or
 ${\rm End}(V^*)$ 
called 
a rotator. Conjugation by a rotator
induces a cyclic permutation of $x,y,z$. 
These rotators exist by 
\cite[Lemma 7.5]{equit}.
For $V$ or $V^*$ and a rotator $\mathcal R$
we compute the matrices that
represent $\mathcal R$ with respect to the twelve bases.

\medskip
\noindent 
Near the end of the paper
we characterize
$x,y,z$ and
$n_x,n_y,n_z$
in terms of their action on the six decompositions of $V$.
We then characterize 
$U_q(\mathfrak{sl}_2)$ itself in the equitable presentation,
 in terms of bidiagonal triples of linear transformations.
This characterization makes heavy use
of the work of 
Darren Funk-Neubauer
\cite{neubauer} concerning bidiagonal pairs of
linear transformations.

\section{Preliminaries}
Our conventions for the paper are as follows.
An algebra is meant to be associative and have a 1.
A subalgebra has the same 1 as the parent algebra.
Throughout the paper fix an integer $d\geq 0$.
Let $\lbrace u_i\rbrace_{i=0}^d$ denote a  sequence.
We call $u_i$ the {\it $i$th component} of the sequence.
By the {\it inversion} of the sequence 
$\lbrace u_i\rbrace_{i=0}^d$ 
we mean the sequence
$\lbrace u_{d-i}\rbrace_{i=0}^d$.
\noindent Fix a field $\F$.
Let $V$ denote a vector space over $\F$ with 
dimension $d+1$. By a {\it decomposition} of $V$
we mean a sequence 
$\lbrace V_i\rbrace_{i=0}^d$ consisting of 
 one-dimensional subspaces  of $V$ such that
 $V=\sum_{i=0}^d V_i$ (direct sum).
Let $\lbrace V_i\rbrace_{i=0}^d$ denote a decomposition
of $V$. 
For notational convenience 
define $V_{-1}=0$
and $V_{d+1}=0$.
Let ${\rm End}(V)$ denote the $\F$-algebra consisting
of the $\F$-linear maps from $V$ to $V$.
An element
${A \in \rm End}(V)$ 
is called
{\it diagonalizable} whenever $V$ is spanned by the
eigenspaces of $A$.
The map $A$ is called {\it multiplicity-free} whenever
$A$ is diagonalizable, and each eigenspace of $A$ has dimension 1.
Note that $A$ is multiplity-free if and only if $A$ has 
$d+1$ mutually distinct eigenvalues in $\F$.
Assume that $A$ is multiplicity-free, 
and let $\lbrace \theta_i\rbrace_{i=0}^d$
denote an ordering of the eigenvalues of $A$.
For $0 \leq i \leq d$ let $V_i$ denote the eigenspace of $A$
for $\theta_i$. Then the sequence
$\lbrace V_i\rbrace_{i=0}^d$ is a decomposition of $V$.
Let $\lbrace v_i \rbrace_{i=0}^d$ denote a basis
for $V$ and let 
$\lbrace V_i\rbrace_{i=0}^d$ denote
a decomposition of $V$. We say that
$\lbrace v_i \rbrace_{i=0}^d$ {\it induces} 
$\lbrace V_i \rbrace_{i=0}^d$ whenever
$v_i \in V_i$ for $0 \leq i \leq d$.

\noindent 
\begin{definition} 
\label{def:lowering}
\rm
Let $\lbrace V_i \rbrace_{i=0}^d$ denote a
decomposition of $V$. An element
 $\phi \in {\rm End}(V)$ is said to be
 {\it diagonal} on  
$\lbrace V_i \rbrace_{i=0}^d$ 
whenever $\phi V_i \subseteq V_{i}$ for $0 \leq i \leq d$.
The map $\phi$ is said to be
{\it lowering} for
$\lbrace V_i \rbrace_{i=0}^d$ 
whenever $\phi V_i \subseteq V_{i-1}$ for $1 \leq i \leq d$
and $\phi V_0=0$.
The map $\phi$ is said to be {\it quasi-lowering}
for $\lbrace V_i \rbrace_{i=0}^d$ 
whenever $\phi V_i \subseteq V_i+ V_{i-1}$ for $1 \leq i \leq d$
and $\phi V_0\subseteq V_0$.
The map $\phi$ is said to be {\it raising} (resp. {\it quasi-raising})
for $\lbrace V_i \rbrace_{i=0}^d$  whenever $\phi$
is lowering (resp. quasi-lowering) for the 
inversion 
$\lbrace V_{d-i} \rbrace_{i=0}^d$.
\end{definition}


\section{The equitable presentation for
$U_q(\mathfrak{sl}_2)$}

\noindent Fix a nonzero $q \in \F$
such that $q^2 \not=1$.
 For an integer $n$  define
\begin{eqnarray*}
\lbrack n \rbrack = \frac{q^n-q^{-n}}{q-q^{-1}}
\end{eqnarray*}
and for $n\geq 0$  define
\begin{eqnarray*}
\lbrack n \rbrack^! 
=
\lbrack n \rbrack 
\lbrack n-1 \rbrack 
\cdots 
\lbrack 2 \rbrack 
\lbrack 1 \rbrack.
\end{eqnarray*}
We interpret 
$\lbrack 0 \rbrack^!=1$.
We now recall the
quantum algebra 
$U_q(\mathfrak{sl}_2)$. 
We will work with the
equitable presentation
\cite{equit,uawe}.

\begin{definition}
\label{def:equit}  \rm
\cite[Definition~1.1]{equit}
For the $\F$-algebra 
$U_q(\mathfrak{sl}_2)$ the equitable presentation
has
generators
$x,y^{\pm 1},z$ and relations
$yy^{-1}=1$, $y^{-1}y=1$,
\begin{eqnarray}
\label{eq:equit}
\frac{qxy-q^{-1}yx}{q-q^{-1}}=1,
\qquad
\frac{qyz-q^{-1}zy}{q-q^{-1}}=1,
\qquad
\frac{qzx-q^{-1}xz}{q-q^{-1}}=1.
\end{eqnarray}
We call $x,y^{\pm 1}, z$ the {\it equitable generators}
for 
$U_q(\mathfrak{sl}_2)$.
\end{definition}

\noindent In the equations
 (\ref{eq:equit}),
 rearrange terms to find that the
 equitable generators $x,y,z$ of 
$U_q(\mathfrak{sl}_2)$ satisfy
\begin{eqnarray*}
&&q(1-yz) = q^{-1}(1-zy),
\\
&&q(1-zx) = q^{-1}(1-xz),
\\
&&q(1-xy) =  q^{-1}(1-yx).
\end{eqnarray*}

\begin{definition} 
\label{def:nnn}
\rm 
\cite[Definition~5.2]{equit}
Let $n_x, n_y, n_z$ denote the following elements in 
$U_q(\mathfrak{sl}_2)$:
\begin{eqnarray*}
&&n_x = \frac{q(1-yz)}{q-q^{-1}} = \frac{q^{-1}(1-zy)}{q-q^{-1}},
\\
&&n_y = \frac{q(1-zx)}{q-q^{-1}} = \frac{q^{-1}(1-xz)}{q-q^{-1}},
\\
&&n_z = \frac{q(1-xy)}{q-q^{-1}} = \frac{q^{-1}(1-yx)}{q-q^{-1}}.
\end{eqnarray*}
\end{definition}

\begin{lemma}
\label{lem:qcom}
{\rm \cite[Lemma~5.4]{equit}}
The following relations hold in 
$U_q(\mathfrak{sl}_2)$:
\begin{eqnarray*}
&&x n_y = q^2 n_y x, \qquad \qquad
x n_z = q^{-2} n_z x,
\label{eq:com1}
\\
&&y n_z = q^2 n_z y, \qquad \qquad
y n_x = q^{-2} n_x y,
\label{eq:com2}
\\
&&z n_x = q^2 n_x z, \qquad \qquad
z n_y = q^{-2} n_y z.
\label{eq:com3}
\end{eqnarray*}
\end{lemma}

\begin{lemma}
\label{lem:nxynzgen}
{\rm \cite[Lemma~6.4]{uawe}}
The algebra 
$U_q(\mathfrak{sl}_2)$ is generated by $n_x, y^{\pm 1}, n_z$. Moreover
\begin{eqnarray}
x = y^{-1} - q^{-1}(q-q^{-1})n_z y^{-1},
\qquad \qquad 
z = y^{-1} - q(q-q^{-1})n_x y^{-1}.
\label{eq:xz}
\end{eqnarray}
\end{lemma}


\section{Comparing
$U_q(\mathfrak{sl}_2)$ and
$U_{q^{-1}}(\mathfrak{sl}_2)$
}

\noindent In this section
we compare the algebras $U_q(\mathfrak{sl}_2)$ and
$U_{q^{-1}}(\mathfrak{sl}_2)$.
For both algebras we use the same notation
$x,y^{\pm 1},z$ for the equitable generators.

\begin{lemma} 
\label{eq:qiequit}
\rm
The equitable presentation for $U_{q^{-1}}(\mathfrak{sl}_2)$ has
generators $x,y^{\pm 1},z$ and relations
$yy^{-1}=1$, $y^{-1}y=1$,
\begin{eqnarray}
\label{eq:qirel}
\frac{qzy-q^{-1}yz}{q-q^{-1}}=1,
\qquad
\frac{qyx-q^{-1}xy}{q-q^{-1}}=1,
\qquad
\frac{qxz-q^{-1}zx}{q-q^{-1}}=1.
\end{eqnarray}
\end{lemma}
\noindent {\it Proof:} 
In Definition
\ref{def:equit}
replace $q$ by $q^{-1}$ and rearrange terms.
\hfill $\Box$ \\

\begin{corollary} There exists an $\F$-algebra isomorphism 
$U_q(\mathfrak{sl}_2) \to
U_{q^{-1}}(\mathfrak{sl}_2)$ that sends
\begin{eqnarray}
x \mapsto z, \qquad \qquad 
y \mapsto y, \qquad \qquad 
z \mapsto x.
\label{eq:isoaction}
\end{eqnarray}
\end{corollary}
\noindent {\it Proof:}
Compare 
(\ref{eq:equit}) and
(\ref{eq:qirel}).
\hfill $\Box$ \\

\medskip
\noindent We just displayed an isomorphism from
$U_q(\mathfrak{sl}_2)$ to
$U_{q^{-1}}(\mathfrak{sl}_2)$.
Next we display an 
antiisomorphism
from $U_q(\mathfrak{sl}_2)$ to
$U_{q^{-1}}(\mathfrak{sl}_2)$.
An antiisomorphism
is defined as follows.
Given
$\F$-algebras
$\mathcal A$, $\mathcal B$
a map $\sigma : 
\mathcal A \to \mathcal B$  is called an
{\it antiisomorphism of $\F$-algebras} whenever
$\sigma$ is an isomorphism of $\F$-vector spaces
and $(ab)^\sigma=b^\sigma a^\sigma$ for
all $a,b \in \mathcal A$.
An antiisomorphism can be interpreted
as follows.
The $\F$-vector space 
 $\mathcal B$ supports an $\F$-algebra structure
 $\mathcal B^{opp}$  such that
 for all $a,b \in \mathcal B$
the product $ab$ (in $\mathcal B^{opp}$)
is equal to 
$ba$ (\rm in  $\mathcal B$).
A map $\sigma : 
\mathcal A \to \mathcal B$  is an antiisomorphism
of $\F$-algebras if and only if 
$\sigma:
\mathcal A \to \mathcal B^{opp}$  is an
isomorphism of $\F$-algebras.


\begin{proposition}
\label{prop:antiiso}
There exists an antiisomorphism
of $\F$-algebras 
$\dagger:
 U_q(\mathfrak{sl}_2) \to
 U_{q^{-1}}(\mathfrak{sl}_2) $
that sends 
\begin{eqnarray}
x \mapsto x, \qquad \qquad 
y \mapsto y, \qquad \qquad 
z \mapsto z.
\label{eq:antiaction}
\end{eqnarray}
\end{proposition}
\noindent {\it Proof:}
In the presentation for
 $U_{q^{-1}}(\mathfrak{sl}_2) $ from
 Lemma
\ref{eq:qiequit}, reverse the order of multiplication
to get a presentation for
 $U_{q^{-1}}(\mathfrak{sl}_2)^{opp}$ that matches the
 presentation for
 $U_{q}(\mathfrak{sl}_2) $ given in
Defininition
\ref{def:equit}. 
Therefore there exists
an $\F$-algebra isomorphism
$\dagger:
 U_q(\mathfrak{sl}_2) \to
 U_{q^{-1}}(\mathfrak{sl}_2)^{opp}$ that satisfies
(\ref{eq:antiaction}). The result follows in view of
the sentence prior to the proposition statement.
\hfill $\Box$ \\

\noindent 
In Definition
\ref{def:nnn} we defined some
elements
$n_x, n_y, n_z$  in  
 $U_q(\mathfrak{sl}_2)$.
 We retain the
notation 
$n_x, n_y, n_z$ for 
 the corresponding elements in
 $U_{q^{-1}}(\mathfrak{sl}_2)$.

\begin{lemma} 
\label{lem:nxnynz}
The antiisomorphism
$\dagger$ from
Proposition
\ref{prop:antiiso}
 sends
\begin{eqnarray*}
n_x \mapsto -n_x, \qquad \qquad 
n_y \mapsto -n_y, \qquad \qquad 
n_z \mapsto -n_z.
\end{eqnarray*}
\end{lemma}
\noindent {\it Proof:}
Use Definition
\ref{def:nnn}.
\hfill $\Box$ \\


\section{The 
 $U_q(\mathfrak{sl}_2)$-module $V$}

\noindent  We turn our attention to the finite-dimensional
irreducible 
$U_q(\mathfrak{sl}_2)$-modules, for $q$ not a root of unity.
These modules are classified up to isomorphism
in \cite[Section~2.6]{jantzen}. The classification 
shows that for any given finite positive 
dimension
there are two isomorphism classes if 
${\rm Char}(\F)\not=2$,
and one 
isomorphism class if 
${\rm Char}(\F)=2$. As we discuss these modules we
will use the following notational assumptions.

\medskip
\noindent
In this paragraph we make some assumptions
that are in effect until
the end of Section 17.
We assume that $q$ is not a root of unity.
We assume that $V$ is an irreducible 
$U_q(\mathfrak{sl}_2)$-module with dimension $d+1$.
By \cite[Lemma~4.2]{equit} 
the element $y$ is multiplicity-free on $V$. Moreover
by 
 \cite[Lemma~4.2]{equit} 
there exists $\varepsilon \in \lbrace 1,-1\rbrace$
such that the eigenvalues of $y$ on $V$
are $\lbrace \varepsilon q^{d-2i}\rbrace_{i=0}^d$.
The scalar $\varepsilon$ is called the {\it type}
of $V$. Replacing 
 $x,y,z$ by $\varepsilon x,\varepsilon y,\varepsilon z$
 the type becomes 1.
For notational convenience we assume
that $V$ has type 1.

\begin{definition}
\label{def:rotator}
\rm
By a {\it rotator} for $V$ we mean an invertible 
$R \in
 {\rm End}(V)$  
such that on $V$,
\begin{eqnarray}
R x R^{-1} = y,\qquad 
R y R^{-1} = z,\qquad 
R z R^{-1} = x.
\label{eq:rotate}
\end{eqnarray}
\end{definition}

\begin{lemma} 
{\rm \cite[Lemma 7.5]{equit}}
\label{lem:om}
There exists a rotator for $V$.
\end{lemma}
\noindent We comment on the uniqueness of a rotator.

\begin{lemma} Let $R $ denote a rotator for
$V$. Then for 
$\Psi \in
 {\rm End}(V)$   the following are equivalent:
\begin{enumerate}
\item[\rm (i)] 
$\Psi$ is a rotator for $V$;
\item[\rm (ii)] there exists $0 \not=\alpha \in \F$
such that $\Psi=\alpha R $.
\end{enumerate}
\end{lemma}
\noindent {\it Proof:} 
${\rm (i) \Rightarrow (ii)}$ 
The composition $G=\Psi R^{-1}$
commutes
with each of $x,y,z$ and therefore
 everything in  
 $U_q(\mathfrak{sl}_2)$.
Recall that $y$ is multiplicity-free on $V$.
The map $G$ commutes with $y$, so
$G$ leaves invariant the eigenspaces of $y$ on $V$.
 Each of these eigenspaces has dimension one,
and is therefore 
contained in an eigenspace of $G$.
Consequently $G$ is diagonalizable on $V$.
Let $W$ denote an eigenspace of $G$,
and let $\alpha$ denote the corresponding
eigenvalue. Note that $\alpha\not=0$ since
$G$ is invertible.
Since $G$ commutes with everything in
 $U_q(\mathfrak{sl}_2)$, we see that
 $W$ is a $U_q(\mathfrak{sl}_2)$-submodule of $V$.
The 
 $U_q(\mathfrak{sl}_2)$-module $V$ is
 irreducible so $W=V$.
  Therefore
 $G=\alpha I$ so
 $\Psi =\alpha R$.
\\
${\rm (ii) \Rightarrow (i)}$ Clear. 
\hfill $\Box$ \\

\begin{lemma}
\label{lem:mft1}
For each of $x,y,z$ the action on
$V$ is multiplicity-free with eigenvalues
$\lbrace q^{d-2i}\rbrace_{i=0}^d$.
\end{lemma}
\noindent {\it Proof:} The assertion applies to
$y$ by construction. The assertion applies to
$x,z$ in view of Lemma
\ref{lem:om}.
\hfill $\Box$ \\


\section{The 
 $U_{q^{-1}}(\mathfrak{sl}_2)$-module $V^*$}

\noindent Recall the 
 $U_q(\mathfrak{sl}_2)$-module $V$ from Section 5.
The {\it dual space} $V^*$ 
is the vector space over $\F$ consisting of the 
$\F$-linear maps $V\to \F$. The vector spaces $V$ and $V^*$ have
the same dimension.
In this section we have two main goals. First we turn 
 $V^*$ into a $U_{q^{-1}}(\mathfrak{sl}_2)$-module.
 Then we show how  the
 $U_q(\mathfrak{sl}_2)$-module $V$ and the
 $U_{q^{-1}}(\mathfrak{sl}_2)$-module $V^*$ are related.

\begin{definition} 
\label{def:bil}
\rm We define a bilinear
form $(\,,\,): V\times V^* \to \F$ such that
$(u,f)=f(u)$ for all $u \in V$ and $f \in V^*$.
The form
$(\,,\,)$
 is nondegenerate.
\end{definition}

\noindent Vectors $u \in V$ and $v \in V^*$
are called {\it orthogonal} whenever
$(u,v)=0$.

\medskip
\noindent We recall the adjoint map
\cite[p.~227]{roman}.
 Let $A \in 
 {\rm End}(V)$. The {\it adjoint} of $A$,
 denoted
 $A^{adj}$,
is the unique element of  
 ${\rm End}(V^*)$ such that
$(Au,v) = (u,A^{adj} v)$ for all $u \in V$ and $v \in V^*$.
The adjoint map
 ${\rm End}(V)\to
 {\rm End}(V^*)$, $A\mapsto A^{adj}$
is an antiisomorphism of $\F$-algebras.

\medskip
\noindent Recall the antiisomorphism $\dagger:
 U_q(\mathfrak{sl}_2) \to
 U_{q^{-1}}(\mathfrak{sl}_2)$
from Proposition
\ref{prop:antiiso}.

\begin{proposition} 
\label{prop:adjointaction}
There exists a
unique $U_{q^{-1}}(\mathfrak{sl}_2)$-module structure on
$V^*$ such that
\begin{eqnarray}
(\zeta u, v) = (u,\zeta^\dagger v) \qquad u \in V, \quad v \in V^*,
\quad \zeta \in 
 U_q(\mathfrak{sl}_2).
\label{eq:requirement}
\end{eqnarray}
\end{proposition}
\noindent {\it Proof:}
The action of
 $U_q(\mathfrak{sl}_2)$ on $V$
induces an $\F$-algebra homomorphism
 $U_q(\mathfrak{sl}_2) \to 
{\rm End}(V)$. Call this homomorphism $\psi$.
The composition
\[
\begin{CD}
 U_{q^{-1}}(\mathfrak{sl}_2)  @>>\dagger^{-1} >  
 U_q(\mathfrak{sl}_2)  @>>\psi> {\rm End}(V) 
                 @>>adj > {\rm End}(V^*)
                  \end{CD} 
              \]
is an $\F$-algebra homomorphism.
This homomorphism gives 
   $V^*$ a
$U_{q^{-1}}(\mathfrak{sl}_2)$-module structure.
By construction 
the $U_{q^{-1}}(\mathfrak{sl}_2)$-module $V^*$
satisfies the requirement (\ref{eq:requirement}).
We have shown that the desired
$U_{q^{-1}}(\mathfrak{sl}_2)$-module structure exists.
One routinely checks that this structure is unique.
\hfill $\Box$ \\

\noindent In the next two propositions we describe
how the 
 $U_{q}(\mathfrak{sl}_2)$-module $V$ is related to the
 $U_{q^{-1}}(\mathfrak{sl}_2)$-module
$V^*$.

\begin{proposition} 
\label{prop:uquqi}
For all $\zeta \in 
 U_{q}(\mathfrak{sl}_2)$,
$\zeta^\dagger$ acts on $V^*$
as the adjoint of the action of $\zeta $ on
$V$.
\end{proposition}
\noindent {\it Proof:}
By 
(\ref{eq:requirement}) and the definition 
of adjoint from above Proposition
\ref{prop:adjointaction}.
\hfill $\Box$ \\

\begin{proposition} 
\label{prop:howuquqirel}
For 
$u \in V$ and $v \in V^*$,
\begin{eqnarray*}
&&(x u,v) = (u,x v),
\qquad \qquad
(yu,v) = (u,yv),
\qquad \qquad
(zu,v) = (u,zv),
\\
&&(n_x u,v) = -(u,n_x v),
\qquad
(n_y u,v) = -(u,n_y v),
\qquad
(n_z u,v) = -(u,n_z v).
\end{eqnarray*}
\end{proposition}
\noindent {\it Proof:}
Evaluate
(\ref{eq:requirement}) using
Proposition
\ref{prop:antiiso}
and
Lemma
\ref{lem:nxnynz}.
\hfill $\Box$ \\

\noindent
Given a subspace $W$ of $V$ (resp. $V^*$)
let $W^\perp$ denote the set of vectors in $V^*$ (resp. $V$)
that are orthogonal to everything
in $W$. We call $W^\perp$ the {\it orthogonal complement} of
$W$. We have $(W^\perp)^\perp = W$ since $(\,,\,)$ is
nondegenerate. For $W, W^\perp$ the sum of the dimensions
is equal to the common dimension of $V, V^*$ which we recall is
$d+1$.

\begin{lemma} 
\label{lem:wperp}
For a subspace $U \subseteq V$ and
an element $\zeta \in 
 U_q(\mathfrak{sl}_2)$,
$U$ is $\zeta$-invariant if and only if
 $U^\perp$ is $\zeta^\dagger$-invariant.
\end{lemma}
\noindent {\it Proof:}
Use 
(\ref{eq:requirement}).
\hfill $\Box$ \\

\begin{lemma}
The 
 $U_{q^{-1}}(\mathfrak{sl}_2)$-module $V^*$
is irreducible.
\end{lemma}
\noindent {\it Proof:}
Let $W$ denote a 
 $U_{q^{-1}}(\mathfrak{sl}_2)$-submodule of $V^*$.
We show that $W=0$ or $W=V^*$. 
 Consider the orthogonal complement $W^\perp\subseteq V$.
By Lemma
\ref{lem:wperp} 
$W^\perp$
is a 
 $U_q(\mathfrak{sl}_2)$-submodule of $V$.
 The 
 $U_q(\mathfrak{sl}_2)$-module $V$ is irreducible
 so
$W^\perp=V$ or  
$W^\perp=0$.
It follows that
$W=0$ or $W=V^*$. 
\hfill $\Box$ \\

\begin{lemma} 
\label{lem:min}
For $\zeta \in 
 U_q(\mathfrak{sl}_2)$ the following coincide:
\begin{enumerate}
\item[\rm (i)] the minimal polynomial for the action of $\zeta$ on  $V$;
\item[\rm (ii)] the minimal polynomial for the action of $\zeta^\dagger$
on  $V^*$.
\end{enumerate}
\end{lemma}
\noindent {\it Proof:}
Use 
(\ref{eq:requirement}).
\hfill $\Box$ \\

\begin{lemma} 
\label{lem:vvs}
For each of $x,y,z$ the action on
$V^*$ is multiplicty-free with eigenvalues
$\lbrace q^{d-2i}\rbrace_{i=0}^d$. Moreover 
the 
 $U_{q^{-1}}(\mathfrak{sl}_2)$-module $V^*$
has type 1.
\end{lemma} 
\noindent {\it Proof:} 
The first assertion follows
from Lemma
\ref{lem:mft1}
and
Lemma
\ref{lem:min}. The last assertion follows from
the first.
\hfill $\Box$ \\


\noindent Lemma
\ref{lem:vvs} implies that for every
result about $V$ there is a corresponding
result about $V^*$, obtained by
replacing $q$ by $q^{-1}$ and adjusting
the notation.

\section{Six decompositions for $V$ and $V^*$}

\medskip
\noindent We continue to discuss the
$U_q(\mathfrak{sl}_2)$-module $V$ and
the 
$U_{q^{-1}}(\mathfrak{sl}_2)$-module $V^*$.
In this section, for $V$ and $V^*$ we will define six decompositions, denoted
\begin{eqnarray}
&&\lbrack x \rbrack,
\qquad \quad 
\,\lbrack y \rbrack,
\qquad   \quad 
\lbrack z \rbrack,
\label{eq:sixdec1}
\\
&&\lbrack x \rbrack^{inv},
\qquad 
\lbrack y \rbrack^{inv},
\qquad 
\lbrack z \rbrack^{inv}.
\label{eq:sixdec2}
\end{eqnarray}
\noindent We will describe these decompositions from
several points of view.

\begin{definition}
\label{def:sixdec}
\rm 
For $\xi \in \lbrace x,y,z\rbrace $  define
the decomposition 
$\lbrack \xi \rbrack$ of $V$ (resp. $V^*$) as follows.
For $0 \leq i \leq d$ the $i$th component of
$\lbrack \xi \rbrack$ is the eigenspace for $\xi$  with eigenvalue
$q^{d-2i}$ (resp. $q^{2i-d}$).
The inversion of $\lbrack \xi \rbrack $ is denoted
by $\lbrack \xi \rbrack^{inv}$.
\end{definition}


\medskip
\noindent
Let $\lbrace V_i\rbrace_{i=0}^d$ denote a decomposition
of $V$ and let 
 $\lbrace V'_i\rbrace_{i=0}^d$ denote a decomposition
of $V^*$. These decompositions are said to be {\it dual}
 whenever
$(V_i,V'_j)=0$ if $i \not=j$ $(0 \leq i,j\leq d)$.
Each decomposition of $V$ (resp. $V^*$) is dual to
a unique decomposition of $V^*$ (resp. $V$).

\medskip

\begin{lemma} 
\label{lem:dualdec}
For the table below, in each row
 we display a decomposition of $V$ and
 a decomposition of $V^*$. These decompositions
 are dual.
\medskip

\centerline{
\begin{tabular}[t]{c c}
{\rm decomp. of $V$} & 
   {\rm decomp. of $V^*$}
   \\ \hline  \hline
  $\lbrack x \rbrack$   & $\lbrack x \rbrack^{inv} $
  \\ 
  $\lbrack x \rbrack^{inv}$   & $\lbrack x \rbrack $
   \\
\hline
  $\lbrack y \rbrack$   & $\lbrack y \rbrack^{inv}$ 
  \\ 
  $\lbrack y \rbrack^{inv}$   & $\lbrack y \rbrack$ 
  \\ 
\hline
  $\lbrack z \rbrack$   & $\lbrack z \rbrack^{inv}$ 
  \\ 
  $\lbrack z \rbrack^{inv}$   & $\lbrack z \rbrack$ 
     \end{tabular}}
     \medskip

\end{lemma}
\noindent {\it Proof:}
We prove the assertion for 
the first row of the table; for the other rows
the proof is similar.
Pick distinct integers $i,j$ $(0 \leq i,j\leq d)$.
Let $u$ (resp. $v$) denote a vector in the
$i$th (resp. $j$th) component of the 
decomposition $\lbrack x \rbrack $ of $V$
(resp. decomposition $\lbrack x \rbrack^{inv}$ of $V^*$).
We show that $u,v$ are orthogonal.  
By Proposition
\ref{prop:howuquqirel}
$(xu,v)=(u,xv)$.
By Definition
\ref{def:sixdec}
$xu = q^{d-2i}u$
and $xv = q^{d-2j}v$. Note that
$q^{d-2i} \not=q^{d-2j}$ since $q$ is not a root of unity.
By these comments 
 $(u,v)=0$. 
\hfill $\Box$ \\


\noindent We now describe
the actions of
$n_x$, 
$n_y$, 
$n_z$ on the decompositions
(\ref{eq:sixdec1}), 
(\ref{eq:sixdec2})  for $V$ and $V^*$.

\begin{theorem}
\label{thm:nxnynz}
Let $\lbrace V_i\rbrace_{i=0}^d$ denote a
decomposition of $V$ or $V^*$ from among
{\rm (\ref{eq:sixdec1}), 
(\ref{eq:sixdec2})}.
Then for $0 \leq i \leq d$ the actions of 
$n_x$, $n_y$, $n_z$ on $V_i$ are given in the table below.
\medskip

\centerline{
\begin{tabular}[t]{c|c c c}
   $\lbrace V_i\rbrace_{i=0}^d$ & 
   {\rm action of $n_x$ on $V_i$}
   &
   {\rm action of $n_y$ on $V_i$}
   &
   {\rm action of $n_z$ on $V_i$}
   \\ \hline  \hline
    $\lbrack x\rbrack$ & $n_xV_i\subseteq V_{i-1}+V_i+V_{i+1}$
  &
     $n_yV_i = V_{i-1}$   
  &
     $n_zV_i = V_{i+1}$   
\\
    $\lbrack x\rbrack^{inv}$
    & $n_xV_i\subseteq V_{i-1}+V_i+V_{i+1}$
  &
     $n_yV_i = V_{i+1}$   
  &
     $n_zV_i = V_{i-1}$   
\\
\hline
      $\lbrack y\rbrack$ &
     $n_xV_i = V_{i+1}$   
    &
     $n_yV_i \subseteq V_{i-1}+V_i+V_{i+1} $
     &
     $n_z V_i=V_{i-1}$   
\\
      $\lbrack y\rbrack^{inv}$ &
     $n_xV_i = V_{i-1}$   
    &
     $n_yV_i \subseteq V_{i-1}+V_i+V_{i+1} $
     &
     $n_z V_i=V_{i+1}$   
\\
\hline
     $\lbrack z\rbrack$ &
      $n_x V_i = V_{i-1}$
    &
     $n_y V_i = V_{i+1}$   
&
      $n_z V_i\subseteq V_{i-1}+V_i+V_{i+1}$
\\
     $\lbrack z\rbrack^{inv}$ &
      $n_x V_i = V_{i+1}$
    &
     $n_y V_i = V_{i-1}$   
&
      $n_z V_i\subseteq V_{i-1}+V_i+V_{i+1}$
\end{tabular}}
     \medskip

\end{theorem}

\noindent {\it Proof:}
First assume that the given
decomposition is $\lbrack y \rbrack $.
Then $V_i$ is an eigenspace for $y$.
We now use two equations from
Lemma
\ref{lem:qcom}.
Using
 $yn_x=q^{-2}n_xy$
we obtain
$n_x V_i \subseteq V_{i+1}$,
and using
 $yn_z=q^{2}n_z y$
we obtain
$n_z V_i \subseteq V_{i-1}$.
We now show
that $n_x V_i= V_{i+1}$.  Suppose
 $n_x V_i \not= V_{i+1}$. Then $i\leq d-1$ since
$V_{d+1}=0$, and now
 $n_xV_i=0$ since $V_{i+1}$ has dimension one.
By our comments so far the sum
$\sum_{j=0}^i V_j$ is invariant under
each of $n_x,y,n_z$. 
By this and
 Lemma
\ref{lem:nxynzgen}
the sum 
$\sum_{j=0}^i V_j$
is a 
$U_q(\mathfrak{sl}_2)$-submodule of $V$.
Since $0 \leq i \leq d-1$ the sum
$\sum_{j=0}^i V_j$ is nonzero and properly contained in $V$.
This contradicts the fact that the
$U_q(\mathfrak{sl}_2)$-module $V$ is irreducible.
Therefore 
 $n_x V_i= V_{i+1}$. 
One similarly shows 
 $n_z V_i= V_{i-1}$.
Now consider the action of
$n_y$ on $V_i$.
By Definition
\ref{def:nnn}
the element $n_y$ is a scalar multiple of
$1-zx$.
By
(\ref{eq:xz})
and our comments so far we have
$zV_i \subseteq V_i+V_{i+1}$ and
$xV_i \subseteq V_i+V_{i-1}$.
Therefore
$n_yV_i \subseteq V_{i-1}+V_i+V_{i+1}$.
We have verified our assertions for
the decomposition $\lbrack y \rbrack$.
For the decomposition 
$\lbrack y \rbrack^{inv}$
our assertions hold by the meaning of inversion.
For the remaining decompositions in the table
our assertions follow from
Lemma
\ref{lem:om}.
\hfill $\Box$ \\


\noindent We now describe
the actions of $x$, $y$, $z$ on the decompositions 
(\ref{eq:sixdec1}), 
(\ref{eq:sixdec2})  for $V$.

\begin{theorem}
\label{thm:xyzV}
Let $\lbrace V_i\rbrace_{i=0}^d$ denote a
decomposition of $V$ from among
{\rm (\ref{eq:sixdec1}), (\ref{eq:sixdec2})}. 
Then for $0 \leq i \leq d$ the actions of 
$x$, $y$, $z$ on $V_i$ are given in the table below.
\medskip

\centerline{
\begin{tabular}[t]{c|c c c}
    $\lbrace V_i\rbrace_{i=0}^d$ & 
   {\rm action of $x$ on $V_i$}
   &
   {\rm action of $y$ on $V_i$}
   &
   {\rm action of $z$ on $V_i$}
   \\ \hline  \hline
    $\lbrack x\rbrack$ & $(x-q^{d-2i}I)V_i=0$
  &
     $(y-q^{2i-d}I)V_i = V_{i+1}$   
  &
     $(z-q^{2i-d}I)V_i = V_{i-1}$   
   \\
    $\lbrack x\rbrack^{inv}$ & $(x-q^{2i-d}I)V_i=0$
  &
     $(y-q^{d-2i}I)V_i = V_{i-1}$   
  &
     $(z-q^{d-2i}I)V_i = V_{i+1}$   
\\
\hline
      $\lbrack y\rbrack$ &
     $(x-q^{2i-d}I)V_i = V_{i-1}$   
    &
     $(y-q^{d-2i}I)V_i =0 $
     &
     $(z-q^{2i-d}I)V_i = V_{i+1}$   
     \\
      $\lbrack y\rbrack^{inv}$ &
     $(x-q^{d-2i}I)V_i = V_{i+1}$   
    &
     $(y-q^{2i-d}I)V_i =0 $
     &
     $(z-q^{d-2i}I)V_i = V_{i-1}$   
\\
\hline
     $\lbrack z\rbrack$ &
      $(x-q^{2i-d}I)V_i = V_{i+1}$
    &
     $(y-q^{2i-d}I)V_i = V_{i-1}$   
&
      $(z-q^{d-2i}I)V_i=0$
\\
     $\lbrack z\rbrack^{inv}$ &
      $(x-q^{d-2i}I)V_i = V_{i-1}$
    &
     $(y-q^{d-2i}I)V_i = V_{i+1}$   
&
      $(z-q^{2i-di}I)V_i=0$
\end{tabular}}
     \medskip

\end{theorem}
\noindent {\it Proof:}
First assume that the given decomposition
is 
$\lbrack y \rbrack$. 
By construction $(y-q^{d-2i}I)V_i=0$.
By Theorem
\ref{thm:nxnynz}
we have
$n_xV_i = V_{i+1}$
and
$n_zV_i = V_{i-1}$.
Now using
(\ref{eq:xz}),
\begin{eqnarray*}
&&(x - q^{2i-d}I)V_i = 
(x - y^{-1})V_i = 
n_z y^{-1} V_i
=
n_z V_i 
= 
V_{i-1},
\\
&&(z - q^{2i-d}I)V_i = 
(z - y^{-1})V_i = 
n_x y^{-1} V_i
=
n_xV_i
= 
V_{i+1}.
\end{eqnarray*}
We have verified our assertions
for the decomposition 
$\lbrack y \rbrack$. 
For the decomposition
$\lbrack y \rbrack^{inv}$
 our assertions 
 follow from 
the meaning of inversion.
For the remaining decompositions in the table
our assertions follow from
Lemma
\ref{lem:om}.
\hfill $\Box$ \\

\noindent We now describe
the actions of $x$, $y$, $z$ on the decompositions 
(\ref{eq:sixdec1}), 
(\ref{eq:sixdec2})  for $V^*$.

\begin{theorem}
\label{thm:xyzVs}
Let $\lbrace V_i\rbrace_{i=0}^d$ denote a
decomposition of $V^*$ from among
{\rm (\ref{eq:sixdec1}), 
(\ref{eq:sixdec2})}.
Then for $0 \leq i \leq d$ the actions of 
$x$, $y$, $z$ on $V_i$ are given in the table below.
\medskip

\centerline{
\begin{tabular}[t]{c|c c c}
   $\lbrace V_i\rbrace_{i=0}^d $ & 
   {\rm action of $x$ on $V_i$}
   &
   {\rm action of $y$ on $V_i$}
   &
   {\rm action of $z$ on $V_i$}
   \\ \hline  \hline
    $\lbrack x\rbrack$ & $(x-q^{2i-d}I)V_i=0$
  &
     $(y-q^{d-2i}I)V_i = V_{i+1}$   
  &
     $(z-q^{d-2i}I)V_i = V_{i-1}$   
   \\ 
    $\lbrack x\rbrack^{inv}$ & $(x-q^{d-2i}I)V_i=0$
  &
     $(y-q^{2i-d}I)V_i = V_{i-1}$   
  &
     $(z-q^{2i-d}I)V_i = V_{i+1}$   
\\
\hline
      $\lbrack y\rbrack$ &
     $(x-q^{d-2i}I)V_i = V_{i-1}$   
    &
     $(y-q^{2i-d}I)V_i =0 $
     &
     $(z-q^{d-2i}I)V_i = V_{i+1}$   
     \\
      $\lbrack y\rbrack^{inv}$ &
 $(x-q^{2i-d}I)V_i = V_{i+1}$   
    &
     $(y-q^{d-2i}I)V_i =0 $
     &
     $(z-q^{2i-d}I)V_i = V_{i-1}$   
\\
\hline
     $\lbrack z\rbrack$ &
      $(x-q^{d-2i}I)V_i = V_{i+1}$
    &
     $(y-q^{d-2i}I)V_i = V_{i-1}$   
&
      $(z-q^{2i-d}I)V_i=0$
\\
     $\lbrack z\rbrack^{inv}$ &
      $(x-q^{2i-d}I)V_i = V_{i-1}$
    &
     $(y-q^{2i-d}I)V_i = V_{i+1}$   
&
      $(z-q^{d-2i}I)V_i=0$
\end{tabular}}
     \medskip

\end{theorem}
\noindent {\it Proof:}
In Theorem
\ref{thm:xyzV}
replace $q$ by $q^{-1}$.
\hfill $\Box$ \\

\noindent We now give some characterizations of
the decomposition 
$\lbrack y \rbrack$; similar characterizations
apply to the
other 
decompositions from among 
(\ref{eq:sixdec1}), 
(\ref{eq:sixdec2}).

\begin{lemma} 
\label{lem:3ver} 
Referring to  $V$ or $V^*$,
the following coincide for $0 \leq i \leq d$:
\begin{enumerate}
\item[\rm (i)]
the $i$th component of the decomposition
$\lbrack y \rbrack $;
\item[\rm (ii)]
 $n^i_x {\rm Ker} (n_z)$;
\item[\rm (iii)]
 $n^{d-i}_z {\rm Ker} (n_x)$.
\end{enumerate}
\end{lemma}
\noindent {\it Proof:} 
Use Theorem
\ref{thm:nxnynz}.
\hfill $\Box$ \\

\begin{lemma} 
\label{lem:5parts}
Let $\lbrace V_i\rbrace_{i=0}^d$ denote a decomposition
of $V$ or $V^*$. Then the following are equivalent:
\begin{enumerate}
\item[\rm (i)] 
$\lbrace V_i\rbrace_{i=0}^d$ is equal to
$\lbrack y \rbrack$;
\item[\rm (ii)] 
$n_zV_0=0$ and
$n_xV_i\subseteq  V_{i+1}$ for $0 \leq i \leq d-1$;
\item[\rm (iii)] 
$n_zV_0=0$ and $n^i_xV_0\subseteq V_i$ for $0 \leq i \leq d$;
\item[\rm (iv)] 
$n_xV_d=0$ and
$n_zV_i\subseteq V_{i-1}$ for $1 \leq i \leq d$;
\item[\rm (v)] 
$n_xV_d=0$ and $n^{d-i}_zV_d\subseteq V_i$ for $0 \leq i \leq d$.
\end{enumerate}
\end{lemma}
\noindent {\it Proof:} 
${\rm (i) \Rightarrow (ii)}$ By Theorem
\ref{thm:nxnynz}.
\\
${\rm (ii) \Rightarrow (iii)}$ Clear.
\\
${\rm (iii) \Rightarrow (i)}$
We invoke Lemma
\ref{lem:3ver}(i),(ii).
For $0 \leq i \leq d$ we have
 $n^i_x {\rm Ker} (n_z) \subseteq V_i$.
In this inclusion each side has dimension one
so we have equality.
\\
\noindent
${\rm (i) \Rightarrow (iv)}$ By Theorem
\ref{thm:nxnynz}.
\\
${\rm (iv) \Rightarrow (v)}$ Clear.
\\
${\rm (v) \Rightarrow (i)}$
We invoke Lemma
\ref{lem:3ver}(i),(iii).
For $0 \leq i \leq d$ we have
 $n^{d-i}_z {\rm Ker} (n_x) \subseteq V_i$.
In this inclusion each side has dimension one
so we have equality.
\hfill $\Box$ \\

\begin{lemma}
\label{lem:raiselower}
Let $\lbrace V_i\rbrace_{i=0}^d$ denote a decomposition
of $V$ or $V^*$. Then 
$\lbrace V_i\rbrace_{i=0}^d$ is equal to 
$\lbrack y \rbrack$ if and only if both
\begin{enumerate}
\item[\rm (i)]
$n_x$ is raising for
$\lbrace V_i\rbrace_{i=0}^d$;
\item[\rm (ii)]
$n_z$ is lowering for
$\lbrace V_i\rbrace_{i=0}^d$.
\end{enumerate}
\end{lemma}
\noindent {\it Proof:}
Use parts (i), (ii), (iv) of Lemma
\ref{lem:5parts}.
\hfill $\Box$ \\

\begin{lemma}
Let $\lbrace V_i\rbrace_{i=0}^d$ denote
a decomposition of $V$ or $V^*$. Then 
$\lbrace V_i\rbrace_{i=0}^d$ is equal to 
$\lbrack y \rbrack$ if and only if the following hold:
\begin{enumerate}
\item[\rm (i)]
$x$ is quasi-lowering for 
$\lbrace V_i\rbrace_{i=0}^d$;
\item[\rm (ii)]
$y$ is diagonal for 
$\lbrace V_i\rbrace_{i=0}^d$;
\item[\rm (iii)]
$z$ is quasi-raising for
$\lbrace V_i\rbrace_{i=0}^d$.
 \end{enumerate}
\end{lemma}
\noindent {\it Proof:}
$ (\Rightarrow)$ By Theorem
\ref{thm:xyzV} and
Theorem \ref{thm:xyzVs}.
\\
$ (\Leftarrow)$
We invoke Lemma
\ref{lem:5parts}(i),(ii). The subspace $V_0$
is invariant under  $x$ and $y$.
The element $n_z$ is a scalar multiple of
$1-xy$, so $V_0$ is invariant under $n_z$.
But $n_z$ is nilpotent and $V_0$ has dimension one,
so 
$n_zV_0=0$.
Similarly 
$n_xV_d=0$.
For $0 \leq i \leq d-1$
we have the inclusions
$zV_i\subseteq V_i+V_{i+1}$,
$yV_i\subseteq V_i$,
$yV_{i+1}\subseteq V_{i+1}$.
Therefore
$yzV_i\subseteq V_i+V_{i+1}$.
The element $n_x$ is a scalar multiple of
$1-yz$, so
$n_x V_i\subseteq V_i+V_{i+1}$.
But $n_x$ is nilpotent and $V_i$ has dimension one,
so in fact 
$n_x V_i\subseteq V_{i+1}$.
Now by  Lemma
\ref{lem:5parts}(i),(ii) the sequence
$\lbrace V_i\rbrace_{i=0}^d$ is equal to 
$\lbrack y \rbrack$.
\hfill $\Box$ \\

\section{Three flags for $V$ and $V^*$}

\noindent We continue to discuss the
$U_q(\mathfrak{sl}_2)$-module $V$ and
the 
$U_{q^{-1}}(\mathfrak{sl}_2)$-module $V^*$.
In this section we consider these modules using the notion of
a flag. Before we get into the details, we comment on the
notation. We will be discussing 
a number of
results that apply to both $V$ and $V^*$. To simplify the
notation we will focus on $V$;
it is understood that similar results 
hold for $V^*$.
By a {\it flag} on $V$ we mean a sequence
$\lbrace U_i\rbrace_{i=0}^d$ of subspaces for $V$
such that
$U_{i-1} \subseteq U_i$ for $1 \leq i \leq d$
and $U_i$ has dimension $i+1$ for $0 \leq i \leq d$.
For the above flag we have $U_d=V$.
Given a decomposition 
$\lbrace V_i \rbrace_{i=0}^d$ 
of $V$ we construct
a flag on $V$ as follows.
 Define
$U_i= V_0+\cdots + V_i$ for $0 \leq i \leq d$.
Then the sequence 
$\lbrace U_i \rbrace_{i=0}^d$ is a flag 
on $V$. This flag is said to be {\it induced} by
the decomposition $\lbrace V_i \rbrace_{i=0}^d$.
Let 
$\lbrace U_i \rbrace_{i=0}^d$ and
$\lbrace U'_i \rbrace_{i=0}^d$ denote flags on
$V$. These flags are called {\it opposite} whenever
$U_i\cap U'_j=0$ if $i+j<d$  $(0 \leq i,j \leq d)$.
The flags
$\lbrace U_i \rbrace_{i=0}^d$ and
$\lbrace U'_i \rbrace_{i=0}^d$ are opposite
if and only if
there exists a decomposition
$\lbrace V_i \rbrace_{i=0}^d$ of $V$
that induces 
$\lbrace U_i \rbrace_{i=0}^d$ and
whose inversion induces 
$\lbrace U'_i \rbrace_{i=0}^d$.
 In this case
$V_i = U_i\cap U'_{d-i}$ for $0 \leq i \leq d$
\cite[Section~7]{ter24}.

\begin{lemma} 
\label{lem:nnnpre}
The following holds for $0 \leq i \leq d+1$.
\begin{enumerate}
\item[\rm (i)]
$n^i_xV$ is the sum of components
$i,i+1,\ldots, d$ of the decomposition
$\lbrack y \rbrack$
and the sum of
 components
$0,1,\ldots, d-i$ of the decomposition
$\lbrack z \rbrack$.
\item[\rm (ii)]
$n^i_yV$ is the sum of components
$i,i+1,\ldots, d$ of the decomposition
$\lbrack z \rbrack$
and the sum of
 components
$0,1,\ldots, d-i$ of the decomposition
$\lbrack x \rbrack$.
\item[\rm (iii)]
$n^i_zV$ is the sum of components
$i,i+1,\ldots, d$ of the decomposition
$\lbrack x \rbrack$
and the sum of
 components
$0,1,\ldots, d-i$ of the decomposition
$\lbrack y \rbrack$.
\end{enumerate}
\end{lemma}
\noindent {\it Proof:}
(i) 
By construction $V$ is the direct sum
of the components of $\lbrack y \rbrack$.
By Theorem
\ref{thm:nxnynz},
for this decomposition
$n_x$ sends component $j$
onto component $j+1$  
for $0 \leq j \leq d-1$. Moreover
$n_x$ sends component $d$ to zero.
By these comments 
$n^i_xV$ is the sum of components
$i,i+1,\ldots, d$ for 
$\lbrack y \rbrack$. 
We have verified our assertion about
$\lbrack y \rbrack $.
Our assertion about
$\lbrack z \rbrack $ is similarly verified.
\\
\noindent (ii), (iii) Apply
Lemma
\ref{lem:om}.
\hfill $\Box$ \\

\noindent The next three lemmas
 follow routinely from Lemma
\ref{lem:nnnpre}.

\begin{lemma} 
\label{lem:dimflag}
Pick $\xi \in \lbrace x,y,z\rbrace $.
Then 
$n^i_\xi V$ has dimension $d-i+1$ for $0 \leq i \leq d$. Moreover
$n^{d+1}_\xi V=0$.
\end{lemma}

\begin{lemma} 
\label{lem:3flags}
Each of the sequences
\begin{eqnarray}
\lbrace n^{d-i}_xV\rbrace_{i=0}^d,\qquad
\lbrace n^{d-i}_yV\rbrace_{i=0}^d,\qquad
\lbrace n^{d-i}_zV\rbrace_{i=0}^d
\label{eq:flags}
\end{eqnarray}
is a flag on $V$.
\end{lemma}

\begin{lemma}
\label{lem:decindflag}
For each row in the table below, we
give a  decomposition of $V$ along with the induced 
flag on $V$.

\medskip

\centerline{
\begin{tabular}[t]{c c}
{\rm decomp. of $V$ } & 
   {\rm induced flag on $V$}
   \\ \hline  \hline
  $\lbrack x \rbrack$   & 
  $\lbrace n_y^{d-i}V \rbrace_{i=0}^d $
  \\ 
  $\lbrack x \rbrack^{inv}$   & 
  $\lbrace n_z^{d-i}V \rbrace_{i=0}^d $
   \\
\hline
  $\lbrack y \rbrack$   & 
  $\lbrace n_z^{d-i}V \rbrace_{i=0}^d $
  \\ 
  $\lbrack y \rbrack^{inv}$   & 
  $\lbrace n_x^{d-i}V \rbrace_{i=0}^d $
  \\ 
\hline
  $\lbrack z \rbrack$   & 
  $\lbrace n_x^{d-i}V \rbrace_{i=0}^d $
  \\ 
  $\lbrack z \rbrack^{inv}$   & 
  $\lbrace n_y^{d-i}V \rbrace_{i=0}^d $
     \end{tabular}}
     \medskip

\end{lemma}

\begin{lemma}
The three flags 
{\rm (\ref{eq:flags})}
are mutually opposite.
\end{lemma}
\noindent {\it Proof:}
This follows from
Lemma
\ref{lem:decindflag} and the comments
about opposite flags from above Lemma
\ref{lem:nnnpre}.
\hfill $\Box$ \\

\begin{lemma}
\label{lem:decompdesc}
For each row of the table below, we give a
decomposition of $V$ along with its
$i$th component 
for $0 \leq i \leq d$.

\medskip

\centerline{
\begin{tabular}[t]{c c}
{\rm decomp. of $V$ } & 
   {\rm $i$th component}
   \\ \hline  \hline
  $\lbrack x \rbrack$   & 
  $n_y^{d-i}V\cap n_z^iV $
  \\ 
  $\lbrack x \rbrack^{inv}$   & 
  $n_y^{i}V\cap n_z^{d-i}V $
   \\
\hline
  $\lbrack y \rbrack$   & 
  $n_z^{d-i}V\cap n_x^iV $
  \\ 
  $\lbrack y \rbrack^{inv}$   & 
  $n_z^{i}V\cap n_x^{d-i}V $
   \\
\hline
  $\lbrack z \rbrack$   & 
  $n_x^{d-i}V\cap n_y^iV $
  \\ 
  $\lbrack z \rbrack^{inv}$   & 
  $n_x^{i}V\cap n_y^{d-i}V $
     \end{tabular}}
     \medskip

\end{lemma}
\noindent {\it Proof:}
Use Lemma
\ref{lem:nnnpre}.
\hfill $\Box$ \\


\begin{lemma}
\label{lem:kernel} 
 Pick $\xi \in \lbrace x,y,z\rbrace $.
Then for $0 \leq i \leq d+1$
the subspace 
$n_\xi^{i}V$ is the  
 kernel of
$n_\xi^{d-i+1}$ on $V$.
\end{lemma}
\noindent {\it Proof:}
Use Theorem
\ref{thm:nxnynz} and
 Lemma
\ref{lem:nnnpre}.
\hfill $\Box$ \\

\begin{lemma}
\label{lem:yzinv}
Pick $\xi \in \lbrace x,y,z\rbrace$.
Then for $0 \leq i \leq d+1$, 
$n^i_\xi V$ is the unique $(d-i+1)$-dimensional subspace
of $V$ that is invariant under 
those elements among
$ x,y,z$ other than $\xi$.
\end{lemma}
\noindent {\it Proof:}
By Lemma \ref{lem:om}, 
we may assume without loss that 
$\xi=x$.
The subspace 
$n^i_xV$ has dimension $d-i+1$ by
Lemma
\ref{lem:dimflag}.
The subspace 
$n^i_xV$ is invariant under $y,z$ 
by Lemma
\ref{lem:nnnpre}(i).
Let $W$ denote a $(d-i+1)$-dimensional subspace of
$V$ that is invariant under $y,z$. 
We show
that $W=n^i_xV$. First assume
$i=d+1$. Then $W=0=n_x^{d+1}V$.
Next assume
$i \leq d$, so that $W\not=0$.
Let $\lbrace V_j\rbrace_{j=0}^d$ denote the 
decomposition $\lbrack y \rbrack$ of $V$.
Note that $y$ is diagonalizable on $W$, since
 $y$ is diagonalizable on $V$ and
 $W$ is
$y$-invariant.
 Therefore $W$ is spanned by the eigenspaces of
 $y$ on $W$.
 Consequently 
 $W=\sum_{j\in S} V_j$ where
 $S = \lbrace j |0 \leq j \leq d,\;V_j \subseteq W\rbrace$.
The subspace $W$ is invariant under $n_x$,
since 
$n_x$ is a scalar multiple of $1-yz$.
Recall from  
Theorem 
\ref{thm:nxnynz} that $n_xV_j=V_{j+1}$ for
$0 \leq j \leq d-1$.
By these comments $j \in S $ implies $j+1 \in S$ for
$0 \leq j \leq d-1$.
The set $S$ is nonempty since $W\not=0$.
Therefore there exists
an integer $t$ $(0 \leq t \leq d)$ such that
$S=\lbrace t,t+1, \ldots, d \rbrace$.
In other words
$W=\sum_{j=t}^d V_j$.
Considering the dimension $t=i$.
Now using
Lemma
\ref{lem:nnnpre}(i) we find
$n^i_xV = \sum_{j=i}^d V_j=
W$.
\hfill $\Box$ \\

\begin{lemma}
\label{lem:dual}
Pick $\xi \in \lbrace x,y,z\rbrace $.
Then for $0 \leq i \leq d+1$ 
the following are orthogonal complements
 with respect to
the bilinear form $(\,,\,)$: 
\begin{eqnarray*}
n_\xi^iV, \qquad \qquad 
n_\xi^{d-i+1}V^*.
\end{eqnarray*}
\end{lemma}
\noindent {\it Proof:}
Combine
Lemma \ref{lem:dualdec}
and
Lemma
\ref{lem:nnnpre}.
\hfill $\Box$ \\


\section{Twelve bases for $V$ and $V^*$}

\noindent 
We continue to work with the 
$U_q(\mathfrak{sl}_2)$-module $V$
and the
$U_{q^{-1}}(\mathfrak{sl}_2)$-module $V^*$.
In this section, for $V$ and $V^*$ we define twelve bases, denoted
\begin{eqnarray}
&&
\lbrack x \rbrack_{row}, \quad 
\lbrack x \rbrack_{col}, \quad  
\lbrack x \rbrack^{inv}_{row}, \quad 
\lbrack x \rbrack^{inv}_{col},
\label{eq:basisr1}
\\
&&
\lbrack y \rbrack_{row}, \quad 
\lbrack y \rbrack_{col}, \quad  
\lbrack y \rbrack^{inv}_{row}, \quad 
\lbrack y \rbrack^{inv}_{col},
\label{eq:basisr2}
\\
&&
\lbrack z \rbrack_{row}, \quad 
\lbrack z \rbrack_{col}, \quad  
\lbrack z \rbrack^{inv}_{row}, \quad 
\lbrack z \rbrack^{inv}_{col}.
\label{eq:basisr4}
\end{eqnarray}

\noindent We will describe how these bases are related to
each other and the decompositions
(\ref{eq:sixdec1}),
(\ref{eq:sixdec2}).
Before we define
(\ref{eq:basisr1})--(\ref{eq:basisr4})
we have some comments.
By Lemma
\ref{lem:dimflag},
for $\xi \in \lbrace x,y,z\rbrace$
the vector spaces $n_\xi^dV$
and $n_\xi^dV^*$ have dimension one. In the next four lemmas
we clarify the meaning of these spaces.

\begin{lemma}
\label{lem:meaning}
The following {\rm (i)--(iii)} hold:
\begin{enumerate}
\item[\rm (i)] $n^d_xV$ is the eigenspace for $y$ (resp. $z$)  on  $V$
with eigenvalue $q^{-d}$ (resp. $q^d$).
\item[\rm (ii)] $n^d_yV$ is the eigenspace for $z$ (resp. $x$)  on $V$
with eigenvalue $q^{-d}$ (resp. $q^d$).
\item[\rm (iii)] $n^d_zV$ is the eigenspace for $x$ (resp. $y$) on $V$
with eigenvalue $q^{-d}$ (resp. $q^d$).
\end{enumerate}
\end{lemma}
\noindent {\it Proof:}
In Lemma
\ref{lem:nnnpre} set $i=d$
and use
Definition
\ref{def:sixdec}.
\hfill $\Box$ \\

\begin{lemma}
The following {\rm (i)--(iii)} hold:
\begin{enumerate}
\item[\rm (i)] $n^d_xV^*$ is the eigenspace for $y$ (resp. $z$) on $V^*$
with eigenvalue $q^{d}$ (resp. $q^{-d}$).
\item[\rm (ii)] $n^d_yV^*$ is the eigenspace for $z$ (resp. $x$) on $V^*$
with eigenvalue $q^{d}$ (resp. $q^{-d}$).
\item[\rm (iii)] $n^d_zV^*$ is the eigenspace for $x$ (resp. $y$) on $V^*$
with eigenvalue $q^{d}$ (resp. $q^{-d}$).
\end{enumerate}
\end{lemma}
\noindent {\it Proof:}
In Lemma
\ref{lem:meaning}
replace $V$ by $V^*$ and $q$ by $q^{-1}$.
\hfill $\Box$ \\

\begin{lemma}
\label{lem:ndvchar2}
The following hold for $\xi \in \lbrace x,y,z \rbrace$:
\begin{enumerate}
\item[\rm (i)]
$n^d_\xi V$ is the unique common eigenspace on $V$
for the two elements among 
$ x,y,z$ other than $\xi$.
\item[\rm (ii)]
$n^d_\xi V^*$ is the unique common eigenspace on $V^*$
for the two elements among
$x,y,z$ other than $\xi$.
\end{enumerate}
\end{lemma}
\noindent {\it Proof:} 
(i) By Lemma
\ref{lem:yzinv}
and since each of $x,y,z$ is multiplicity-free on $V$.
\\
\noindent (ii) Similar to the proof of (i).
\hfill $\Box$ \\

\begin{lemma}
\label{lem:ndvasker}
The following hold for $\xi \in \lbrace x,y,z \rbrace$:
\begin{enumerate}
\item[\rm (i)]
$n^d_\xi V$ is the kernel of $n_\xi$ on $V$.
\item[\rm (ii)]
$n^d_\xi V^*$ is the kernel of $n_\xi$ on $V^*$.
\end{enumerate}
\end{lemma}
\noindent {\it Proof:} To obtain part (i) 
set $i=d$ in Lemma
\ref{lem:kernel}.
Part (ii)  is similarly obtained.
\hfill $\Box$ \\



\begin{definition}
\label{def:urow}
\rm
Pick $\xi \in \lbrace x,y,z\rbrace$.
A basis
$\lbrace v_i \rbrace_{i=0}^d$ 
for $V$ is said to be
{\it $\lbrack \xi \rbrack_{row}$}
whenever:
\begin{enumerate}
\item[\rm (i)] For $0 \leq i \leq d$ the vector
$v_i$ is contained in component $i$ of the decomposition
$\lbrack \xi \rbrack$;
\item[\rm (ii)] $\sum_{i=0}^d v_i \in n^d_\xi V$.
\end{enumerate}
A 
$\lbrack \xi \rbrack_{row}$ 
basis for $V^*$ 
is similarly defined, with $V$ replaced by
$V^*$ 
in (ii) above.
By a 
$\lbrack \xi \rbrack^{inv}_{row}$ basis we mean
the inversion of a 
$\lbrack \xi \rbrack_{row}$ basis.
\end{definition}

\noindent Consider the bases for $V$ and $V^*$ from 
Definition \ref{def:urow}.
Shortly we will discuss the existence and
uniqueness of these bases.

\begin{lemma} 
\label{lem:rowbasispre}
Consider the decomposition 
$\lbrack y \rbrack$ of $V$.
For $0 \leq i \leq d$ let $v_i$ denote a vector in
the $i$th component. Then the following 
{\rm (i)--(v)} 
are equivalent:
\begin{enumerate}
\item[\rm (i)] 
$\sum_{i=0}^d v_i \in n^d_yV$;
\item[\rm (ii)]
$(z-q^{2i-d}) v_i = (q^{-d}-q^{2i+2-d})
v_{i+1}$ for $0 \leq i \leq d-1$;
\item[\rm (iii)]
 $n_x v_i = q^{-i}\lbrack i+1\rbrack 
v_{i+1}$ for $0 \leq i \leq d-1$;
\item[\rm (iv)]
$(x-q^{2i-d}) v_i = (q^d-q^{2i-2-d})v_{i-1}$
 for $1 \leq i \leq d$;
\item[\rm (v)]
$n_z v_i = -q^{d-i}\lbrack d-i+1\rbrack
v_{i-1}$ for $1 \leq i \leq d$.
\end{enumerate}
Now assume that {\rm (i)--(v)} hold.
Then 
$\lbrace v_i\rbrace_{i=0}^d$ are all zero or all nonzero.
\end{lemma}
\noindent {\it Proof:}
By construction $yv_i=q^{d-2i}v_i$ for $0 \leq i \leq d$.
Also $xv_0=q^{-d}v_0$
by Lemma
\ref{lem:meaning}(iii) and
$zv_d=q^dv_d$
by Lemma
\ref{lem:meaning}(i).
Abbreviate $\eta = \sum_{i=0}^d v_i$.
\\
\noindent
${\rm (i) \Leftrightarrow (ii)}$ 
By Lemma
\ref{lem:meaning}(ii),
$\eta\in n^d_yV$ if and only if $z\eta =q^{-d}\eta $. 
Using $zv_d=q^dv_d$ we obtain
$(z-q^{-d})\eta = \sum_{i=0}^{d-1} w_i$, where
\begin{eqnarray*}
w_i = (z-q^{2i-d})v_{i} + (q^{2i+2-d}-q^{-d})v_{i+1}
\qquad (0 \leq i \leq d-1).
\end{eqnarray*}
By Theorem
\ref{thm:xyzV}, 
for $0 \leq i \leq d-1$
the vector $w_i$ is contained in component
$i+1$ of
$\lbrack y \rbrack$.
Thus 
$(z-q^{-d})\eta =0$ if and only if
$w_i=0$ for $0 \leq i \leq d-1$. The result follows.
\\
${\rm (ii) \Leftrightarrow (iii)}$ 
Using the equation on the right in
(\ref{eq:xz}),
\begin{eqnarray*}
(z-q^{2i-d})v_i = -(q-q^{-1})q^{2i-d+1}n_x v_i
\qquad (0 \leq i \leq d-1).
\end{eqnarray*} 
The result follows.
\\
\noindent
${\rm (i) \Leftrightarrow (iv)}$ 
By Lemma
\ref{lem:meaning}(ii),
$\eta \in n^d_yV$ if and only if $x\eta =q^{d}\eta $. 
Using $xv_0=q^{-d}v_0$ we obtain
$(x-q^{d})\eta = \sum_{i=1}^{d} u_i$ where
\begin{eqnarray*}
u_i = (x-q^{2i-d})v_{i} + (q^{2i-2-d}-q^{d})v_{i-1}
\qquad (1 \leq i \leq d).
\end{eqnarray*}
By Theorem
\ref{thm:xyzV}, 
for $1 \leq i \leq d$
the vector $u_i$ is contained in component
$i-1$ of
$\lbrack y \rbrack$.
Thus 
$(x-q^{d})\eta =0$ if and only if
$u_i=0$ for $1 \leq i \leq d$. The result follows.
\\
${\rm (iv) \Leftrightarrow (v)}$ 
Using the equation on the left in
(\ref{eq:xz}),
\begin{eqnarray*}
(x-q^{2i-d})v_i = -(q-q^{-1})q^{2i-d-1}n_z v_i
\qquad (1 \leq i \leq d).
\end{eqnarray*} 
The result follows.
\\
\noindent Now assume that (i)--(v) hold. 
By condition
(iii), $v_i=0$ implies $v_{i+1}=0$ for
$0 \leq i \leq d-1$. 
By condition
(v), $v_i=0$ implies $v_{i-1}=0$ for
$1 \leq i \leq d$. Therefore
$\lbrace v_i\rbrace_{i=0}^d$ are all zero or all nonzero.
\hfill $\Box$ \\

\begin{lemma} 
\label{lem:pp}
Let $\lbrace v_i\rbrace_{i=0}^d$ 
denote vectors in $V$, not all zero. Then
the following are equivalent:
\begin{enumerate}
\item[\rm (i)]
$\lbrace v_i\rbrace_{i=0}^d$ is a 
$\lbrack y \rbrack_{row} $ basis for $V$;
\item[\rm (ii)]
$yv_0=q^dv_0$ and
$(z-q^{2i-d}) v_i = (q^{-d}-q^{2i+2-d})
v_{i+1}$ for $0 \leq i \leq d-1$;
\item[\rm (iii)]
$yv_0=q^dv_0$ and
 $n_x v_i = q^{-i}\lbrack i+1\rbrack 
v_{i+1}$ for $0 \leq i \leq d-1$;
\item[\rm (iv)]
$yv_d=q^{-d}v_d$ and
$(x-q^{2i-d}) v_i = (q^d-q^{2i-2-d})v_{i-1}$
 for $1 \leq i \leq d$;
\item[\rm (v)]
$y v_d=q^{-d}v_d$ and
$n_z v_i = -q^{d-i}\lbrack d-i+1\rbrack
v_{i-1}$ for $1 \leq i \leq d$.
\end{enumerate}
Now assume that {\rm (i)--(v)} hold. Then
\begin{eqnarray}
zv_d = q^dv_d, \qquad
n_xv_d = 0, \qquad
xv_0 = q^{-d}v_0, \qquad
n_zv_0 = 0.
\label{eq:extra}
\end{eqnarray}
\end{lemma}
\noindent {\it Proof:}
Each condition (i)--(v) implies that
that for $0 \leq i \leq d$ the vector
$v_i$ is contained in component $i$ of
$\lbrack  y \rbrack$. 
Now these conditions are equivalent
in view of Lemma
\ref{lem:rowbasispre}.
Next assume that (i)--(v) hold. Then
the equations (\ref{eq:extra}) hold by
Lemma
\ref{lem:meaning}
and Lemma \ref{lem:ndvasker}.
\hfill $\Box$ \\

\begin{lemma}
\label{lem:urow}
Pick $\xi \in \lbrace x,y,z\rbrace $.
There exists a $\lbrack \xi \rbrack_{row}$ basis for $V$
and $V^*$. 
\end{lemma}
\noindent {\it Proof:}
Without loss we may assume that
the underlying vector space is $V$.
First suppose that $\xi=y$. Let $v_0$ denote
a nonzero vector in component $0$ of
the decomposition $\lbrack y \rbrack$ of $V$.
Thus $yv_0=q^dv_0$.
For $0 \leq i \leq d-1$ define
$v_{i+1}$ to satisfy Lemma
\ref{lem:pp}(iii).
By construction 
the sequence
$\lbrace v_i\rbrace_{i=0}^d$ satisfies
 Lemma
\ref{lem:pp}(iii).
By that lemma
$\lbrace v_i\rbrace_{i=0}^d$ is
a $\lbrack y \rbrack_{row} $ basis for $V$.
We have proven the result for $\xi=y$.
To get the result for the remaining values of
$\xi$ use
 Lemma
\ref{lem:om}.
\hfill $\Box$ \\

\noindent In Definition
\ref{def:urow} we defined some bases for $V$ and $V^*$.
These bases are not unique; we will discuss this issue
in Lemma
\ref{lem:notunique}.

\begin{lemma} 
\label{lem:inducerow}
Pick $\xi \in \lbrace x,y,z\rbrace$.
Then for $V$ and $V^*$,
the decomposition $\lbrack \xi \rbrack$
is induced by each  
$\lbrack \xi \rbrack_{row}$ basis.
Moreover the
 decomposition $\lbrack \xi \rbrack^{inv}$
is induced by each 
$\lbrack \xi \rbrack^{inv}_{row}$ basis.
\end{lemma}
\noindent {\it Proof:}
The first assertion follows from Definition
\ref{def:urow}(i).
The second assertion follows by the meaning of inversion.
\hfill $\Box$ \\

%

\noindent
Let $\lbrace u_i\rbrace_{i=0}^d$ denote a basis for
$V$ and let
 $\lbrace v_i\rbrace_{i=0}^d$
denote a basis for $V^*$.
These bases are said to be {\it dual}
 whenever
$(u_i, v_j)=\delta_{ij}$ for $0 \leq i,j\leq d$.
Each basis for $V$ (resp. $V^*$) is dual to a unique
basis for $V^*$ (resp. $V$).

\begin{definition}
\label{def:ucol}
\rm
Pick $\xi \in \lbrace x,y,z\rbrace$.
A basis for $V$ (resp. $V^*$) is called
$\lbrack \xi \rbrack_{col}$ whenever it is dual to
a $\lbrack \xi \rbrack^{inv}_{row}$ basis for $V^*$ (resp. $V$).
By a 
$\lbrack \xi \rbrack^{inv}_{col}$ basis we mean
the inversion of a 
$\lbrack \xi \rbrack_{col}$ basis.
\end{definition}

\begin{lemma} 
\label{lem:inducecol}
Pick $\xi \in \lbrace x,y,z\rbrace$.
Then for $V$ and $V^*$,
the decomposition $\lbrack \xi \rbrack$
is induced by each 
$\lbrack \xi \rbrack_{col}$ basis.
Moreover the
 decomposition $\lbrack \xi \rbrack^{inv}$
is induced by each 
$\lbrack \xi \rbrack^{inv}_{col}$ basis.
\end{lemma}
\noindent {\it Proof:}
Use
Lemma
\ref{lem:dualdec},
Lemma \ref{lem:inducerow},
and 
Definition
\ref{def:ucol}.
\hfill $\Box$ \\


\noindent In Definition
\ref{def:urow}
and
Definition
\ref{def:ucol}
we defined the bases 
(\ref{eq:basisr1})--(\ref{eq:basisr4})
for $V$ and $V^*$. We now discuss the uniqueness
of these bases. For notational convenience
we will focus on the 
$\lbrack y \rbrack_{row}$  basis for $V$; similar
results apply to the remaining bases.

\begin{lemma}
\label{lem:notunique}
Let
$\lbrace v_i\rbrace_{i=0}^d$ denote a
$\lbrack y \rbrack_{row}$ basis for $V$.
Let 
$\lbrace v'_i\rbrace_{i=0}^d$ denote any vectors
in $V$. Then the following are equivalent:
\begin{enumerate}
\item[\rm (i)] the sequence 
$\lbrace v'_i\rbrace_{i=0}^d$ is a 
$\lbrack y \rbrack_{row}$ basis for $V$;
\item[\rm (ii)] there exists 
$0 \not=\alpha  \in \F$ such that $v'_i=\alpha v_i$
for $0 \leq i \leq d$.
\end{enumerate}
\end{lemma}
\noindent {\it Proof:}
Use Lemma
\ref{lem:pp}.
\hfill $\Box$ \\


%

\begin{lemma} 
\label{lem:whatdualwhat}
Pick $\xi \in \lbrace x,y,z\rbrace$.
For the table below, in each row
 we display a basis for $V$ and its dual basis for $V^*$.
\medskip

\centerline{
\begin{tabular}[t]{c c}
   {\rm basis for $V$} & 
   {\rm dual basis for $V^*$}
   \\ \hline  \hline
  $\lbrack \xi \rbrack_{row}$   & $\lbrack \xi \rbrack^{inv}_{col} $
  \\ 
  $\lbrack \xi \rbrack_{col}$   & $\lbrack \xi \rbrack^{inv}_{row}$ 
   \\ 
  $\lbrack \xi \rbrack^{inv}_{row}$   & $\lbrack \xi \rbrack_{col} $
  \\ 
  $\lbrack \xi \rbrack^{inv}_{col}$   & $\lbrack \xi \rbrack_{row}$ 
     \end{tabular}}
     \medskip

\end{lemma}
\noindent {\it Proof:}
By Definition
\ref{def:ucol}
and the meaning of inversion.
\hfill $\Box$ \\


\section{The matrices representing $x,y,z$ with respect to the twelve bases }

\noindent
We continue to discuss the
$U_q(\mathfrak{sl}_2)$-module $V$ and
the 
$U_{q^{-1}}(\mathfrak{sl}_2)$-module $V^*$.
Recall the twelve bases
(\ref{eq:basisr1})--(\ref{eq:basisr4})
for $V$ and $V^*$. In this section we find the
matrices that represent $x,y,z$ with respect to these bases.

\medskip
\noindent We will use the following notation.
Let ${\rm Mat}_{d+1}(\F)$ denote the $\F$-algebra consisting of
the $d+1$ by $d+1$ matrices that have all entries in $\F$.
We index the rows and columns by $0,1,\ldots, d$.
Let $\lbrace v_i\rbrace_{i=0}^d$
denote a basis for $V$. For $A \in {\rm End}(V)$
and $B\in 
{\rm Mat}_{d+1}(\F)$, we say that {\it $B$ represents $A$
with respect to 
 $\lbrace v_i\rbrace_{i=0}^d$} whenever
$Av_j = \sum_{i=0}^d B_{ij}v_i$ for $0 \leq j \leq d$.

\medskip
\noindent We have a comment.
Let $\lbrace u_i\rbrace_{i=0}^d$ denote a basis for
$V$ and let
 $\lbrace v_i\rbrace_{i=0}^d$
denote the basis for $V^*$ that is dual to 
$\lbrace u_i\rbrace_{i=0}^d$.
 Pick $A\in 
 {\rm End}(V)$ and let $B$
denote the matrix in
${\rm Mat}_{d+1}(\F)$ that represents $A$ with
respect to 
 $\lbrace u_i\rbrace_{i=0}^d$.
Then the transpose $B^t$ represents 
the adjoint $A^{adj}$ with respect to 
 $\lbrace v_i\rbrace_{i=0}^d$.

\begin{lemma}
\label{lem:btrans}
Let
 $\lbrace u_i\rbrace_{i=0}^d$ denote a basis for
$V$ and let
 $\lbrace v_i\rbrace_{i=0}^d$
denote the basis for $V^*$ that is dual to 
$\lbrace u_i\rbrace_{i=0}^d$.
Pick
 $\zeta \in U_q(\mathfrak{sl}_2)$ 
and let
$B$ denote the matrix 
in ${\rm Mat}_{d+1}(\F)$ that represents $\zeta$ with
respect to 
 $\lbrace u_i\rbrace_{i=0}^d$.
Then $B^t$ represents  $\zeta^\dagger$
with respect to 
 $\lbrace v_i\rbrace_{i=0}^d$.
\end{lemma} 
\noindent {\it Proof:} 
By Proposition 
\ref{prop:uquqi}
and the comment above this lemma.
\hfill $\Box$ \\

\medskip
\noindent  We now define some matrices in
${\rm Mat}_{d+1}(\F)$.

\begin{definition} \rm
Let $K_q$ denote the 
diagonal matrix
in ${\rm Mat}_{d+1}(\F)$ with $(i,i)$-entry
$q^{d-2i}$ for $0 \leq i \leq d$.
\end{definition}

\begin{example} \rm
For $d=3$,
\begin{eqnarray*}
K_q = \mbox{\rm diag}( q^3, q, q^{-1}, q^{-3}).
\end{eqnarray*}
\end{example}

\begin{definition}
\label{def:z}
\rm
We define a matrix $Z \in 
{\rm Mat}_{d+1}(\F)$ as follows.
For $0 \leq i,j\leq d$ the $(i,j)$-entry is
$\delta_{i+j,d}$.
Note that $Z^2=I$.
\end{definition}

\begin{example} \rm
For  $d=3$,
\begin{eqnarray*}
&&
Z = \quad 
\left(
\begin{array}{ cccc}
0&0&0&1  \\
  0 &0&1&0\\
   0& 1 &0&0\\
 1&0  &0   &0
\end{array}
\right).
\end{eqnarray*}
\end{example}

\begin{lemma}
\label{lem:ZBZ}
For $B \in 
{\rm Mat}_{d+1}(\F)$  and 
$0 \leq i,j\leq d$ the following coincide:
\begin{enumerate}
\item[\rm (i)] the $(i,j)$-entry of $ZBZ$;
\item[\rm (i)] the $(d-i,d-j)$-entry of $B$.
\end{enumerate}
\end{lemma}
\noindent {\it Proof:} Use matrix multiplication.
\hfill $\Box$ \\

\noindent Let $B$ denote a matrix in 
${\rm Mat}_{d+1}(\F)$. Then $B$ is called
{\it lower bidiagonal} whenever both
(i) each nonzero entry is on the diagonal or the subdiagonal;
(ii) each entry on the subdiagonal is  nonzero.
The matrix $B$ is called 
{\it upper bidiagonal} whenever $B^t$ is lower
bidiagonal.

\begin{definition}
\label{def:eq}
\rm
Let 
$E_q$  denote the upper bidiagonal matrix
in 
${\rm Mat}_{d+1}(\F)$ 
with $(i,i)$-entry
$q^{2i-d}$ for $0 \leq i \leq d$ and
$(i-1,i)$-entry $q^d-q^{2i-2-d}$
for $1 \leq i \leq d$.
\end{definition}

\noindent We will be discussing the following eight matrices:
\begin{eqnarray}
&&E_q, \qquad \quad \;\;
E_{q^{-1}},  \qquad \quad\;\; 
E^t_q, \qquad  \quad\;\;
E^t_{q^{-1}}, 
\label{eq:ematrow1}
\\
&&ZE_qZ, \qquad 
ZE_{q^{-1}}Z,  \qquad
ZE^t_qZ, \qquad 
ZE^t_{q^{-1}}Z.
\label{eq:ematrow2}
\end{eqnarray}

\begin{lemma}
For the matrices
{\rm (\ref{eq:ematrow1})},
{\rm (\ref{eq:ematrow2})} we display the entries in the
table below. Each entry not shown is zero.

\medskip

\centerline{
\begin{tabular}[t]{c |ccc}
  {\rm matrix} & 
  {\rm $(i,i-1)$-entry} &
  {\rm $(i,i)$-entry} &
  {\rm $(i-1,i)$-entry} 
   \\ \hline  \hline
  $E_q$ & 
  $0$ 
  &
  $q^{2i-d}$
  &
  $q^d-q^{2i-2-d}$
  \\
  $E_{q^{-1}}$ & 
  $0$ 
  &
  $q^{d-2i}$
  &
  $q^{-d}-q^{d-2i+2}$
  \\
  $E^t_q$ & 
  $q^d-q^{2i-2-d}$
  &
  $q^{2i-d}$
  &
  $0$ 
  \\
  $E^t_{q^{-1}}$ & 
  $q^{-d}-q^{d-2i+2}$
  &
  $q^{d-2i}$
  &
  $0$ 
  \\
 \hline
  $ZE_qZ$ & 
  $q^d-q^{d-2i}$
  &
  $q^{d-2i}$
  &
  $0$ 
  \\
  $ZE_{q^{-1}}Z$ & 
  $q^{-d}-q^{2i-d}$
  &
  $q^{2i-d}$
  &
  $0$ 
  \\
  $ZE^t_qZ$ & 
  $0$ 
  &
  $q^{d-2i}$
  &
  $q^d-q^{d-2i}$
  \\
  $ZE^t_{q^{-1}}Z$ & 
  $0$ 
  &
  $q^{2i-d}$
  &
  $q^{-d}-q^{2i-d}$
     \end{tabular}}
     \medskip

\end{lemma}
\noindent {\it Proof:} 
Use Lemma \ref{lem:ZBZ}.
\hfill $\Box$ \\

\noindent Let $B$ denote a matrix in 
               ${\rm Mat}_{d+1}(\F)$.
For $\alpha \in \F$, $B$ is said to have {\it constant row sum $\alpha$}
whenever $\alpha = \sum_{j=0}^d B_{ij} $ for $0 \leq i \leq d$.
The matrix $B$ is said to have {\it constant column sum $\alpha$}
wheneve $B^t$ has constant row sum $\alpha$.

\begin{lemma} 
Each of the matrices 
{\rm (\ref{eq:ematrow1})},
{\rm (\ref{eq:ematrow2})} is described as follows:
{\rm (i)} it is upper or lower bidiagonal;
{\rm (ii)} the diagonal part is $K_q$ or $K^{-1}_q$;
{\rm (iii)} it has constant row sum or constant column sum. 
The details
are given in the table below.

\medskip

\centerline{
\begin{tabular}[t]{c |ccc}
  {\rm matrix} & 
  {\rm upper/lower bidiag.} &
  {\rm diagonal part} &
  {\rm row/colum sum} 
   \\ \hline  \hline
  $E_q$ & 
  {\rm upper bidiag.} & {\rm $K^{-1}_q$} & 
  {\rm const. row sum $q^d$}
  \\
  $E_{q^{-1}}$ & 
  {\rm upper bidiag.} & {\rm  $K_q$} & 
  {\rm const. row sum $q^{-d}$}
  \\
  $E^t_q$ & 
  {\rm lower bidiag.} & {\rm  $K^{-1}_q$} & 
  {\rm const. column sum $q^{d}$}
  \\
  $E^t_{q^{-1}}$ & 
  {\rm lower bidiag.} & {\rm  $K_q$} & 
  {\rm const. column sum $q^{-d}$}
  \\
 \hline
  $ZE_qZ$ & 
  {\rm lower bidiag.} & {\rm  $K_q$} & 
  {\rm const. row sum $q^{d}$}
  \\
  $ZE_{q^{-1}}Z$ & 
  {\rm lower bidiag.} & {\rm  $K^{-1}_q$} & 
  {\rm const. row sum $q^{-d}$}
  \\
  $ZE^t_qZ$ & 
  {\rm upper bidiag.} & {\rm  $K_q$} & 
  {\rm const. column sum $q^{d}$}
  \\
  $ZE^t_{q^{-1}}Z$ & 
  {\rm upper bidiag.} & {\rm $K^{-1}_q$} & 
  {\rm const. column sum $q^{-d}$}
     \end{tabular}}
     \medskip

\end{lemma}

\begin{example}\rm For $d=3$,

\begin{eqnarray*}
&&
E_q=
\left(
\begin{array}{ cccc}
q^{-3} & q^3-q^{-3}  &0&0  \\
0 &   q^{-1} & q^3-q^{-1} &0\\
0 &0&  q &  q^3-q \\
 0&0&0&  q^3 
\end{array}
\right),
\\
&&
\qquad \qquad \qquad 
 \qquad 
E_{q^{-1}}=
\left(
\begin{array}{ cccc}
q^{3} & q^{-3}-q^{3}  &0& 0 \\
 0&   q & q^{-3}-q &0\\
 0&0&  q^{-1} &  q^{-3}-q^{-1} \\
 0&0&0&  q^{-3} 
\end{array}
\right),
\end{eqnarray*}
\begin{eqnarray*}
&&
E^t_q = 
\left(
\begin{array}{ ccc c}
q^{-3} &  0 &0& 0 \\
 q^{3}-q^{-3} &  q^{-1} &0&0\\
 0&q^{3} - q^{-1}&  q &0 \\
 0&0&q^{3}-q&  q^3 
\end{array}
\right),
\\
&&
\qquad \qquad \qquad 
\qquad 
E^t_{q^{-1}}=
\left(
\begin{array}{ ccc c}
q^{3} & 0  &0&0 \\
 q^{-3}-q^{3} &  q &0&0\\
0 &q^{-3} - q&  q^{-1} & 0\\
 0&0&q^{-3}-q^{-1}&  q^{-3} 
\end{array}
\right),
\end{eqnarray*}
\begin{eqnarray*}
&&Z E_q Z=
\left(
\begin{array}{ ccc c}
q^{3} &  0 &0& 0 \\
 q^{3}-q &  q &0&0\\
 0&q^{3} - q^{-1}&  q^{-1} &0 \\
 0&0&q^{3}-q^{-3}&  q^{-3} 
\end{array}
\right),
\\
&&
\qquad \qquad \qquad 
 \qquad 
ZE_{q^{-1}}Z=
\left(
\begin{array}{ ccc c}
q^{-3} &  0 &0&0  \\
 q^{-3}-q^{-1} &  q^{-1} &0&0\\
0 &q^{-3} - q&  q & 0\\
0 &0&q^{-3}-q^3&  q^3 
\end{array}
\right),
\end{eqnarray*}
\begin{eqnarray*}
&&Z E^t_q Z =
\left(
\begin{array}{ ccc c}
q^{3} & q^{3}-q  &0&0  \\
0 &   q & q^{3}-q^{-1} &0\\
0 &0&  q^{-1} &  q^{3}-q^{-3} \\
0 &0&0&  q^{-3}  
\end{array}
\right),
\\
&&
\qquad \qquad \qquad 
 \qquad 
ZE^t_{q^{-1}} Z=
\left(
\begin{array}{ ccc c}
q^{-3} & q^{-3}-q^{-1}  &0&0  \\
0 &   q^{-1} & q^{-3}-q &0\\
 0&0&  q &  q^{-3}-q^3 \\
 0&0&0&  q^3   
\end{array}
\right).
\end{eqnarray*}

\end{example}

\begin{note}
\label{zzzaction}
\rm Consider the set of eight matrices
{\rm (\ref{eq:ematrow1})},
{\rm (\ref{eq:ematrow2})}.
The set is closed under each of the following
maps:
\begin{enumerate}
\item[\rm (i)] the transpose map;
\item[\rm (ii)] replace $q$ by $q^{-1}$;
\item[\rm (iii)] conjugation by $Z$.
\end{enumerate}
Each of the maps (i)--(iii) has order 2, and
these maps 
mutually commute.  This gives an action of the
group $\Z_2 \times \Z_2 \times \Z_2$ on
the set of eight matrices
{\rm (\ref{eq:ematrow1})},
{\rm (\ref{eq:ematrow2})}.  This action is
transitive.
\end{note}

\begin{theorem}
\label{thm:vmat}
Consider the elements 
$x$, $y$, $z$ of
 $U_q(\mathfrak{sl}_2)$.
In the table below we display the matrices that represent
these elements with respect to the twelve bases for $V$ from 
{\rm (\ref{eq:basisr1})--(\ref{eq:basisr4})}.

\medskip

\centerline{
\begin{tabular}[t]{c |ccc}
  {\rm basis} & 
  {\rm $x$} &
  {\rm  $y$} &
  {\rm  $z$} 
   \\ \hline  \hline
 $\lbrack x \rbrack_{row}$
 &
 $K_q$
 &
 $ZE_{q^{-1}}Z$
 &
 $E_q$
 \\
 $\lbrack x \rbrack_{col}$
 &
 $K_q$
 &
 $E^t_q$
 &
 $ZE^t_{q^{-1}}Z$
 \\
 $\lbrack x \rbrack^{inv}_{row}$
 &
 $K^{-1}_q$
 &
 $E_{q^{-1}}$
 & 
 $ZE_qZ$
 \\
 $\lbrack x \rbrack^{inv}_{col}$
 &
 $K^{-1}_q$
 &
 $ZE^t_qZ$
  & 
 $E^t_{q^{-1}}$
 \\
 \hline
 $\lbrack y \rbrack_{row}$
 & 
 $E_q$
 &
 $K_q$
 &
 $ZE_{q^{-1}}Z$
 \\
 $\lbrack y \rbrack_{col}$
 &
 $ZE^t_{q^{-1}}Z$
& 
 $K_q$
 &
 $E^t_q$
 \\
 $\lbrack y \rbrack^{inv}_{row}$
 & 
 $ZE_qZ$
 &
 $K^{-1}_q$
 &
 $E_{q^{-1}}$
 \\
 $\lbrack y \rbrack^{inv}_{col}$
  & 
 $E^t_{q^{-1}}$
 &
 $K^{-1}_q$
 &
 $ZE^t_qZ$
 \\
 \hline
 $\lbrack z \rbrack_{row}$
 &
 $ZE_{q^{-1}}Z$
 & 
 $E_q$
 &
 $K_q$
 \\
 $\lbrack z \rbrack_{col}$
 &
 $E^t_q$
 &
 $ZE^t_{q^{-1}}Z$
& 
 $K_q$
 \\
 $\lbrack z \rbrack^{inv}_{row}$
 &
 $E_{q^{-1}}$
 & 
 $ZE_qZ$
 &
 $K^{-1}_q$
 \\
 $\lbrack z \rbrack^{inv}_{col}$
 &
 $ZE^t_qZ$
  & 
 $E^t_{q^{-1}}$
 &
 $K^{-1}_q$
     \end{tabular}}
     \medskip

\end{theorem}
\noindent {\it Proof:} 
We first verify the data for the middle third of the table.
Using Lemma
\ref{lem:pp} and the construction,
we get the matrices that represent
 $x,y,z$ with respect to a 
 $\lbrack y \rbrack_{row}$ basis for $V$.
For these matrices conjugate by $Z$ to get the
matrices that represent $x,y,z$ with respect to
 a $\lbrack y \rbrack^{inv}_{row}$ basis for $V$.
For these matrices replace $q$ by $q^{-1}$ to get
the matrices that represent $x,y,z$ with respect to
a 
  $\lbrack y \rbrack^{inv}_{row}$ basis for $V^*$.
For these matrices take the transpose and invoke
Lemma 
\ref{lem:btrans}
to get the matrices that represent 
 $x,y,z$ with respect to
a 
  $\lbrack y \rbrack_{col}$ basis for $V$.
For these matrices conjugate by $Z$ to get the
 matrices that represent 
 $x,y,z$ with respect to
a 
  $\lbrack y \rbrack^{inv}_{col}$ basis for $V$.
We have now verified the data for the middle third of
the table. To verify the rest of the table use Lemma
\ref{lem:om}.
\hfill $\Box$ \\

\begin{theorem}
\label{thm:vsmat}
Consider the elements 
$x$, $y$, $z$ of
 $U_{q^{-1}}(\mathfrak{sl}_2)$.
In the table below we display the matrices that represent
these elements with respect to the twelve bases for $V^*$ from 
{\rm (\ref{eq:basisr1})--(\ref{eq:basisr4})}.

\medskip

\centerline{
\begin{tabular}[t]{c |ccc}
  {\rm basis} & 
  {\rm $x$} &
  {\rm  $y$} &
  {\rm  $z$} 
   \\ \hline  \hline
 $\lbrack x \rbrack_{row}$
 &
 $K^{-1}_{q}$
 &
 $ZE_{q}Z$
 &
 $E_{q^{-1}}$
 \\
 $\lbrack x \rbrack_{col}$
 &
 $K^{-1}_{q}$
 &
 $E^t_{q^{-1}}$
 &
 $ZE^t_{q}Z$
 \\
 $\lbrack x \rbrack^{inv}_{row}$
 &
 $K_q$
 &
 $E_{q}$
 & 
 $ZE_{q^{-1}}Z$
 \\
 $\lbrack x \rbrack^{inv}_{col}$
 &
 $K_q$
 &
 $ZE^t_{q^{-1}}Z$
  & 
 $E^t_{q}$
 \\
 \hline
 $\lbrack y \rbrack_{row}$
 & 
 $E_{q^{-1}}$
 &
 $K^{-1}_{q}$
 &
 $ZE_{q}Z$
 \\
 $\lbrack y \rbrack_{col}$
 &
 $ZE^t_{q}Z$
& 
 $K^{-1}_{q}$
 &
 $E^t_{q^{-1}}$
 \\
 $\lbrack y \rbrack^{inv}_{row}$
 & 
 $ZE_{q^{-1}}Z$
 &
 $K_q$
 &
 $E_{q}$
 \\
 $\lbrack y \rbrack^{inv}_{col}$
  & 
 $E^t_{q}$
 &
 $K_q$
 &
 $ZE^t_{q^{-1}}Z$
 \\
 \hline
 $\lbrack z \rbrack_{row}$
 &
 $ZE_{q}Z$
 & 
 $E_{q^{-1}}$
 &
 $K^{-1}_q$
 \\
 $\lbrack z \rbrack_{col}$
 &
 $E^t_{q^{-1}}$
 &
 $ZE^t_{q}Z$
& 
 $K^{-1}_q$
 \\
 $\lbrack z \rbrack^{inv}_{row}$
 &
 $E_{q}$
 & 
 $ZE_{q^{-1}}Z$
 &
 $K_q$
 \\
 $\lbrack z \rbrack^{inv}_{col}$
 &
 $ZE^t_{q^{-1}}Z$
  & 
 $E^t_{q}$
 &
 $K_q$
     \end{tabular}}
     \medskip

\end{theorem}
\noindent {\it Proof:}
In the table of Theorem
\ref{thm:vmat} replace $q$ by $q^{-1}$.
\hfill $\Box$ \\

\section{The matrices representing $n_x, n_y, n_z$ with respect to the
twelve bases}

\noindent We continue to discuss the
$U_q(\mathfrak{sl}_2)$-module $V$ and
the 
$U_{q^{-1}}(\mathfrak{sl}_2)$-module $V^*$.
Recall the twelve bases
(\ref{eq:basisr1})--(\ref{eq:basisr4}) for $V$ and $V^*$.
In the previous section we found the matrices that represent
$x,y,z$ with respect to these bases. In the present section we
find the matrices that represent
$n_x,n_y,n_z$ with respect to these bases.

\begin{definition} 
\label{def:Nq}
\rm
Let $N_q$ denote the matrix in 
${\rm Mat}_{d+1}(\F)$ 
with $(i,i-1)$-entry
$q^{1-i}\lbrack i \rbrack$ for
$1 \leq i \leq d$, and all other
entries 0.
\end{definition}

\noindent Recall the matrix 
 $Z$
from
Definition \ref{def:z}.
We will be discussing the following
eight matrices:
\begin{eqnarray}
&&N_q, \qquad \quad \;\;
N_{q^{-1}},  \qquad \quad\;\; 
N^t_q, \qquad  \quad\;\;
N^t_{q^{-1}}, 
\label{eq:nmatrow1}
\\
&&ZN_qZ, \qquad 
ZN_{q^{-1}}Z,  \qquad
ZN^t_qZ, \qquad 
ZN^t_{q^{-1}}Z.
\label{eq:nmatrow2}
\end{eqnarray}

\begin{lemma}
\label{lem:Nqentries}
For the matrices
{\rm (\ref{eq:nmatrow1})},
{\rm (\ref{eq:nmatrow2})} we display the entries in the
table below. Each entry not shown is zero.

\medskip

\centerline{
\begin{tabular}[t]{c |cc}
  {\rm matrix} & 
  {\rm $(i,i-1)$-entry} &
  {\rm $(i-1,i)$-entry} 
   \\ \hline  \hline
  $N_q$ & 
  $q^{1-i}\lbrack i \rbrack $ 
  &
  $0$
  \\
  $N_{q^{-1}}$ & 
  $q^{i-1}\lbrack i \rbrack $ 
  &
  $0$
  \\
  $N^t_{q}$ & 
  $0$ & 
  $q^{1-i}\lbrack i \rbrack$
  \\
  $N^t_{q^{-1}}$ & 
  $0$ & 
  $q^{i-1}\lbrack i \rbrack $ 
  \\
 \hline
  $ZN_qZ$ & 
  $0$ 
   & 
  $q^{i-d}\lbrack d-i+1 \rbrack $ 
  \\
  $ZN_{q^{-1}}Z$ & 
  $0$ 
   & 
  $q^{d-i}\lbrack d-i+1 \rbrack $ 
  \\
  $ZN^t_qZ$ & 
  $q^{i-d}\lbrack d-i+1 \rbrack $ 
 & 
  $0$ 
  \\
  $ZN^t_{q^{-1}}Z$ & 
  $q^{d-i}\lbrack d-i+1 \rbrack $ 
 & 
  $0$ 
     \end{tabular}}
     \medskip

\end{lemma}
\noindent {\it Proof:} 
Use Lemma \ref{lem:ZBZ}.
\hfill $\Box$ \\


\begin{example}\rm For $d=3$,
\begin{eqnarray*}
&&
N_q=
\left(
\begin{array}{ cccc}
0 & 0  &0&0  \\
\lbrack 1 \rbrack &  0& 0 &0\\
0 &q^{-1}\lbrack 2 \rbrack &  0 &  0 \\
 0&0&q^{-2}\lbrack 3 \rbrack&  0 
\end{array}
\right),
  \qquad 
N_{q^{-1}}=
\left(
\begin{array}{ cccc}
0 & 0  &0&0  \\
\lbrack 1 \rbrack &  0& 0 &0\\
0 &q\lbrack 2 \rbrack &  0 &  0 \\
 0&0&q^{2}\lbrack 3 \rbrack&  0 
\end{array}
\right),
\end{eqnarray*}
\begin{eqnarray*}
&&
N^t_q = 
\left(
\begin{array}{ cccc}
0 & \lbrack 1 \rbrack  &0&0  \\
0 &  0& q^{-1}\lbrack 2 \rbrack &0\\
0 &0 &  0 &  q^{-2}\lbrack 3 \rbrack \\
 0&0& 0&  0 
\end{array}
\right),
 \qquad 
N^t_{q^{-1}}=
\left(
\begin{array}{ cccc}
0 & \lbrack 1 \rbrack  &0&0  \\
0 &  0& q\lbrack 2 \rbrack &0\\
0 &0 &  0 &  q^{2}\lbrack 3 \rbrack \\
 0&0& 0&  0 
\end{array}
\right),
\end{eqnarray*}
\begin{eqnarray*}
&&Z N_q Z=
\left(
\begin{array}{ cccc}
0 & q^{-2}\lbrack 3 \rbrack  &0&0  \\
0 &  0& q^{-1}\lbrack 2 \rbrack &0\\
0 &0 &  0 &  \lbrack 1 \rbrack \\
 0&0& 0&  0 
\end{array}
\right),
 \qquad 
ZN_{q^{-1}}Z=
\left(
\begin{array}{ cccc}
0 & q^{2}\lbrack 3 \rbrack  &0&0  \\
0 &  0& q\lbrack 2 \rbrack &0\\
0 &0 &  0 &  \lbrack 1 \rbrack \\
 0&0& 0&  0 
\end{array}
\right),
\end{eqnarray*}
\begin{eqnarray*}
&&Z N^t_q Z =
\left(
\begin{array}{ cccc}
0 & 0  &0&0  \\
q^{-2}\lbrack 3 \rbrack &  0& 0 &0\\
0 &q^{-1}\lbrack 2 \rbrack &  0 &  0 \\
 0&0&\lbrack 1 \rbrack&  0 
\end{array}
\right),
 \qquad 
ZN^t_{q^{-1}} Z=
\left(
\begin{array}{ cccc}
0 & 0  &0&0  \\
q^{2}\lbrack 3 \rbrack &  0& 0 &0\\
0 &q\lbrack 2 \rbrack &  0 &  0 \\
 0&0&\lbrack 1 \rbrack&  0 
\end{array}
\right).
\end{eqnarray*}

\end{example}

\begin{note}\rm
Consider the set of eight matrices
{\rm (\ref{eq:nmatrow1})},
{\rm (\ref{eq:nmatrow2})}.
This set is closed under each of the 
three maps from Note
\ref{zzzaction}. This gives an action of
the group $\Z_2 \times \Z_2 \times \Z_2$ on
the set of eight matrices
{\rm (\ref{eq:nmatrow1})},
{\rm (\ref{eq:nmatrow2})}.  This action is
transitive.
\end{note}

\begin{definition} 
\label{def:tq}
\rm
Let $T_q$ denote the tridiagonal matrix in 
${\rm Mat}_{d+1}(\F)$ with the following entries.
For $1 \leq i \leq d$ the $(i,i-1)$-entry
is $q^{3i-2d-1}\lbrack i \rbrack $ and
the $(i-1,i)$-entry
is $-q^{3i-d-2}\lbrack d-i+1 \rbrack $.
For $0 \leq i \leq d$
the $(i,i)$-entry is 
\begin{eqnarray*}
q^{2i-d}\lbrack i \rbrack \lbrack d-i+1\rbrack (q-q^{-1})-
q^{2i-d+1} \lbrack 2i-d\rbrack.
\end{eqnarray*}
\end{definition}
\begin{note} 
\label{note:ztz}
\rm We have $Z T_q Z = -T_{q^{-1}}$. This is
routinely checked using
Lemma
\ref{lem:ZBZ}.
\end{note}

\begin{theorem}
\label{thm:nxmat}
Consider the elements
$n_x$, $n_y$, $n_z$ of
 $U_{q}(\mathfrak{sl}_2)$.
In the table below we display the matrices that represent
these elements with respect to the twelve bases
for $V$ from 
{\rm (\ref{eq:basisr1})--(\ref{eq:basisr4})}.

\medskip

\centerline{
\begin{tabular}[t]{c |ccc}
  {\rm basis} & 
  {\rm $n_x$} &
  {\rm  $n_y$} &
  {\rm  $n_z$} 
   \\ \hline  \hline
 $\lbrack x \rbrack_{row}$
 &
 $T_q$
 & 
 $-ZN_{q^{-1}}Z$
 & 
 $N_q$
 \\
 $\lbrack x \rbrack_{col}$
 &
 $T_q^t$
 &
 $N^t_q$
 &
 $-ZN^t_{q^{-1}}Z$
 \\
 $\lbrack x \rbrack^{inv}_{row}$
 &
 $-T_{q^{-1}}$
 &
 $-N_{q^{-1}}$
 & 
 $ZN_qZ$
 \\
 $\lbrack x \rbrack^{inv}_{col}$
 &
 $-T^t_{q^{-1}}$
 &
 $ZN^t_qZ$
  & 
 $-N^t_{q^{-1}}$
 \\
 \hline
 $\lbrack y \rbrack_{row}$
 & 
 $N_q$
 &
 $T_q$ 
 &
 $-ZN_{q^{-1}}Z$
 \\
 $\lbrack y \rbrack_{col}$
 &
 $-ZN^t_{q^{-1}}Z$
& 
 $T^t_q$ 
 &
 $N^t_q$
 \\
 $\lbrack y \rbrack^{inv}_{row}$
 & 
 $ZN_qZ$
 &
 $-T_{q^{-1}}$ 
 &
 $-N_{q^{-1}}$
 \\
 $\lbrack y \rbrack^{inv}_{col}$
  & 
 $-N^t_{q^{-1}}$
 &
 $-T^t_{q^{-1}}$ 
 &
 $ZN^t_qZ$
 \\
 \hline
 $\lbrack z \rbrack_{row}$
 & 
 $-ZN_{q^{-1}}Z$
 & 
 $N_q$
 &
 $T_q$ 
 \\
 $\lbrack z \rbrack_{col}$
 &
 $N^t_q$
 &
 $-ZN^t_{q^{-1}}Z$
& 
 $T^t_q$ 
 \\
 $\lbrack z \rbrack^{inv}_{row}$
 &
 $-N_{q^{-1}}$
 & 
 $ZN_qZ$
 &
 $-T_{q^{-1}}$ 
 \\
 $\lbrack z \rbrack^{inv}_{col}$
 &
 $ZN^t_qZ$
  & 
 $-N^t_{q^{-1}}$
 &
 $-T^t_{q^{-1}}$ 
     \end{tabular}}
     \medskip

\end{theorem}
\noindent {\it Proof:} We first verify the data for the
middle third of the table.
Consider the matrices that represent
 $n_x,n_y, n_z$ with respect to
a $\lbrack y \rbrack_{row} $ basis for $V$.
For $n_x,n_z$ these matrices are obtained using 
Lemma
\ref{lem:pp} and the construction.
Concerning $n_y$,
recall from Definition
\ref{def:nnn} that $n_y = q^{-1}(1-xz)(q-q^{-1})^{-1}$.
By Theorem
\ref{thm:vmat}  
the matrix $E_q$ (resp. $ZE_{q^{-1}}Z$) represents 
$x$ (resp. $z$) with respect to a 
 $\lbrack y \rbrack_{row} $ basis for $V$.
One verifies using Definition
\ref{def:tq}
that 
 $T_q = q^{-1}(1-E_qZE_{q^{-1}}Z)(q-q^{-1})^{-1}$.
By these comments the matrix $T_q$ 
represents $n_y$ with respect to a
 $\lbrack y \rbrack_{row} $ basis for $V$.
We have obtained the matrices that represent
 $n_x,n_y, n_z$ with respect to
a $\lbrack y \rbrack_{row} $ basis for $V$.
For these matrices conjugate by $Z$ 
and use Note
\ref{note:ztz}
to get the
matrices that represent 
 $n_x,n_y, n_z$ with respect to
a $\lbrack y \rbrack^{inv}_{row} $ basis for $V$.
For these matrices replace $q$ by $q^{-1}$ to get
the matrices that represent 
 $n_x,n_y, n_z$ with respect to
a $\lbrack y \rbrack^{inv}_{row} $ basis for $V^*$.
For these matrices take $-1$ times the transpose and invoke
Lemmas
\ref{lem:nxnynz},
\ref{lem:btrans} to get the matrices that represent
 $n_x,n_y, n_z$ with respect to
a $\lbrack y \rbrack_{col} $ basis for $V$.
For these matrices conjugate by $Z$ and use
Note
\ref{note:ztz}
to get the matrices
that represent 
 $n_x,n_y, n_z$ with respect to
a $\lbrack y \rbrack^{inv}_{col} $ basis for $V$.
We have now verified the data for the middle third of
the table. To verify the rest of the table use Lemma
\ref{lem:om}.
\hfill $\Box$ \\

\begin{theorem}
\label{thm:nxmats}
Consider the elements
$n_x$, $n_y$, $n_z$ of
 $U_{q^{-1}}(\mathfrak{sl}_2)$.
In the table below we display the matrices that represent
these elements with respect to the twelve bases for $V^*$
from
{\rm (\ref{eq:basisr1})--(\ref{eq:basisr4})}.
\medskip

\centerline{
\begin{tabular}[t]{c |ccc}
  {\rm basis} & 
  {\rm $n_x$} &
  {\rm  $n_y$} &
  {\rm  $n_z$} 
   \\ \hline  \hline
 $\lbrack x \rbrack_{row}$
 &
 $T_{q^{-1}}$
 & 
 $-ZN_{q}Z$
 & 
 $N_{q^{-1}}$
 \\
 $\lbrack x \rbrack_{col}$
 &
 $T^t_{q^{-1}}$
 &
 $N^t_{q^{-1}}$
 &
 $-ZN^t_{q}Z$
 \\
 $\lbrack x \rbrack^{inv}_{row}$
 &
 $-T_q$
 &
 $-N_{q}$
 & 
 $ZN_{q^{-1}}Z$
 \\
 $\lbrack x \rbrack^{inv}_{col}$
 &
 $-T^t_q$
 &
 $ZN^t_{q^{-1}}Z$
  & 
 $-N^t_{q}$
 \\
 \hline
 $\lbrack y \rbrack_{row}$
 & 
 $N_{q^{-1}}$
 &
 $T_{q^{-1}}$
 &
 $-ZN_{q}Z$
 \\
 $\lbrack y \rbrack_{col}$
 &
 $-ZN^t_{q}Z$
& 
 $T^t_{q^{-1}}$
 &
 $N^t_{q^{-1}}$
 \\
 $\lbrack y \rbrack^{inv}_{row}$
 & 
 $ZN_{q^{-1}}Z$
 &
 $-T_q$
 &
 $-N_{q}$
 \\
 $\lbrack y \rbrack^{inv}_{col}$
  & 
 $-N^t_{q}$
 &
 $-T^t_q$
 &
 $ZN^t_{q^{-1}}Z$
 \\
 \hline
 $\lbrack z \rbrack_{row}$
 & 
 $-ZN_{q}Z$
 & 
 $N_{q^{-1}}$
 &
 $T_{q^{-1}}$
 \\
 $\lbrack z \rbrack_{col}$
 &
 $N^t_{q^{-1}}$
 &
 $-ZN^t_{q}Z$
& 
 $T^t_{q^{-1}}$
 \\
 $\lbrack z \rbrack^{inv}_{row}$
 &
 $-N_{q}$
 & 
 $ZN_{q^{-1}}Z$
 &
 $-T_q$
 \\
 $\lbrack z \rbrack^{inv}_{col}$
 &
 $ZN^t_{q^{-1}}Z$
  & 
 $-N^t_{q}$
 &
 $-T^t_q$
     \end{tabular}}
     \medskip

\end{theorem}
\noindent {\it Proof:} In the table of Theorem
\ref{thm:nxmat}
replace $q$ by $q^{-1}$.
\hfill $\Box$ \\



\section{Comments on the bilinear form}

\noindent We continue to discuss the
$U_q(\mathfrak{sl}_2)$-module $V$ and
the 
$U_{q^{-1}}(\mathfrak{sl}_2)$-module $V^*$.
Recall the twelve bases
(\ref{eq:basisr1})--(\ref{eq:basisr4})
for $V$ and $V^*$. In Section 15 we will compute the
transition matrices between certain pairs of bases among
these twelve. Before we get to this, it is convenient
to establish a few facts about the bilinear form
$(\,,\,)$ from Definition
\ref{def:bil}.

\medskip
\noindent We recall some notation. 
For integers $n\geq i \geq 0$ define
\begin{eqnarray*}
\left[\begin{array}{c} n  \\
i \end{array} \right] &=& 
\frac{\lbrack n \rbrack^!}{\lbrack i \rbrack^! \lbrack n-i\rbrack^! }.
\end{eqnarray*}

\begin{lemma}
\label{lem:innerp}
Pick $\xi \in \lbrace x,y,z\rbrace$.
Let
$\lbrace u_i\rbrace_{i=0}^d$ 
denote a  $\lbrack \xi \rbrack_{row} $ basis for
$V$
and let $\lbrace v_i\rbrace_{i=0}^d$ 
denote a $\lbrack \xi \rbrack_{row}$ basis for
$V^*$. Then 
\begin{eqnarray}
(u_r, v_s) = \delta_{r+s,d}(-1)^r q^{r(d-1)}
\left[\begin{array}{c} d  \\
r \end{array} \right] (u_0, v_d)
\label{eq:udv0}
\end{eqnarray}
for $0 \leq r,s\leq d$.
\end{lemma}
\noindent {\it Proof:}  By Lemma
\ref{lem:om}, without loss we may assume $\xi=y$.
If $r+s\not=d$ then
$(u_r, v_s) = 0$ by
Lemma \ref{lem:dualdec}.
By Proposition
\ref{prop:howuquqirel}
we have
\begin{eqnarray}
(n_z u_i, v_{d-i+1})= -(u_i, n_z v_{d-i+1})
\label{eq:bilrec}
\end{eqnarray}
for $1 \leq i \leq d$. The action of $n_z$
on 
$\lbrace u_i\rbrace_{i=0}^d$  is given in
Theorem
\ref{thm:nxmat},
and the action of 
 $n_z$
on 
$\lbrace v_i\rbrace_{i=0}^d$  is given in
Theorem
\ref{thm:nxmats}.
Evaluating 
(\ref{eq:bilrec}) using this data we find
\begin{eqnarray*}
q^{d-i}\lbrack d-i+1\rbrack ( u_{i-1}, v_{d-i+1})= 
-q^{1-i}\lbrack i \rbrack (u_i, v_{d-i})
\qquad \qquad (1 \leq i \leq d).
\end{eqnarray*}
Solving this recursion we find
\begin{eqnarray*}
(u_r, v_{d-r}) = (-1)^{r}q^{r(d-1)}
\left[\begin{array}{c} d  \\
r \end{array} \right] (u_0, v_d)
\qquad \qquad (0 \leq r \leq d).
\end{eqnarray*}
The result follows.
\hfill $\Box$ \\

\begin{corollary}
\label{cor:0d}
With reference to Lemma
\ref{lem:innerp},
\begin{eqnarray}
(u_d, v_0) = (-1)^d q^{d(d-1)}
(u_0, v_d).
\label{eq:0d}
\end{eqnarray}
\end{corollary}
\noindent {\it Proof:} In 
(\ref{eq:udv0})
set $r=d$ and $s=0$.
\hfill $\Box$ \\

\section{A normalization for the twelve bases}

\noindent 
We continue to discuss the
$U_q(\mathfrak{sl}_2)$-module $V$ and
the 
$U_{q^{-1}}(\mathfrak{sl}_2)$-module $V^*$.
Recall the twelve bases
(\ref{eq:basisr1})--(\ref{eq:basisr4})
for $V$ and $V^*$. In Section 15 we will compute the
transition matrices between certain pairs of bases among
these twelve. In order to do this efficiently
we first normalize our bases.

\begin{definition}
\label{def:eta}
\rm
For $\xi \in \lbrace x,y,z \rbrace$ let 
$\eta_\xi$ (resp. $\eta^*_\xi$) denote a nonzero vector in
$n_\xi^dV$ (resp. 
$n_\xi^dV^*$).
\end{definition}

\begin{lemma}  The following {\rm (i), (ii)} hold.
\begin{enumerate}
\item[\rm (i)]
For distinct $u,v \in \lbrace x,y,z\rbrace$
we have $(\eta_u, \eta^*_v)\not=0$.
\item[\rm (ii)]
Assume $d\geq 1$.
Then for $u \in \lbrace x,y,z\rbrace$
we have
$(\eta_u, \eta^*_u)=0$.
\end{enumerate}
\end{lemma}
\noindent {\it Proof:} (i)
The vector $\eta_u$ is a basis for $n_u^dV$.
By Lemma \ref{lem:dual}
the orthogonal complement of
$n_u^dV$ is $n_uV^*$. By
Lemma
\ref{lem:nnnpre} and
Definition
\ref{def:eta}
$\eta^*_v \not\in n_uV^*$. Therefore
$(\eta_u,\eta^*_v)\not=0$.
\\
\noindent (ii)
We mentioned above that
the orthogonal complement of
$n_u^dV$ is $n_uV^*$.
We assume $d\geq 1$ so
$n_uV^*$ contains 
$n^d_uV^*$.
Therefore 
$n_u^dV$ and
$n_u^dV^*$ are orthogonal
so
$(\eta_u, \eta^*_u)=0$.
\hfill $\Box$ \\

\begin{lemma}
\label{lem:urownorm}
Pick $\xi \in \lbrace x,y,z\rbrace $.
There exists a unique basis
$\lbrace v_i\rbrace_{i=0}^d$ for $V$ such that:
\begin{enumerate}
\item[\rm (i)] for $0 \leq i \leq d$ the vector
$v_i$ is contained in component $i$ of the decomposition
$\lbrack \xi \rbrack$;
\item[\rm (ii)] $\eta_\xi = \sum_{i=0}^d v_i$.
\end{enumerate}
\end{lemma}
\noindent {\it Proof:}
Concerning existence,
let $\lbrace u_i\rbrace_{i=0}^d$ denote a 
 $\lbrack \xi \rbrack_{row}$ basis for $V$.
Then $\sum_{i=0}^d u_i$ is contained in
$n^d_\xi V$ and is therefore a scalar multiple of
$\eta_\xi $. Call this scalar $\kappa $ and observe
that $\kappa \not=0$. Define
$v_i=u_i/\kappa$ for $0 \leq i \leq d$.
Then
$\lbrace v_i\rbrace_{i=0}^d$ is the desired basis.
We have shown that the desired basis exists.
The uniqueness assertion is readily verified.
\hfill $\Box$ \\

\begin{definition}
\label{def:urownorm}
\rm
Pick $\xi \in \lbrace x,y,z\rbrace$.
Let $\lbrack \xi \rbrack_{row}$ denote the basis for $V$
that satisfies conditions
(i), (ii) of
Lemma
\ref{lem:urownorm}.
The basis 
$\lbrack \xi \rbrack_{row}$ for $V^*$
similarly defined, with $\eta_\xi$ replaced by
$\eta^*_\xi$ 
in Lemma
\ref{lem:urownorm}(ii).
The inversion of
$\lbrack \xi \rbrack_{row}$ is denoted
$\lbrack \xi \rbrack^{inv}_{row}$.
\end{definition}

\begin{lemma} 
\label{lem:endpt}
In the table below we give three bases for
$V$. For each basis we describe the components $0$ and $d$.
\medskip

\centerline{
\begin{tabular}[t]{c|c c}
   {\rm basis for $V$} &  {\rm component $0$}& {\rm component $d$}
   \\ \hline  \hline
$\lbrack x \rbrack_{row}$ &
$
\eta_y \frac{(\eta_x, \eta^*_z)}{(\eta_y, \eta^*_z)}
$
&
$
\eta_z \frac{(\eta_x, \eta^*_y)}{(\eta_z, \eta^*_y)}
$
\\
$\lbrack y \rbrack_{row}$ &
$
\eta_z \frac{(\eta_y, \eta^*_x)}{(\eta_z, \eta^*_x)}
$
&
$
\eta_x \frac{(\eta_y, \eta^*_z)}{(\eta_x, \eta^*_z)}
$
\\
$\lbrack z \rbrack_{row}$ &
$
\eta_x \frac{(\eta_z, \eta^*_y)}{(\eta_x, \eta^*_y)}
$
&
$
\eta_y \frac{(\eta_z, \eta^*_x)}{(\eta_y, \eta^*_x)}
$
     \end{tabular}}
     \medskip

\end{lemma}
\noindent {\it Proof:}
Denote the basis $\lbrack x \rbrack_{row} $
by $\lbrace v_i\rbrace_{i=0}^d$.
Recall from Lemma 
\ref{lem:urownorm}(i)  that
for $0 \leq i \leq d$ the vector $v_i$ is
contained in component $i$ of 
the decomposition $\lbrack x \rbrack$.
 Component $0$ of 
 $\lbrack x \rbrack$
(resp. 
 component $d$ of 
 $\lbrack x \rbrack$)
 is equal to 
$n_y^dV$ (resp. 
$n_z^dV$)
and is therefore spanned by
 $\eta_y$ (resp. $\eta_z$).
Consequently there exist $\alpha,\beta \in \F$ such that
$v_0 = \alpha \eta_y$ and
$v_d = \beta \eta_z$.
By Lemma 
\ref{lem:urownorm}(ii) 
$\eta_x = \sum_{i=0}^d v_i$.
Using Lemma
\ref{lem:nnnpre} and
Lemma \ref{lem:dual}
we find
$(v_i, \eta^*_z)=0$ for $1 \leq i \leq d$.
Therefore 
\begin{eqnarray*}
(\eta_x, \eta^*_z) = \sum_{i=0}^d (v_i, \eta^*_z) = 
(v_0,\eta^*_z) = \alpha (\eta_y,\eta^*_z)
\end{eqnarray*}
so $\alpha = 
(\eta_x, \eta^*_z)/
(\eta_y,\eta^*_z)$. 
Using Lemma
\ref{lem:nnnpre} and
Lemma \ref{lem:dual}
we find
$(v_i, \eta^*_y)=0$ for $0 \leq i \leq d-1$.
Therefore
\begin{eqnarray*}
(\eta_x, \eta^*_y) = \sum_{i=0}^d (v_i, \eta^*_y) = 
(v_d,\eta^*_y) = \beta (\eta_z,\eta^*_y)
\end{eqnarray*}
so $\beta = 
(\eta_x, \eta^*_y)/
(\eta_z,\eta^*_y)$.
We have verified our assertions for the basis 
$\lbrack x \rbrack_{row}$. To verify our remaining assertions
use Lemma \ref{lem:om}.
\hfill $\Box$ \\

\begin{lemma} 
\label{lem:endptdual}
In the table below we give three bases for
$V^*$. For each basis we describe the components $0$ and $d$.
\medskip

\centerline{
\begin{tabular}[t]{c|c c}
   {\rm basis for $V^*$} & {\rm  component $0$}&{\rm component $d$}
   \\ \hline  \hline
$\lbrack x \rbrack_{row}$ &
$
\eta^*_y \frac{(\eta_z, \eta^*_x)}{(\eta_z, \eta^*_y)}
$
&
$
\eta^*_z \frac{(\eta_y, \eta^*_x)}{(\eta_y, \eta^*_z)}
$
\\
$\lbrack y \rbrack_{row}$ &
$
\eta^*_z \frac{(\eta_x, \eta^*_y)}{(\eta_x, \eta^*_z)}
$
&
$
\eta^*_x \frac{(\eta_z, \eta^*_y)}{(\eta_z, \eta^*_x)}
$
\\
$\lbrack z \rbrack_{row}$ &
$
\eta^*_x \frac{(\eta_y, \eta^*_z)}{(\eta_y, \eta^*_x)}
$
&
$
\eta^*_y \frac{(\eta_x, \eta^*_z)}{(\eta_x, \eta^*_y)}
$
     \end{tabular}}
     \medskip

\end{lemma}
\noindent {\it Proof:}
Similar to the proof of
Lemma \ref{lem:endpt}.
\hfill $\Box$ \\

\begin{definition}
\label{def:ucolnorm}
\rm
For $\xi \in \lbrace x,y,z\rbrace$ let
$\lbrack \xi \rbrack_{col}$ denote the
basis for $V$ (resp. $V^*$) that is dual to
the basis 
$\lbrack \xi \rbrack^{inv}_{row}$ for $V^*$ (resp. $V$).
The inversion of
$\lbrack \xi \rbrack_{col}$ is denoted
$\lbrack \xi \rbrack^{inv}_{col}$.
\end{definition}


\begin{lemma}
\label{lem:endptcolnorm}
In the table below we give three bases for
$V$. For each basis we describe the components $0$ and $d$.
\medskip

\centerline{
\begin{tabular}[t]{c|c c}
   {\rm basis for $V$} &  {\rm component $0$}& {\rm component $d$}
   \\ \hline  \hline
$\lbrack x \rbrack_{col}$ &
$
\frac{\eta_y}{(\eta_y, \eta^*_x)}
$
&
$
\frac{\eta_z }{(\eta_z, \eta^*_x)}
$
\\
$\lbrack y \rbrack_{col}$ &
$
\frac{\eta_z}{(\eta_z, \eta^*_y)}
$
&
$
\frac{\eta_x }{(\eta_x, \eta^*_y)}
$
\\
$\lbrack z \rbrack_{col}$ &
$
\frac{\eta_x}{(\eta_x, \eta^*_z)}
$
&
$
\frac{\eta_y }{(\eta_y, \eta^*_z)}
$
     \end{tabular}}
     \medskip

\end{lemma}
\noindent {\it Proof:}
For the vector space $V$ consider the basis 
 $\lbrack x \rbrack_{col}$ and the decomposition
 $\lbrack x \rbrack$.
By Lemma
\ref{lem:inducecol},
 $\lbrack x \rbrack_{col}$ induces
 $\lbrack x \rbrack$.
Component $0$ (resp. component $d$) of
 $\lbrack x \rbrack$ is equal to
$n^d_yV$ (resp.
$n^d_zV$)
and is therefore spanned by
$\eta_y$ (resp. $\eta_z$).
Therefore,
component $0$ (resp. component $d$) of
 $\lbrack x \rbrack_{col}$ is a scalar multiple
of 
$\eta_y$ (resp. $\eta_z$).
To find the scalars, use
the fact that component 
$0$ (resp. component $d$)
of 
 $\lbrack x \rbrack_{col}$ 
has inner product 1 with 
component $d$ (resp. component $0$) 
of the basis 
 $\lbrack x \rbrack_{row}$ for $V^*$.
These components of the basis 
 $\lbrack x \rbrack_{row}$ for $V^*$
are given in 
Lemma
\ref{lem:endptdual}.
By these comments we routinely
verify our assertions for
the basis $\lbrack x \rbrack_{col}$. To verify
our remaining
assertions use
Lemma
\ref{lem:om}.
\hfill $\Box$ \\

\begin{lemma} In the table below we give three bases for
$V^*$. For each basis we describe the components $0$ and $d$.
\medskip

\centerline{
\begin{tabular}[t]{c|c c}
   {\rm basis for $V^*$} &  {\rm component $0$}& {\rm component $d$}
   \\ \hline  \hline
$\lbrack x \rbrack_{col}$ &
$
\frac{\eta^*_y}{(\eta_x, \eta^*_y)}
$
&
$
\frac{\eta^*_z }{(\eta_x, \eta^*_z)}
$
\\
$\lbrack y \rbrack_{col}$ &
$
\frac{\eta^*_z}{(\eta_y, \eta^*_z)}
$
&
$
\frac{\eta^*_x }{(\eta_y, \eta^*_x)}
$
\\
$\lbrack z \rbrack_{col}$ &
$
\frac{\eta^*_x}{(\eta_z, \eta^*_x)}
$
&
$
\frac{\eta^*_y }{(\eta_z, \eta^*_y)}
$
     \end{tabular}}
     \medskip

\end{lemma}
\noindent {\it Proof:}
Similar to the proof of
Lemma
\ref{lem:endptcolnorm}.
\hfill $\Box$ \\



\begin{lemma} 
\label{lem:whatdualwhat2}
Pick $\xi \in \lbrace x,y,z\rbrace$.
For the table below, in each row
 we display a basis for $V$ and a basis for $V^*$.
These bases are dual.
\medskip

\centerline{
\begin{tabular}[t]{c c}
   {\rm basis for $V$} & 
   {\rm  basis for $V^*$}
   \\ \hline  \hline
  $\lbrack \xi \rbrack_{row}$   & $\lbrack \xi \rbrack^{inv}_{col} $
  \\ 
  $\lbrack \xi \rbrack_{col}$   & $\lbrack \xi \rbrack^{inv}_{row}$ 
   \\ 
  $\lbrack \xi \rbrack^{inv}_{row}$   & $\lbrack \xi \rbrack_{col} $
  \\ 
  $\lbrack \xi \rbrack^{inv}_{col}$   & $\lbrack \xi \rbrack_{row}$ 
     \end{tabular}}
     \medskip

\end{lemma}
\noindent {\it Proof:}
By Definition
\ref{def:ucolnorm}
and the meaning of inversion.
\hfill $\Box$ \\


\noindent We now consider how the scalars
\begin{eqnarray*}
(\eta_u, \eta^*_v) \qquad \qquad u,v \in \lbrace x,y,z\rbrace,
\quad u\not=v
\end{eqnarray*}
are related.

\begin{proposition} 
\label{prop:sixfactors}
We have
\begin{eqnarray*}
\frac{
(\eta_x, \eta^*_y)
(\eta_y, \eta^*_z)
(\eta_z, \eta^*_x)
}
{
(\eta_x, \eta^*_z)
(\eta_y, \eta^*_x)
(\eta_z, \eta^*_y)
}
= (-1)^d q^{d(d-1)}.
\end{eqnarray*}
\end{proposition}
\noindent {\it Proof:} 
Let $\lbrace u_i\rbrace_{i=0}^d$ denote the basis
$\lbrack y \rbrack_{row}$ for $V$
and let 
 $\lbrace v_i\rbrace_{i=0}^d$ denote the basis
$\lbrack y \rbrack_{row}$ for $V^*$. 
These bases satisfy
(\ref{eq:0d}).
By 
Lemma
\ref{lem:endpt},
\begin{eqnarray}
u_0 =
\eta_z \frac{(\eta_y, \eta^*_x)}{(\eta_z, \eta^*_x)},
\qquad \qquad 
u_d = 
\eta_x \frac{(\eta_y, \eta^*_z)}{(\eta_x, \eta^*_z)}.
\label{eq:u0ud}
\end{eqnarray}
By 
Lemma
\ref{lem:endptdual},
\begin{eqnarray}
v_0 = 
\eta^*_z \frac{(\eta_x, \eta^*_y)}{(\eta_x, \eta^*_z)},
\qquad \qquad 
v_d = 
\eta^*_x \frac{(\eta_z, \eta^*_y)}{(\eta_z, \eta^*_x)}.
\label{eq:v0vd}
\end{eqnarray}
In the equation 
(\ref{eq:0d}), eliminate
$u_0, u_d$ using
(\ref{eq:u0ud}) and 
eliminate
$v_0, v_d$ using
(\ref{eq:v0vd}).
The result follows
after a routine simplification.
\hfill $\Box$ \\

\begin{note}\rm
By Proposition
\ref{prop:sixfactors} the scalars
\begin{eqnarray*}
(\eta_u, \eta^*_v) \qquad \qquad u,v \in \lbrace x,y,z\rbrace,
\quad u\not=v
\end{eqnarray*}
are determined by the sequence
\begin{eqnarray}
(\eta_x,\eta^*_y), \quad 
(\eta_y,\eta^*_z), \quad 
(\eta_z,\eta^*_x), \quad 
(\eta_y,\eta^*_x), \quad 
(\eta_z,\eta^*_y).
\label{eq:5list}
\end{eqnarray}
The scalars 
(\ref{eq:5list}) are ``free'' in the following sense.
Given a sequence $\theta$ of five nonzero scalars in $\F$,
there exist vectors $\eta_x, \eta_y, \eta_z$ and
$\eta^*_x, \eta^*_y, \eta^*_z$ as in Definition
\ref{def:eta} such that the sequence 
(\ref{eq:5list}) is equal to $\theta$.
\end{note}


\section{The twelve normalized bases in closed form}

\noindent We continue to discuss the
$U_q(\mathfrak{sl}_2)$-module $V$ and
the 
$U_{q^{-1}}(\mathfrak{sl}_2)$-module $V^*$.
Recall the 
twelve bases
(\ref{eq:basisr1})--(\ref{eq:basisr4}) for $V$ and
$V^*$, normalized as in Section 13.
In this section we display these normalized
bases in closed form.

\newpage
\begin{theorem}
\label{thm:closedform}
In the table below we list twelve bases for
$V$.
For each basis 
we
display component $i$ for $0 \leq i \leq d$. We give two versions.
\medskip

\centerline{
\begin{tabular}[t]{c|c c}
   {\rm basis }
   &  {\rm component $i$ (version 1)}
   &  {\rm component $i$ (version 2)}
   \\ \hline  \hline
$
\lbrack x \rbrack_{row} 
$
&
$
\frac{q^{\binom{i}{2}}}{\lbrack i \rbrack^!}
\frac{
(\eta_x,\eta^*_z)}
{(\eta_y,\eta^*_z)}
n_z^i \eta_y
$
&
$
\frac{(-1)^{d-i}q^{-\binom{d-i}{2}}}{\lbrack d-i \rbrack^!}
\frac{
(\eta_x,\eta^*_y)}
{(\eta_z,\eta^*_y)}
n_y^{d-i} \eta_z
$
\\
\\
$
\lbrack x \rbrack_{col} 
$
&
$
\frac{(-1)^i 
\lbrack d-i\rbrack^!
q^{i(1-d)+\binom{i}{2}} 
}
{\lbrack d \rbrack^! (\eta_y,\eta^*_x)}
n_z^i \eta_y
$
&
$
\frac{\lbrack i \rbrack^! q^{(d-i)(d-1)-\binom{d-i}{2}}
}
{\lbrack d \rbrack^! (\eta_z, \eta^*_x)
}
n_y^{d-i} \eta_z
$
\\
\\
$\lbrack x \rbrack^{inv}_{row}$
&
$
\frac{
(-1)^i q^{-\binom{i}{2}}
}
{\lbrack i \rbrack^!}
\frac{(\eta_x,\eta_y^*)}
{(\eta_z,\eta^*_y)} n_y^i \eta_z
$
&
$
\frac{q^{\binom{d-i}{2}}}{\lbrack d-i \rbrack^!}
\frac{
(\eta_x,\eta^*_z)}
{(\eta_y,\eta^*_z)}
n_z^{d-i} \eta_y
$
\\
\\
$\lbrack x \rbrack^{inv}_{col}$
&
$
\frac{\lbrack d-i \rbrack^! q^{i(d-1)-\binom{i}{2}}
}
{\lbrack d \rbrack^! (\eta_z, \eta^*_x)
}
n_y^i \eta_z
$
&
$
\frac{(-1)^{d-i} 
\lbrack i\rbrack^!
q^{(d-i)(1-d)+\binom{d-i}{2}} 
}
{\lbrack d \rbrack^! (\eta_y,\eta^*_x)}
n_z^{d-i} \eta_y
$
\\
\\
\hline
$
\lbrack y \rbrack_{row} 
$
&
$
\frac{q^{\binom{i}{2}}}{\lbrack i \rbrack^!}
\frac{
(\eta_y,\eta^*_x)}
{(\eta_z,\eta^*_x)}
n_x^i \eta_z
$
&
$
\frac{(-1)^{d-i}q^{-\binom{d-i}{2}}}{\lbrack d-i \rbrack^!}
\frac{
(\eta_y,\eta^*_z)}
{(\eta_x,\eta^*_z)}
n_z^{d-i} \eta_x
$
\\
\\
$
\lbrack y \rbrack_{col} 
$
&
$
\frac{(-1)^i 
\lbrack d-i\rbrack^!
q^{i(1-d)+\binom{i}{2}} 
}
{\lbrack d \rbrack^! (\eta_z,\eta^*_y)}
n_x^i \eta_z
$
&
$
\frac{\lbrack i \rbrack^! q^{(d-i)(d-1)-\binom{d-i}{2}}
}
{\lbrack d \rbrack^! (\eta_x, \eta^*_y)
}
n_z^{d-i} \eta_x
$
\\
\\
$\lbrack y \rbrack^{inv}_{row}$
&
$
\frac{
(-1)^i q^{-\binom{i}{2}}
}
{\lbrack i \rbrack^!}
\frac{(\eta_y,\eta_z^*)}
{(\eta_x,\eta^*_z)} n_z^i \eta_x
$
&
$
\frac{q^{\binom{d-i}{2}}}{\lbrack d-i \rbrack^!}
\frac{
(\eta_y,\eta^*_x)}
{(\eta_z,\eta^*_x)}
n_x^{d-i} \eta_z
$
\\
\\
$\lbrack y \rbrack^{inv}_{col}$
&
$
\frac{\lbrack d-i \rbrack^! q^{i(d-1)-\binom{i}{2}}
}
{\lbrack d \rbrack^! (\eta_x, \eta^*_y)
}
n_z^i \eta_x
$
&
$
\frac{(-1)^{d-i} 
\lbrack i\rbrack^!
q^{(d-i)(1-d)+\binom{d-i}{2}} 
}
{\lbrack d \rbrack^! (\eta_z,\eta^*_y)}
n_x^{d-i} \eta_z
$
\\
\\
\hline
$
\lbrack z \rbrack_{row} 
$
&
$
\frac{q^{\binom{i}{2}}}{\lbrack i \rbrack^!}
\frac{
(\eta_z,\eta^*_y)}
{(\eta_x,\eta^*_y)}
n_y^i \eta_x
$
&
$
\frac{(-1)^{d-i}q^{-\binom{d-i}{2}}}{\lbrack d-i \rbrack^!}
\frac{
(\eta_z,\eta^*_x)}
{(\eta_y,\eta^*_x)}
n_x^{d-i} \eta_y
$
\\
\\
$
\lbrack z \rbrack_{col} 
$
&
$
\frac{(-1)^i 
\lbrack d-i\rbrack^!
q^{i(1-d)+\binom{i}{2}} 
}
{\lbrack d \rbrack^! (\eta_x,\eta^*_z)}
n_y^i \eta_x
$
&
$
\frac{\lbrack i \rbrack^! q^{(d-i)(d-1)-\binom{d-i}{2}}
}
{\lbrack d \rbrack^! (\eta_y, \eta^*_z)
}
n_x^{d-i} \eta_y
$
\\
\\
$\lbrack z \rbrack^{inv}_{row}$
&
$
\frac{
(-1)^i q^{-\binom{i}{2}}
}
{\lbrack i \rbrack^!}
\frac{(\eta_z,\eta_x^*)}
{(\eta_y,\eta^*_x)} n_x^i \eta_y
$
&
$
\frac{q^{\binom{d-i}{2}}}{\lbrack d-i \rbrack^!}
\frac{
(\eta_z,\eta^*_y)}
{(\eta_x,\eta^*_y)}
n_y^{d-i} \eta_x
$
\\
\\
$\lbrack z \rbrack^{inv}_{col}$
&
$
\frac{\lbrack d-i \rbrack^! q^{i(d-1)-\binom{i}{2}}
}
{\lbrack d \rbrack^! (\eta_y, \eta^*_z)
}
n_x^i \eta_y
$
&
$
\frac{(-1)^{d-i} 
\lbrack i\rbrack^!
q^{(d-i)(1-d)+\binom{d-i}{2}} 
}
{\lbrack d \rbrack^! (\eta_x,\eta^*_z)}
n_y^{d-i} \eta_x
$
\\
\end{tabular}}
     \medskip
\end{theorem}
\noindent {\it Proof:}
We first verify the data for the middle third of the table.
Consider the basis 
$\lbrack y \rbrack_{row} $ for $V$.
Denote this basis by
$\lbrace v_i\rbrace_{i=0}^d$.
The actions of $n_x$ and $n_z$  on 
$\lbrace v_i\rbrace_{i=0}^d$ 
are given in 
Lemma
\ref{lem:pp}(iii),(v).
 The information shows that
$v_i = q^{i-1}\lbrack i \rbrack^{-1} n_x v_{i-1}$ for 
$1 \leq i \leq d$, and
$v_i = -q^{i-d+1}\lbrack d-i\rbrack^{-1} n_z v_{i+1}$
for $0 \leq i \leq d-1$.
Therefore 
 both
\begin{eqnarray}
v_i = \frac{q^{\binom{i}{2}}
}
{
\lbrack i \rbrack^!
}
n^i _xv_0,
\qquad \qquad
v_i = (-1)^{d-i} \frac{q^{-\binom{d-i}{2}}}
{
\lbrack d-i \rbrack^!
}
n^{d-i}_z v_d
\label{eq:middle}
\end{eqnarray}
for $0 \leq i \leq d$.
By Lemma
\ref{lem:endpt},
\begin{eqnarray}
\label{eq:v0vd2}
v_0 = \frac{(\eta_y,\eta^*_x)}{(\eta_z,\eta^*_x)}\eta_z,
\qquad \qquad
v_d = \frac{(\eta_y,\eta^*_z)}{(\eta_x,\eta^*_z)}\eta_x.
\end{eqnarray}
In line
(\ref{eq:middle}), eliminate $v_0$ and $v_d$ using 
(\ref{eq:v0vd2}) to obtain the  
two descriptions for $v_i$ 
given in the table.
\\
\noindent
Next consider the basis 
$\lbrack y \rbrack_{col} $ for $V$.
Denote this basis by
$\lbrace v_i\rbrace_{i=0}^d$.
By Theorem
\ref{thm:nxmat}  the matrix
$-ZN^t_{q^{-1}}Z$ (resp. $N^t_q$)
represents 
$n_x$ (resp. $n_z$) with respect to
$\lbrace v_i\rbrace_{i=0}^d$. 
The entries of
$ZN^t_{q^{-1}}Z$ and $N^t_q$ are given in
Lemma \ref{lem:Nqentries}.
By these comments
$n_x v_{i-1} = -q^{d-i}\lbrack d-i+1\rbrack v_i$
and
$n_z v_i = q^{1-i} \lbrack i \rbrack v_{i-1}$
for $1 \leq i \leq d$.
Consequently
$v_i
=
-q^{i-d}\lbrack d-i+1 \rbrack^{-1}   n_x v_{i-1} $
for $1 \leq i \leq d$, and
$
v_{i} 
=
q^{i} \lbrack i+1\rbrack^{-1} n_z v_{i+1}
$
for $0 \leq i \leq d-1$.
Therefore both
\begin{eqnarray}
v_i = 
\frac{(-1)^i 
\lbrack d-i\rbrack^!
q^{i(1-d)+\binom{i}{2}} 
}
{\lbrack d \rbrack^! }
n_x^i v_0,
\qquad \qquad
v_i = 
\frac{\lbrack i \rbrack^! q^{(d-i)(d-1)-\binom{d-i}{2}}
}
{\lbrack d \rbrack^!
}
n_z^{d-i} v_d
\label{eq:middle2}
\end{eqnarray}
for $0 \leq i \leq d$.
By Lemma
\ref{lem:endptcolnorm},
\begin{eqnarray}
v_0 = \frac{\eta_z}{(\eta_z,\eta^*_y)},
\qquad \qquad
v_d = \frac{\eta_x}{(\eta_x,\eta^*_y)}.
\label{eq:v0vd3}
\end{eqnarray}
In line
(\ref{eq:middle2}), eliminate $v_0$ and $v_d$ using 
(\ref{eq:v0vd3}) to obtain the  
two descriptions for $v_i$ 
given in the table.
\\
\noindent Next consider the basis 
$\lbrack y\rbrack^{inv}_{row}$ for $V$.
For this basis component $i$ is equal to
component $d-i$ of the basis
$\lbrack y\rbrack_{row}$ for $V$.
\\
\noindent Next consider the basis 
$\lbrack y\rbrack^{inv}_{col}$ for $V$.
For this basis component $i$ is equal to
component $d-i$ of the basis
$\lbrack y\rbrack_{col}$ for $V$.
\\
\noindent We have now verified the data for the middle third
of the table. To verify the rest of the table
use Lemma
\ref{lem:om}.
\hfill $\Box$ \\

\newpage
\begin{theorem}
\label{thm:closedformdual}
In the table below we list twelve bases for
$V^*$.
For each basis we 
display component $i$ for
$0 \leq i \leq d$.
We give two versions.
\medskip

\centerline{
\begin{tabular}[t]{c|c c}
   {\rm basis }
   &  {\rm component $i$ (version 1)}
   &  {\rm component $i$ (version 2)}
   \\ \hline  \hline
$
\lbrack x \rbrack_{row} 
$
&
$
\frac{q^{-\binom{i}{2}}}{\lbrack i \rbrack^!}
\frac{
(\eta_z,\eta^*_x)}
{(\eta_z,\eta^*_y)}
n_z^i \eta^*_y
$
&
$
\frac{(-1)^{d-i}q^{\binom{d-i}{2}}}{\lbrack d-i \rbrack^!}
\frac{
(\eta_y,\eta^*_x)}
{(\eta_y,\eta^*_z)}
n_y^{d-i} \eta^*_z
$
\\
\\
$
\lbrack x \rbrack_{col} 
$
&
$
\frac{(-1)^i 
\lbrack d-i\rbrack^!
q^{i(d-1)-\binom{i}{2}} 
}
{\lbrack d \rbrack^! (\eta_x,\eta^*_y)}
n_z^i \eta^*_y
$
&
$
\frac{\lbrack i \rbrack^! q^{(d-i)(1-d)+\binom{d-i}{2}}
}
{\lbrack d \rbrack^! (\eta_x, \eta^*_z)
}
n_y^{d-i} \eta^*_z
$
\\
\\
$\lbrack x \rbrack^{inv}_{row}$
&
$
\frac{
(-1)^i q^{\binom{i}{2}}
}
{\lbrack i \rbrack^!}
\frac{(\eta_y,\eta_x^*)}
{(\eta_y,\eta^*_z)} n_y^i \eta^*_z
$
&
$
\frac{q^{-\binom{d-i}{2}}}{\lbrack d-i \rbrack^!}
\frac{
(\eta_z,\eta^*_x)}
{(\eta_z,\eta^*_y)}
n_z^{d-i} \eta^*_y
$
\\
\\
$\lbrack x \rbrack^{inv}_{col}$
&
$
\frac{\lbrack d-i \rbrack^! q^{i(1-d)+\binom{i}{2}}
}
{\lbrack d \rbrack^! (\eta_x, \eta^*_z)
}
n_y^i \eta^*_z
$
&
$
\frac{(-1)^{d-i} 
\lbrack i\rbrack^!
q^{(d-i)(d-1)-\binom{d-i}{2}} 
}
{\lbrack d \rbrack^! (\eta_x,\eta^*_y)}
n_z^{d-i} \eta^*_y
$
\\
\\
\hline
$
\lbrack y \rbrack_{row} 
$
&
$
\frac{q^{-\binom{i}{2}}}{\lbrack i \rbrack^!}
\frac{
(\eta_x,\eta^*_y)}
{(\eta_x,\eta^*_z)}
n_x^i \eta^*_z
$
&
$
\frac{(-1)^{d-i}q^{\binom{d-i}{2}}}{\lbrack d-i \rbrack^!}
\frac{
(\eta_z,\eta^*_y)}
{(\eta_z,\eta^*_x)}
n_z^{d-i} \eta^*_x
$
\\
\\
$
\lbrack y \rbrack_{col} 
$
&
$
\frac{(-1)^i 
\lbrack d-i\rbrack^!
q^{i(d-1)-\binom{i}{2}} 
}
{\lbrack d \rbrack^! (\eta_y,\eta^*_z)}
n_x^i \eta^*_z
$
&
$
\frac{\lbrack i \rbrack^! q^{(d-i)(1-d)+\binom{d-i}{2}}
}
{\lbrack d \rbrack^! (\eta_y, \eta^*_x)
}
n_z^{d-i} \eta^*_x
$
\\
\\
$\lbrack y \rbrack^{inv}_{row}$
&
$
\frac{
(-1)^i q^{\binom{i}{2}}
}
{\lbrack i \rbrack^!}
\frac{(\eta_z,\eta_y^*)}
{(\eta_z,\eta^*_x)} n_z^i \eta^*_x
$
&
$
\frac{q^{-\binom{d-i}{2}}}{\lbrack d-i \rbrack^!}
\frac{
(\eta_x,\eta^*_y)}
{(\eta_x,\eta^*_z)}
n_x^{d-i} \eta^*_z
$
\\
\\
$\lbrack y \rbrack^{inv}_{col}$
&
$
\frac{\lbrack d-i \rbrack^! q^{i(1-d)+\binom{i}{2}}
}
{\lbrack d \rbrack^! (\eta_y, \eta^*_x)
}
n_z^i \eta^*_x
$
&
$
\frac{(-1)^{d-i} 
\lbrack i\rbrack^!
q^{(d-i)(d-1)-\binom{d-i}{2}} 
}
{\lbrack d \rbrack^! (\eta_y,\eta^*_z)}
n_x^{d-i} \eta^*_z
$
\\
\\
\hline
$
\lbrack z \rbrack_{row} 
$
&
$
\frac{q^{-\binom{i}{2}}}{\lbrack i \rbrack^!}
\frac{
(\eta_y,\eta^*_z)}
{(\eta_y,\eta^*_x)}
n_y^i \eta^*_x
$
&
$
\frac{(-1)^{d-i}q^{\binom{d-i}{2}}}{\lbrack d-i \rbrack^!}
\frac{
(\eta_x,\eta^*_z)}
{(\eta_x,\eta^*_y)}
n_x^{d-i} \eta^*_y
$
\\
\\
$
\lbrack z \rbrack_{col} 
$
&
$
\frac{(-1)^i 
\lbrack d-i\rbrack^!
q^{i(d-1)-\binom{i}{2}} 
}
{\lbrack d \rbrack^! (\eta_z,\eta^*_x)}
n_y^i \eta^*_x
$
&
$
\frac{\lbrack i \rbrack^! q^{(d-i)(1-d)+\binom{d-i}{2}}
}
{\lbrack d \rbrack^! (\eta_z, \eta^*_y)
}
n_x^{d-i} \eta^*_y
$
\\
\\
$\lbrack z \rbrack^{inv}_{row}$
&
$
\frac{
(-1)^i q^{\binom{i}{2}}
}
{\lbrack i \rbrack^!}
\frac{(\eta_x,\eta_z^*)}
{(\eta_x,\eta^*_y)} n_x^i \eta^*_y
$
&
$
\frac{q^{-\binom{d-i}{2}}}{\lbrack d-i \rbrack^!}
\frac{
(\eta_y,\eta^*_z)}
{(\eta_y,\eta^*_x)}
n_y^{d-i} \eta^*_x
$
\\
\\
$\lbrack z \rbrack^{inv}_{col}$
&
$
\frac{\lbrack d-i \rbrack^! q^{i(1-d)+\binom{i}{2}}
}
{\lbrack d \rbrack^! (\eta_z, \eta^*_y)
}
n_x^i \eta^*_y
$
&
$
\frac{(-1)^{d-i} 
\lbrack i\rbrack^!
q^{(d-i)(d-1)-\binom{d-i}{2}} 
}
{\lbrack d \rbrack^! (\eta_z,\eta^*_x)}
n_y^{d-i} \eta^*_x
$
\\
\end{tabular}}
     \medskip
\end{theorem}
\noindent {\it Proof:}
In Theorem
\ref{thm:closedform} replace $q$ by $q^{-1}$.
Also replace 
$\eta_\xi$ by $\eta^*_\xi$ for 
$\xi \in \lbrace x,y,z\rbrace$,
and 
replace $(\eta_u,\eta^*_v)$
by
$(\eta_v,\eta^*_u)$
for distinct
$u,v \in \lbrace x,y,z\rbrace$.
\hfill $\Box$ \\

\noindent We finish this section with some comments.

\begin{corollary}
\label{cor:curious1}
The following hold:
\begin{eqnarray*}
&&
n^d_x \eta_y = \lbrack d \rbrack^! q^{-\binom{d}{2}}
\frac{(\eta_y,\eta^*_z)}{(\eta_x,\eta^*_z)} \eta_x,
\qquad \qquad
n^d_z \eta_y = (-1)^d \lbrack d \rbrack^! q^{\binom{d}{2}}
\frac{(\eta_y,\eta^*_x)}{(\eta_z,\eta^*_x)} \eta_z,
\\
&&
n^d_y \eta_z = \lbrack d \rbrack^! q^{-\binom{d}{2}}
\frac{(\eta_z,\eta^*_x)}{(\eta_y,\eta^*_x)} \eta_y,
\qquad \qquad
n^d_x \eta_z = (-1)^d \lbrack d \rbrack^! q^{\binom{d}{2}}
\frac{(\eta_z,\eta^*_y)}{(\eta_x,\eta^*_y)} \eta_x,
\\
&&
n^d_z \eta_x = \lbrack d \rbrack^! q^{-\binom{d}{2}}
\frac{(\eta_x,\eta^*_y)}{(\eta_z,\eta^*_y)} \eta_z,
\qquad \qquad
n^d_y \eta_x = (-1)^d \lbrack d \rbrack^! q^{\binom{d}{2}}
\frac{(\eta_x,\eta^*_z)}{(\eta_y,\eta^*_z)} \eta_y.
\end{eqnarray*}
\end{corollary}
\noindent {\it Proof:}
In the table of Theorem
\ref{thm:closedform}, set $i=0$ and
compare the two versions using
 Proposition
\ref{prop:sixfactors}.
\hfill $\Box$ \\

\begin{corollary}
The following hold:
\begin{eqnarray*}
&&
n^d_x \eta^*_y = \lbrack d \rbrack^! q^{\binom{d}{2}}
\frac{(\eta_z,\eta^*_y)}{(\eta_z,\eta^*_x)} \eta^*_x,
\qquad \qquad
n^d_z \eta^*_y = (-1)^d \lbrack d \rbrack^! q^{-\binom{d}{2}}
\frac{(\eta_x,\eta^*_y)}{(\eta_x,\eta^*_z)} \eta^*_z,
\\
&&
n^d_y \eta^*_z = \lbrack d \rbrack^! q^{\binom{d}{2}}
\frac{(\eta_x,\eta^*_z)}{(\eta_x,\eta^*_y)} \eta^*_y,
\qquad \qquad
n^d_x \eta^*_z = (-1)^d \lbrack d \rbrack^! q^{-\binom{d}{2}}
\frac{(\eta_y,\eta^*_z)}{(\eta_y,\eta^*_x)} \eta^*_x,
\\
&&
n^d_z \eta^*_x = \lbrack d \rbrack^! q^{\binom{d}{2}}
\frac{(\eta_y,\eta^*_x)}{(\eta_y,\eta^*_z)} \eta^*_z,
\qquad \qquad
n^d_y \eta^*_x = (-1)^d \lbrack d \rbrack^! q^{-\binom{d}{2}}
\frac{(\eta_z,\eta^*_x)}{(\eta_z,\eta^*_y)} \eta^*_y.
\end{eqnarray*}
\end{corollary}
\noindent {\it Proof:}
Similar to the proof of 
Corollary \ref{cor:curious1}.
\hfill $\Box$ \\


\section{Transition matrices between the twelve normalized bases}

\noindent We continue to discuss the
$U_q(\mathfrak{sl}_2)$-module $V$ and
the 
$U_{q^{-1}}(\mathfrak{sl}_2)$-module $V^*$.
Recall the twelve bases
(\ref{eq:basisr1})--(\ref{eq:basisr4})
for $V$ and $V^*$,  normalized as in
Section 13.
In this section we will compute the
transition matrices between certain pairs of bases among
these twelve. First we discuss a few terms. In this
discussion we focus on $V$; similar comments apply to $V^*$.

\medskip
\noindent 
Suppose we are given two bases for $V$, denoted
$\lbrace u_i\rbrace_{i=0}^d$
and 
$\lbrace v_i\rbrace_{i=0}^d$. By the {\it transition matrix from
$\lbrace u_i\rbrace_{i=0}^d$
to
$\lbrace v_i\rbrace_{i=0}^d$} we mean the
matrix $S \in {\rm Mat}_{d+1}(\F)$ such that
$v_j = \sum_{i=0}^d S_{ij} u_i$ for
$0 \leq j \leq d$.
Let $S$ denote the transition matrix from 
$\lbrace u_i\rbrace_{i=0}^d$
to
$\lbrace v_i\rbrace_{i=0}^d$.
Then $S^{-1}$ exists and equals the transition matrix
from
$\lbrace v_i\rbrace_{i=0}^d$
to
$\lbrace u_i\rbrace_{i=0}^d$.

\medskip
\noindent 
Let $\lbrace w_i\rbrace_{i=0}^d$ denote a basis for
$V$ and let $T$ denote the transition matrix from
$\lbrace v_i\rbrace_{i=0}^d$
to
$\lbrace w_i\rbrace_{i=0}^d$.
Then $ST$ is the transition matrix from
$\lbrace u_i\rbrace_{i=0}^d$
to
$\lbrace w_i\rbrace_{i=0}^d$.

\medskip
\noindent 
Let ${ A \in \rm End}(V)$ 
and let 
$ B$ denote the matrix in ${\rm Mat}_{d+1}(\F)$
that represents $A$ with respect to
$\lbrace u_i\rbrace_{i=0}^d$. Then
the matrix $S^{-1}BS$  represents $A$
with respect to 
$\lbrace v_i\rbrace_{i=0}^d$.

\medskip
\noindent 
Let $\lbrace u_i\rbrace_{i=0}^d$ and
 $\lbrace v_i\rbrace_{i=0}^d$
denote bases for
$V$.
Let $\lbrace u^*_i\rbrace_{i=0}^d$ 
 (resp. $\lbrace v^*_i\rbrace_{i=0}^d$)
denote the basis for
$V^*$ that is dual to
$\lbrace u_i\rbrace_{i=0}^d$ 
 (resp. $\lbrace v_i\rbrace_{i=0}^d$)
with respect to $(\,,\,)$.
Let $S$
denote the transition matrix from 
 $\lbrace u_i\rbrace_{i=0}^d$ to
 $\lbrace v_i\rbrace_{i=0}^d$.
Then $S^t$ is the transition matrix from
 $\lbrace v^*_i\rbrace_{i=0}^d$
to
$\lbrace u^*_i\rbrace_{i=0}^d$.

\medskip

\noindent Recall the matrix 
$Z$ 
from
Definition
\ref{def:z}.
Let $\lbrace v_i\rbrace_{i=0}^d$ denote a basis
for $V$ and consider the inverted basis
$\lbrace v_{d-i}\rbrace_{i=0}^d$.
 Then $Z$ is the transition matrix from
 $\lbrace v_i\rbrace_{i=0}^d$ to
$\lbrace v_{d-i}\rbrace_{i=0}^d$.

\begin{lemma} Consider the twelve bases
{\rm (\ref{eq:basisr1})--(\ref{eq:basisr4})}
for $V$ and $V^*$. For each basis, the transition matrix
to its inversion is equal to $Z$. In other words,
each of the following transition matrices is equal to $Z$:
\begin{eqnarray*}
&&\lbrack x \rbrack_{row} \to
 \lbrack x \rbrack^{inv}_{row},
\qquad 
\lbrack x \rbrack_{col} \to
 \lbrack x \rbrack^{inv}_{col},
\qquad 
 \lbrack x \rbrack^{inv}_{row} \to
\lbrack x \rbrack_{row},
\qquad 
\lbrack x \rbrack^{inv}_{col} \to
 \lbrack x \rbrack_{col},
\\
&&\lbrack y \rbrack_{row} \to
 \lbrack y \rbrack^{inv}_{row},
\qquad 
\lbrack y \rbrack_{col} \to
 \lbrack y \rbrack^{inv}_{col},
\qquad 
 \lbrack y \rbrack^{inv}_{row} \to
\lbrack y \rbrack_{row},
\qquad 
\lbrack y \rbrack^{inv}_{col} \to
 \lbrack y \rbrack_{col},
\\
&&\lbrack z \rbrack_{row} \to
 \lbrack z \rbrack^{inv}_{row},
\qquad 
\lbrack z \rbrack_{col} \to
 \lbrack z \rbrack^{inv}_{col},
\qquad 
 \lbrack z \rbrack^{inv}_{row} \to
\lbrack z \rbrack_{row},
\qquad 
\lbrack z \rbrack^{inv}_{col} \to
 \lbrack z \rbrack_{col}.
\end{eqnarray*}
\end{lemma}



\noindent Next we display some diagonal transition matrices.

\begin{theorem} 
\label{thm:diagtrans}
In the table below we display
some transition matrices between bases for $V$.
Each transition matrix is diagonal.
For $0 \leq i \leq d$ the $(i,i)$-entry is given.
\medskip

\centerline{
\begin{tabular}[t]{c|c}
   {\rm transition matrix}
   &  {\rm $(i,i)$-entry for $0 \leq i \leq d$ }
   \\ \hline  \hline
$\lbrack x \rbrack_{row} \to
\lbrack x \rbrack_{col}$
&
$(-1)^i q^{i(1-d)}
\left[\begin{array}{c} d  \\
i \end{array} \right]^{-1} 
\frac{(\eta_y, \eta^*_z)}{(\eta_y, \eta^*_x)(\eta_x, \eta^*_z)}
$
\\
$\lbrack x \rbrack_{col} \to
\lbrack x \rbrack_{row}$ &
$
(-1)^i q^{i(d-1)}
\left[\begin{array}{c} d  \\
i \end{array} \right] 
\frac{
(\eta_y, \eta^*_x)(\eta_x, \eta^*_z)
}{
(\eta_y, \eta^*_z)
}
$
\\
$\lbrack x \rbrack^{inv}_{row} \to
\lbrack x \rbrack^{inv}_{col}$
&
$
(-1)^{d-i} q^{(d-i)(1-d)}
\left[\begin{array}{c} d  \\
i \end{array} \right]^{-1} 
\frac{(\eta_y, \eta^*_z)}{(\eta_y, \eta^*_x)(\eta_x, \eta^*_z)}
$
\\
$\lbrack x \rbrack^{inv}_{col} \to
\lbrack x \rbrack^{inv}_{row}
$
&
$
(-1)^{d-i} q^{(d-i)(d-1)}
\left[\begin{array}{c} d  \\
i \end{array} \right] 
\frac{
(\eta_y, \eta^*_x)(\eta_x, \eta^*_z)
}{
(\eta_y, \eta^*_z)
}
$
\\
\hline
$\lbrack y \rbrack_{row} \to
\lbrack y \rbrack_{col}$
&
$(-1)^i q^{i(1-d)}
\left[\begin{array}{c} d  \\
i \end{array} \right]^{-1} 
\frac{(\eta_z, \eta^*_x)}{(\eta_z, \eta^*_y)(\eta_y, \eta^*_x)}
$
\\
$\lbrack y \rbrack_{col} \to
\lbrack y \rbrack_{row}$ &
$
(-1)^i q^{i(d-1)}
\left[\begin{array}{c} d  \\
i \end{array} \right] 
\frac{
(\eta_z, \eta^*_y)(\eta_y, \eta^*_x)
}{
(\eta_z, \eta^*_x)
}
$
\\
$\lbrack y \rbrack^{inv}_{row} \to
\lbrack y \rbrack^{inv}_{col}$
&
$
(-1)^{d-i} q^{(d-i)(1-d)}
\left[\begin{array}{c} d  \\
i \end{array} \right]^{-1} 
\frac{(\eta_z, \eta^*_x)}{(\eta_z, \eta^*_y)(\eta_y, \eta^*_x)}
$
\\
$\lbrack y \rbrack^{inv}_{col} \to
\lbrack y \rbrack^{inv}_{row}
$
&
$
(-1)^{d-i} q^{(d-i)(d-1)}
\left[\begin{array}{c} d  \\
i \end{array} \right] 
\frac{
(\eta_z, \eta^*_y)(\eta_y, \eta^*_x)
}{
(\eta_z, \eta^*_x)
}
$
\\
\hline
$\lbrack z \rbrack_{row} \to
\lbrack z \rbrack_{col}$
&
$(-1)^i q^{i(1-d)}
\left[\begin{array}{c} d  \\
i \end{array} \right]^{-1} 
\frac{(\eta_x, \eta^*_y)}{(\eta_x, \eta^*_z)(\eta_z, \eta^*_y)}
$
\\
$\lbrack z \rbrack_{col} \to
\lbrack z \rbrack_{row}$ &
$
(-1)^i q^{i(d-1)}
\left[\begin{array}{c} d  \\
i \end{array} \right] 
\frac{
(\eta_x, \eta^*_z)(\eta_z, \eta^*_y)
}{
(\eta_x, \eta^*_y)
}
$
\\
$\lbrack z \rbrack^{inv}_{row} \to
\lbrack z \rbrack^{inv}_{col}$
&
$
(-1)^{d-i} q^{(d-i)(1-d)}
\left[\begin{array}{c} d  \\
i \end{array} \right]^{-1} 
\frac{(\eta_x, \eta^*_y)}{(\eta_x, \eta^*_z)(\eta_z, \eta^*_y)}
$
\\
$\lbrack z \rbrack^{inv}_{col} \to
\lbrack z \rbrack^{inv}_{row}
$
&
$
(-1)^{d-i} q^{(d-i)(d-1)}
\left[\begin{array}{c} d  \\
i \end{array} \right] 
\frac{
(\eta_x, \eta^*_z)(\eta_z, \eta^*_y)
}{
(\eta_x, \eta^*_y)
}
$
\end{tabular}}
     \medskip

\end{theorem}
\noindent {\it Proof:}
We first verify the data for the middle third of the table.
Let
$\lbrace u_i\rbrace_{i=0}^d$ 
(resp. $\lbrace v_i\rbrace_{i=0}^d$) denote the
basis
$\lbrack y \rbrack_{row}$ 
(resp. $\lbrack y \rbrack_{col}$) for $V$.
By Theorem
\ref{thm:closedform},
\begin{eqnarray*}
u_i = 
\frac{q^{\binom{i}{2}}}{\lbrack i \rbrack^!}
\frac{
(\eta_y,\eta^*_x)}
{(\eta_z,\eta^*_x)}
n_x^i \eta_z,
\qquad \qquad 
v_i =
\frac{(-1)^i 
\lbrack d-i\rbrack^!
q^{i(1-d)+\binom{i}{2}} 
}
{\lbrack d \rbrack^! (\eta_z,\eta^*_y)}
n_x^i \eta_z.
\end{eqnarray*}
Comparing these we find
\begin{eqnarray*}
v_i = u_i  
(-1)^i q^{i(1-d)}
\left[\begin{array}{c} d  \\
i \end{array} \right]^{-1} 
\frac{(\eta_z, \eta^*_x)}{(\eta_z, \eta^*_y)(\eta_y, \eta^*_x)}.
\end{eqnarray*}
Therefore the transition matrix
$\lbrack y \rbrack_{row} \to
\lbrack y \rbrack_{col}$ is as claimed.
For this matrix take the inverse to get the
transition matrix
$\lbrack y \rbrack_{col} \to
\lbrack y \rbrack_{row}$.
To get the transition matrix
$\lbrack y \rbrack^{inv}_{row} \to
\lbrack y \rbrack^{inv}_{col}$, conjugate the
transition matrix
$\lbrack y \rbrack_{row} \to
\lbrack y \rbrack_{col}$ by the matrix $Z$ 
from Definition
\ref{def:z}.
To get the transition matrix
$\lbrack y \rbrack^{inv}_{col} \to
\lbrack y \rbrack^{inv}_{row}$, take the inverse of
the transition matrix
$\lbrack y \rbrack^{inv}_{row} \to
\lbrack y \rbrack^{inv}_{col}$.
We have verified the data for the middle third of the
table. To verify the rest of the table use Lemma
\ref{lem:om}.
\hfill $\Box$ \\

\begin{theorem} In the table below we display
some transition matrices between bases for $V^*$.
Each transition matrix is diagonal.
For $0 \leq i \leq d$ the $(i,i)$-entry is given.
\medskip

\centerline{
\begin{tabular}[t]{c|c}
   {\rm transition matrix}
   &  {\rm $(i,i)$-entry for $0 \leq i \leq d$ }
   \\ \hline  \hline
$\lbrack x \rbrack_{row} \to
\lbrack x \rbrack_{col}$
&
$(-1)^i q^{i(d-1)}
\left[\begin{array}{c} d  \\
i \end{array} \right]^{-1} 
\frac{(\eta_z, \eta^*_y)}
{
(\eta_z, \eta^*_x)
(\eta_x, \eta^*_y)
}
$
\\
$\lbrack x \rbrack_{col} \to
\lbrack x \rbrack_{row}$ &
$
(-1)^i q^{i(1-d)}
\left[\begin{array}{c} d  \\
i \end{array} \right] 
\frac{
(\eta_z, \eta^*_x)
(\eta_x, \eta^*_y)
}{
(\eta_z, \eta^*_y)
}
$
\\
$\lbrack x \rbrack^{inv}_{row} \to
\lbrack x \rbrack^{inv}_{col}$
&
$
(-1)^{d-i} q^{(d-i)(d-1)}
\left[\begin{array}{c} d  \\
i \end{array} \right]^{-1} 
\frac{(\eta_z, \eta^*_y)}{
(\eta_z, \eta^*_x)
(\eta_x, \eta^*_y)
}
$
\\
$\lbrack x \rbrack^{inv}_{col} \to
\lbrack x \rbrack^{inv}_{row}
$
&
$
(-1)^{d-i} q^{(d-i)(1-d)}
\left[\begin{array}{c} d  \\
i \end{array} \right] 
\frac{
(\eta_z, \eta^*_x)
(\eta_x, \eta^*_y)
}{
(\eta_z, \eta^*_y)
}
$
\\
\hline
$\lbrack y \rbrack_{row} \to
\lbrack y \rbrack_{col}$
&
$(-1)^i q^{i(d-1)}
\left[\begin{array}{c} d  \\
i \end{array} \right]^{-1} 
\frac{(\eta_x, \eta^*_z)}{
(\eta_x, \eta^*_y)
(\eta_y, \eta^*_z)
}
$
\\
$\lbrack y \rbrack_{col} \to
\lbrack y \rbrack_{row}$ &
$
(-1)^i q^{i(1-d)}
\left[\begin{array}{c} d  \\
i \end{array} \right] 
\frac{
(\eta_x, \eta^*_y)
(\eta_y, \eta^*_z)
}{
(\eta_x, \eta^*_z)
}
$
\\
$\lbrack y \rbrack^{inv}_{row} \to
\lbrack y \rbrack^{inv}_{col}$
&
$
(-1)^{d-i} q^{(d-i)(d-1)}
\left[\begin{array}{c} d  \\
i \end{array} \right]^{-1} 
\frac{(\eta_x, \eta^*_z)}{
(\eta_x, \eta^*_y)
(\eta_y, \eta^*_z)
}
$
\\
$\lbrack y \rbrack^{inv}_{col} \to
\lbrack y \rbrack^{inv}_{row}
$
&
$
(-1)^{d-i} q^{(d-i)(1-d)}
\left[\begin{array}{c} d  \\
i \end{array} \right] 
\frac{
(\eta_x, \eta^*_y)
(\eta_y, \eta^*_z)
}{
(\eta_x, \eta^*_z)
}
$
\\
\hline
$\lbrack z \rbrack_{row} \to
\lbrack z \rbrack_{col}$
&
$(-1)^i q^{i(d-1)}
\left[\begin{array}{c} d  \\
i \end{array} \right]^{-1} 
\frac{(\eta_y, \eta^*_x)}{
(\eta_y, \eta^*_z)
(\eta_z, \eta^*_x)
}
$
\\
$\lbrack z \rbrack_{col} \to
\lbrack z \rbrack_{row}$ &
$
(-1)^i q^{i(1-d)}
\left[\begin{array}{c} d  \\
i \end{array} \right] 
\frac{
(\eta_y, \eta^*_z)
(\eta_z, \eta^*_x)
}{
(\eta_y, \eta^*_x)
}
$
\\
$\lbrack z \rbrack^{inv}_{row} \to
\lbrack z \rbrack^{inv}_{col}$
&
$
(-1)^{d-i} q^{(d-i)(d-1)}
\left[\begin{array}{c} d  \\
i \end{array} \right]^{-1} 
\frac{(\eta_y, \eta^*_x)}{
(\eta_y, \eta^*_z)
(\eta_z, \eta^*_x)
}
$
\\
$\lbrack z \rbrack^{inv}_{col} \to
\lbrack z \rbrack^{inv}_{row}
$
&
$
(-1)^{d-i} q^{(d-i)(1-d)}
\left[\begin{array}{c} d  \\
i \end{array} \right] 
\frac{
(\eta_y, \eta^*_z)
(\eta_z, \eta^*_x)
}{
(\eta_y, \eta^*_x)
}
$
\end{tabular}}
     \medskip

\end{theorem}
\noindent {\it Proof:}
In Theorem
\ref{thm:diagtrans} replace $q$ by $q^{-1}$,
and also
replace $(\eta_u,\eta^*_v)$
by
 $(\eta_v,\eta^*_u)$
for distinct  
$u,v \in \lbrace x,y,z\rbrace$.
\hfill $\Box$ \\


\noindent Next we display some lower triangular transition 
matrices.

\newpage
\begin{theorem} 
\label{thm:lowert}
In the table below we display some transition matrices
between bases for $V$.
Each transition matrix is lower triangular.
For $0 \leq j \leq i \leq d$ the $(i,j)$-entry is given.
\medskip

\centerline{
\begin{tabular}[t]{c|c}
   {\rm transition matrix}
   &  {\rm $(i,j)$-entry for $0 \leq j\leq i\leq d$  }
   \\ \hline  \hline
$
\lbrack x \rbrack_{row} \to
\lbrack y \rbrack^{inv}_{row}
$
&
$
(-1)^j q^{j(1-i)} 
\left[\begin{array}{c} i  \\
j \end{array} \right] 
\frac{(\eta_y,\eta^*_z)}{(\eta_x, \eta^*_z)}
$
\\
$\lbrack x \rbrack_{col} \to
\lbrack y \rbrack^{inv}_{col}
$
&
$
(-1)^{d-i} q^{(i-d)(d-j-1)} 
\left[\begin{array}{c} d-j  \\
i-j \end{array} \right] 
\frac{(\eta_z,\eta^*_x)}{(\eta_z, \eta^*_y)}
$
\\
$\lbrack x \rbrack^{inv}_{row} \to 
\lbrack z \rbrack_{row}$
&
$
(-1)^j q^{j(i-1)} 
\left[\begin{array}{c} i  \\
j \end{array} \right] 
\frac{(\eta_z,\eta^*_y)}{(\eta_x, \eta^*_y)}
$
\\
$\lbrack x \rbrack^{inv}_{col} \to
\lbrack z \rbrack_{col}$
&
$
(-1)^{d-i} q^{(d-i)(d-j-1)} 
\left[\begin{array}{c} d-j  \\
i-j \end{array} \right] 
\frac{(\eta_y,\eta^*_x)}{(\eta_y, \eta^*_z)}
$
\\
\hline
$
\lbrack y \rbrack_{row} \to
\lbrack z \rbrack^{inv}_{row}
$
&
$
(-1)^j q^{j(1-i)} 
\left[\begin{array}{c} i  \\
j \end{array} \right] 
\frac{(\eta_z,\eta^*_x)}{(\eta_y, \eta^*_x)}
$
\\
$\lbrack y \rbrack_{col} \to
\lbrack z \rbrack^{inv}_{col}
$
&
$
(-1)^{d-i} q^{(i-d)(d-j-1)} 
\left[\begin{array}{c} d-j  \\
i-j \end{array} \right] 
\frac{(\eta_x,\eta^*_y)}{(\eta_x, \eta^*_z)}
$
\\
$\lbrack y \rbrack^{inv}_{row} \to 
\lbrack x \rbrack_{row}$
&
$
(-1)^j q^{j(i-1)} 
\left[\begin{array}{c} i  \\
j \end{array} \right] 
\frac{(\eta_x,\eta^*_z)}{(\eta_y, \eta^*_z)}
$
\\
$\lbrack y \rbrack^{inv}_{col} \to
\lbrack x \rbrack_{col}$
&
$
(-1)^{d-i} q^{(d-i)(d-j-1)} 
\left[\begin{array}{c} d-j  \\
i-j \end{array} \right] 
\frac{(\eta_z,\eta^*_y)}{(\eta_z, \eta^*_x)}
$
\\
\hline
$
\lbrack z \rbrack_{row} \to
\lbrack x \rbrack^{inv}_{row}
$
&
$
(-1)^j q^{j(1-i)} 
\left[\begin{array}{c} i  \\
j \end{array} \right] 
\frac{(\eta_x,\eta^*_y)}{(\eta_z, \eta^*_y)}
$
\\
$\lbrack z \rbrack_{col} \to
\lbrack x \rbrack^{inv}_{col}
$
&
$
(-1)^{d-i} q^{(i-d)(d-j-1)} 
\left[\begin{array}{c} d-j  \\
i-j \end{array} \right] 
\frac{(\eta_y,\eta^*_z)}{(\eta_y, \eta^*_x)}
$
\\
$\lbrack z \rbrack^{inv}_{row} \to 
\lbrack y \rbrack_{row}$
&
$
(-1)^j q^{j(i-1)} 
\left[\begin{array}{c} i  \\
j \end{array} \right] 
\frac{(\eta_y,\eta^*_x)}{(\eta_z, \eta^*_x)}
$
\\
$\lbrack z \rbrack^{inv}_{col} \to
\lbrack y \rbrack_{col}$
&
$
(-1)^{d-i} q^{(d-i)(d-j-1)} 
\left[\begin{array}{c} d-j  \\
i-j \end{array} \right] 
\frac{(\eta_x,\eta^*_z)}{(\eta_x, \eta^*_y)}
$
\end{tabular}}
     \medskip

\end{theorem}
\noindent {\it Proof:}
We first verify the data for the middle third of the table.
\\
\noindent 
$\lbrack y \rbrack_{row} \to 
\lbrack z \rbrack^{inv}_{row}$.
Let $\lbrace u_i\rbrace_{i=0}^d$ denote the basis
$\lbrack y \rbrack_{row}$ for $V$, and let
 $\lbrace v_i\rbrace_{i=0}^d$ denote the basis
$\lbrack z \rbrack^{inv}_{row}$ for $V$.
By Theorem
\ref{thm:closedform},
\begin{eqnarray}
v_j = \frac{(-1)^j 
q^{-\binom{j}{2}}
 }{
\lbrack j 
\rbrack^!} \frac{
(\eta_z,\eta^*_x) }
{(\eta_y,\eta^*_x)}
\,n_x^j \eta_y
\qquad \qquad (0 \leq j \leq d).
\label{eq:vjform}
\end{eqnarray}
In line
(\ref{eq:vjform}) we evaluate the right-hand side.
We have $\eta_y = \sum_{i=0}^d u_i$
by 
Definition \ref{def:urownorm}.
By 
Lemma
\ref{lem:pp}
$n_x u_i = 
q^{-i}\lbrack i+1 \rbrack u_{i+1}$
for $0 \leq i \leq d-1$ and
$n_x u_d =0$.
Evaluating the right-hand side of
(\ref{eq:vjform})
using these comments,
we find 
that the transition matrix
$\lbrack y \rbrack_{row} \to
\lbrack z \rbrack^{inv}_{row}$ is as claimed.
\\
\noindent
$\lbrack y \rbrack_{col} \to 
\lbrack z \rbrack^{inv}_{col}$. 
Compute the product of  transition matrices
\begin{eqnarray*}
\lbrack y \rbrack_{col} \to 
\lbrack y \rbrack_{row} \to 
\lbrack z \rbrack^{inv}_{row} \to 
\lbrack z \rbrack^{inv}_{col}.
\end{eqnarray*}
In this product the first and last factors
are from Theorem
\ref{thm:diagtrans}, and the middle factor is 
from earlier in this proof.
\\
\noindent 
$\lbrack y \rbrack^{inv}_{row} \to 
\lbrack x \rbrack_{row}$. 
Let $\lbrace u_i\rbrace_{i=0}^d$ denote the basis
$\lbrack y \rbrack^{inv}_{row}$ for $V$, and let
 $\lbrace v_i\rbrace_{i=0}^d$ denote the basis
$\lbrack x \rbrack_{row}$ for $V$.
By
Theorem
\ref{thm:closedform},
\begin{eqnarray}
v_j = \frac{
q^{\binom{j}{2}}
 }{
\lbrack j 
\rbrack^!} 
\frac{(\eta_x,\eta^*_z) }{
(\eta_y,\eta^*_z)}
\,n_z^j 
\eta_y
\qquad \qquad (0 \leq j \leq d).
\label{eq:vjform2}
\end{eqnarray}
In line
(\ref{eq:vjform2}) we evaluate the right-hand side.
We have $\eta_y = \sum_{i=0}^d u_i$
by 
Definition \ref{def:urownorm}.
By 
Lemma
\ref{lem:pp}
$n_z u_i = 
-q^i\lbrack i+1 \rbrack u_{i+1}$
for $0 \leq i \leq d-1$ and
$n_z u_d =0$.
Evaluating the right-hand side of
(\ref{eq:vjform2})
using these comments,
we find 
that the transition matrix
$\lbrack y \rbrack^{inv}_{row}  \to
\lbrack x \rbrack_{row}$ is as claimed.
\\
\noindent
$\lbrack y \rbrack^{inv}_{col}  \to
\lbrack x \rbrack_{col}$.
Compute the product of transition matrices 
\begin{eqnarray*}
\lbrack y \rbrack^{inv}_{col}
\to
\lbrack y \rbrack^{inv}_{row}
\to
\lbrack x \rbrack_{row}
\to
\lbrack x \rbrack_{col}.
\end{eqnarray*}
In this product the first and third factors are from Theorem
\ref{thm:diagtrans}, and the middle factor is 
from earlier in this proof.
\\
\noindent We have now verified the data for the middle
third of the table. To verify the rest of the table use
Lemma
\ref{lem:om}.
\hfill $\Box$ \\

\begin{theorem} 
In the table below we display some transition matrices
between bases for $V^*$.
Each transition matrix is lower triangular.
For $0 \leq j \leq i \leq d$ the $(i,j)$-entry is given.
\medskip

\centerline{
\begin{tabular}[t]{c|c}
   {\rm transition matrix}
   &  {\rm $(i,j)$-entry  for $0 \leq j\leq i\leq d$ }
   \\ \hline  \hline
$
\lbrack x \rbrack_{row} \to
\lbrack y \rbrack^{inv}_{row}
$
&
$
(-1)^j q^{j(i-1)} 
\left[\begin{array}{c} i  \\
j \end{array} \right] 
\frac{(\eta_z,\eta^*_y)}{(\eta_z, \eta^*_x)}
$
\\
$\lbrack x \rbrack_{col} \to
\lbrack y \rbrack^{inv}_{col}
$
&
$
(-1)^{d-i} q^{(d-i)(d-j-1)} 
\left[\begin{array}{c} d-j  \\
i-j \end{array} \right] 
\frac{(\eta_x,\eta^*_z)}{(\eta_y, \eta^*_z)}
$
\\
$\lbrack x \rbrack^{inv}_{row} \to 
\lbrack z \rbrack_{row}$
&
$
(-1)^j q^{j(1-i)} 
\left[\begin{array}{c} i  \\
j \end{array} \right] 
\frac{(\eta_y,\eta^*_z)}{(\eta_y, \eta^*_x)}
$
\\
$\lbrack x \rbrack^{inv}_{col} \to
\lbrack z \rbrack_{col}$
&
$
(-1)^{d-i} q^{(i-d)(d-j-1)} 
\left[\begin{array}{c} d-j  \\
i-j \end{array} \right] 
\frac{(\eta_x,\eta^*_y)}{(\eta_z, \eta^*_y)}
$
\\
\hline
$
\lbrack y \rbrack_{row} \to
\lbrack z \rbrack^{inv}_{row}
$
&
$
(-1)^j q^{j(i-1)} 
\left[\begin{array}{c} i  \\
j \end{array} \right] 
\frac{(\eta_x,\eta^*_z)}{(\eta_x, \eta^*_y)}
$
\\
$\lbrack y \rbrack_{col} \to
\lbrack z \rbrack^{inv}_{col}
$
&
$
(-1)^{d-i} q^{(d-i)(d-j-1)} 
\left[\begin{array}{c} d-j  \\
i-j \end{array} \right] 
\frac{(\eta_y,\eta^*_x)}{(\eta_z, \eta^*_x)}
$
\\
$\lbrack y \rbrack^{inv}_{row} \to 
\lbrack x \rbrack_{row}$
&
$
(-1)^j q^{j(1-i)} 
\left[\begin{array}{c} i  \\
j \end{array} \right] 
\frac{(\eta_z,\eta^*_x)}{(\eta_z, \eta^*_y)}
$
\\
$\lbrack y \rbrack^{inv}_{col} \to
\lbrack x \rbrack_{col}$
&
$
(-1)^{d-i} q^{(i-d)(d-j-1)} 
\left[\begin{array}{c} d-j  \\
i-j \end{array} \right] 
\frac{(\eta_y,\eta^*_z)}{(\eta_x, \eta^*_z)}
$
\\
\hline
$
\lbrack z \rbrack_{row} \to
\lbrack x \rbrack^{inv}_{row}
$
&
$
(-1)^j q^{j(i-1)} 
\left[\begin{array}{c} i  \\
j \end{array} \right] 
\frac{(\eta_y,\eta^*_x)}{(\eta_y, \eta^*_z)}
$
\\
$\lbrack z \rbrack_{col} \to
\lbrack x \rbrack^{inv}_{col}
$
&
$
(-1)^{d-i} q^{(d-i)(d-j-1)} 
\left[\begin{array}{c} d-j  \\
i-j \end{array} \right] 
\frac{(\eta_z,\eta^*_y)}{(\eta_x, \eta^*_y)}
$
\\
$\lbrack z \rbrack^{inv}_{row} \to 
\lbrack y \rbrack_{row}$
&
$
(-1)^j q^{j(1-i)} 
\left[\begin{array}{c} i  \\
j \end{array} \right] 
\frac{(\eta_x,\eta^*_y)}{(\eta_x, \eta^*_z)}
$
\\
$\lbrack z \rbrack^{inv}_{col} \to
\lbrack y \rbrack_{col}$
&
$
(-1)^{d-i} q^{(i-d)(d-j-1)} 
\left[\begin{array}{c} d-j  \\
i-j \end{array} \right] 
\frac{(\eta_z,\eta^*_x)}{(\eta_y, \eta^*_x)}
$
\end{tabular}}
     \medskip

\end{theorem}
In Theorem
\ref{thm:lowert}
replace $q$ by $q^{-1}$, and also
replace 
$(\eta_u,\eta^*_v)$ by
$(\eta_v,\eta^*_u)$
for distinct 
$u,v \in \lbrace x,y,z\rbrace$.
\hfill $\Box$ \\


\section{Rotators }

\noindent In this section we discuss the mathematics
involving the rotators
from 
Definition
\ref{def:rotator}.

\begin{proposition}
Given a rotator for $V$,
the inverse of the adjoint
is a rotator for $V^*$.
\end{proposition}
\noindent {\it Proof:}
Let $R $ denote the rotator for $V$ in question,
and note that $R$ satisfies
(\ref{eq:rotate}). In these equations apply
the adjoint map to each side. The result shows that on $V^*$,
\begin{eqnarray*}
(R^{adj})^{-1} x R^{adj} = y,
\qquad \qquad 
(R^{adj})^{-1} y R^{adj} = z,
\qquad \qquad 
(R^{adj})^{-1} z R^{adj} = x.
\end{eqnarray*}
Let $\Psi$ denote the inverse of 
$R^{adj}$.
In terms of $\Psi $ the above equations become
\begin{eqnarray*}
\Psi x \Psi^{-1} = y,
\qquad \qquad 
\Psi y \Psi^{-1} = z,
\qquad \qquad 
\Psi z \Psi^{-1} = x.
\end{eqnarray*}
Therefore $\Psi$ is a rotator for $V^*$.
\hfill $\Box$ \\

\begin{definition} 
\label{lem:pmat}
Define $P_q\in 
{\rm Mat}_{d+1}(\F)$
 to have the following
 $(i,j)$-entry 
for $0 \leq i,j\leq d$.
For  $i+j<d$ this entry is 0. For 
 $i+j\geq d$ this entry is
\begin{eqnarray*}
(-1)^{d-j}q^{(d-j)(1-i)} 
\left[\begin{array}{c} i  \\
d-j \end{array} \right] .
\end{eqnarray*}
\end{definition}

\begin{theorem}
\label{thm:omegaP}
In the table below we display
some transition matrices between bases for $V$.
\medskip

\centerline{
\begin{tabular}[t]{c|c}
   {\rm transition }
   &  {\rm transition matrix}
   \\ \hline  \hline
$\lbrack x \rbrack_{row} \to
\lbrack y \rbrack_{row}$
& 
$
P_q
\frac{(\eta_y, \eta^*_z)}
{
(\eta_x, \eta^*_z)
}
$
\\
$\lbrack x \rbrack_{col} \to
\lbrack y \rbrack_{col}$
& 
$
P^t_q
\frac{(\eta_z, \eta^*_x)}
{
(\eta_z, \eta^*_y)
}
$
\\
$\lbrack x \rbrack^{inv}_{row} \to
\lbrack y \rbrack^{inv}_{row}$
& 
$
ZP_qZ
\frac{(\eta_y, \eta^*_z)}
{
(\eta_x, \eta^*_z)
}
$
\\
$\lbrack x \rbrack^{inv}_{col} \to
\lbrack y \rbrack^{inv}_{col}$
& 
$
ZP^t_qZ
\frac{(\eta_z, \eta^*_x)}
{
(\eta_z, \eta^*_y)
}
$
\\
\hline
$\lbrack y \rbrack_{row} \to
\lbrack z \rbrack_{row}$
& 
$
P_q
\frac{(\eta_z, \eta^*_x)}
{
(\eta_y, \eta^*_x)
}
$
\\
$\lbrack y \rbrack_{col} \to
\lbrack z \rbrack_{col}$
& 
$
P^t_q
\frac{(\eta_x, \eta^*_y)}
{
(\eta_x, \eta^*_z)
}
$
\\
$\lbrack y \rbrack^{inv}_{row} \to
\lbrack z \rbrack^{inv}_{row}$
& 
$
ZP_qZ
\frac{(\eta_z, \eta^*_x)}
{
(\eta_y, \eta^*_x)
}
$
\\
$\lbrack y \rbrack^{inv}_{col} \to
\lbrack z \rbrack^{inv}_{col}$
& 
$
ZP^t_qZ
\frac{(\eta_x, \eta^*_y)}
{
(\eta_x, \eta^*_z)
}
$
\\
\hline
$\lbrack z \rbrack_{row} \to
\lbrack x \rbrack_{row}$
& 
$
P_q
\frac{(\eta_x, \eta^*_y)}
{
(\eta_z, \eta^*_y)
}
$
\\
$\lbrack z \rbrack_{col} \to
\lbrack x \rbrack_{col}$
& 
$
P^t_q
\frac{(\eta_y, \eta^*_z)}
{
(\eta_y, \eta^*_x)
}
$
\\
$\lbrack z \rbrack^{inv}_{row} \to
\lbrack x \rbrack^{inv}_{row}$
& 
$
ZP_qZ
\frac{(\eta_x, \eta^*_y)}
{
(\eta_z, \eta^*_y)
}
$
\\
$\lbrack z \rbrack^{inv}_{col} \to
\lbrack x \rbrack^{inv}_{col}$
& 
$
ZP^t_qZ
\frac{(\eta_y, \eta^*_z)}
{
(\eta_y, \eta^*_x)
}
$
\end{tabular}}
     \medskip

\end{theorem}
\noindent {\it Proof:}
We first verify the data for the middle third of the table.
\\
\noindent 
$\lbrack y \rbrack_{row} \to
\lbrack z \rbrack_{row}$.
Compute the  product of transition
matrices
\begin{eqnarray*}
\lbrack y \rbrack_{row} \to 
\lbrack z \rbrack^{inv}_{row}
\to
\lbrack z \rbrack_{row}.
\end{eqnarray*}
In this product the first factor is from
Theorem
\ref{thm:lowert} and the second factor is $Z$.
\\
\noindent 
$\lbrack y \rbrack_{col} \to
\lbrack z \rbrack_{col}$.
Compute the  product of transition
matrices
\begin{eqnarray*}
\lbrack y \rbrack_{col} \to 
\lbrack z \rbrack^{inv}_{col}
\to
\lbrack z \rbrack_{col}.
\end{eqnarray*}
In this product the first factor is from
Theorem
\ref{thm:lowert} and the second factor is $Z$.
\\
\noindent
$\lbrack y \rbrack^{inv}_{row} \to
\lbrack z \rbrack^{inv}_{row}$.
Conjugate the transition
matrix
$\lbrack y \rbrack_{row} \to
\lbrack z \rbrack_{row}$ by $Z$.
\\
\noindent 
$\lbrack y \rbrack^{inv}_{col} \to
\lbrack z \rbrack^{inv}_{col}$.
Conjugate the transition
matrix
$\lbrack y \rbrack_{col} \to
\lbrack z \rbrack_{col}$ by $Z$.
\\
\noindent We have now verified the data for
the middle third of the table. To verify
the rest of the table use Lemma
\ref{lem:om}.
\hfill $\Box$ \\

\begin{theorem}
\label{thm:omegaPdual}
In the table below we display
some transition matrices between bases for $V^*$.
\medskip

\centerline{
\begin{tabular}[t]{c|c}
   {\rm transition }
   &  {\rm transition matrix}
   \\ \hline  \hline
$\lbrack x \rbrack_{row} \to
\lbrack y \rbrack_{row}$
& 
$
P_{q^{-1}}
\frac{(\eta_z, \eta^*_y)}
{
(\eta_z, \eta^*_x)
}
$
\\
$\lbrack x \rbrack_{col} \to
\lbrack y \rbrack_{col}$
& 
$
P^t_{q^{-1}}
\frac{(\eta_x, \eta^*_z)}
{
(\eta_y, \eta^*_z)
}
$
\\
$\lbrack x \rbrack^{inv}_{row} \to
\lbrack y \rbrack^{inv}_{row}$
& 
$
ZP_{q^{-1}}Z
\frac{(\eta_z, \eta^*_y)}
{
(\eta_z, \eta^*_x)
}
$
\\
$\lbrack x \rbrack^{inv}_{col} \to
\lbrack y \rbrack^{inv}_{col}$
& 
$
ZP^t_{q^{-1}}Z
\frac{(\eta_x, \eta^*_z)}
{
(\eta_y, \eta^*_z)
}
$
\\
\hline
$\lbrack y \rbrack_{row} \to
\lbrack z \rbrack_{row}$
& 
$
P_{q^{-1}}
\frac{(\eta_x, \eta^*_z)}
{
(\eta_x, \eta^*_y)
}
$
\\
$\lbrack y \rbrack_{col} \to
\lbrack z \rbrack_{col}$
& 
$
P^t_{q^{-1}}
\frac{(\eta_y, \eta^*_x)}
{
(\eta_z, \eta^*_x)
}
$
\\
$\lbrack y \rbrack^{inv}_{row} \to
\lbrack z \rbrack^{inv}_{row}$
& 
$
ZP_{q^{-1}}Z
\frac{(\eta_x, \eta^*_z)}
{
(\eta_x, \eta^*_y)
}
$
\\
$\lbrack y \rbrack^{inv}_{col} \to
\lbrack z \rbrack^{inv}_{col}$
& 
$
ZP^t_{q^{-1}}Z
\frac{(\eta_y, \eta^*_x)}
{
(\eta_z, \eta^*_x)
}
$
\\
\hline
$\lbrack z \rbrack_{row} \to
\lbrack x \rbrack_{row}$
& 
$
P_{q^{-1}}
\frac{(\eta_y, \eta^*_x)}
{
(\eta_y, \eta^*_z)
}
$
\\
$\lbrack z \rbrack_{col} \to
\lbrack x \rbrack_{col}$
& 
$
P^t_{q^{-1}}
\frac{(\eta_z, \eta^*_y)}
{
(\eta_x, \eta^*_y)
}
$
\\
$\lbrack z \rbrack^{inv}_{row} \to
\lbrack x \rbrack^{inv}_{row}$
& 
$
ZP_{q^{-1}}Z
\frac{(\eta_y, \eta^*_x)}
{
(\eta_y, \eta^*_z)
}
$
\\
$\lbrack z \rbrack^{inv}_{col} \to
\lbrack x \rbrack^{inv}_{col}$
& 
$
ZP^t_{q^{-1}}Z
\frac{(\eta_z, \eta^*_y)}
{
(\eta_x, \eta^*_y)
}
$
\end{tabular}}
     \medskip
\end{theorem}
\noindent {\it Proof:}
In Theorem
\ref{thm:omegaP},
replace $q$ by $q^{-1}$
and also
replace $(\eta_u,\eta^*_v)$
by
 $(\eta_v,\eta^*_u)$
for distinct  
$u,v \in \lbrace x,y,z\rbrace$.
\hfill $\Box$ \\

\begin{lemma} 
For the matrix $P_q$ in Definition
\ref{lem:pmat},
\begin{eqnarray*}
P_q^3 = (-1)^d q^{-d(d-1)}I.
\end{eqnarray*}
\end{lemma}
\noindent {\it Proof:}
Consider the bases
$
\lbrack x \rbrack_{row}
$,
$
\lbrack y \rbrack_{row}$,
$
\lbrack z \rbrack_{row}
$
for $V$.
The identity matrix is equal to the product
of transition matrices
\begin{eqnarray*}
\lbrack x \rbrack_{row} \to
\lbrack y \rbrack_{row} \to
\lbrack z \rbrack_{row} \to
\lbrack x \rbrack_{row}.
\end{eqnarray*}
In this product the three factors are given
in Theorem
\ref{thm:omegaP}. Simplify the product
using Proposition
\ref{prop:sixfactors} to get the result.
\hfill $\Box$ \\

\begin{lemma}
\label{lem:inverseP}
We have $P^{-1}_q = Z P_{q^{-1}}Z$.
For $0 \leq i,j\leq d$ the
$(i,j)$-entry of $P^{-1}_q$ is given as follows.
For $i+j>d$ this entry is 0. For
 $i+j\leq d$ this entry is
\begin{eqnarray*}
(-1)^{j}q^{j(d-i-1)} 
\left[\begin{array}{c} d-i  \\
j \end{array} \right] .
\end{eqnarray*}
\end{lemma}
\noindent {\it Proof:} Using
Definition
\ref{lem:pmat}
one checks that the
entries of $ZP_{q^{-1}}Z$ are as shown.
It remains to verify that
 $P^{-1}_q=ZP_{q^{-1}}Z$.
To this end
we consider the transition matrices between
some bases for $V$.
Let $T$ denote the transition matrix 
$\lbrack z \rbrack_{row} \to
\lbrack y \rbrack_{row}$. On
one hand, $T$
is the inverse of the transition matrix 
$\lbrack y \rbrack_{row} \to
\lbrack z \rbrack_{row}$. The transition matrix
for $\lbrack y \rbrack_{row} \to
\lbrack z \rbrack_{row}$
 can be found
in Theorem
\ref{thm:omegaP}.
On the other hand, $T$
is the product of transition matrices
\begin{eqnarray*}
\lbrack z \rbrack_{row} \to
\lbrack z \rbrack^{inv}_{row}
\to
\lbrack y \rbrack_{row}.
\end{eqnarray*}
In this product the first factor is
$Z$ and the second factor is from
Theorem
\ref{thm:lowert}. By these comments one verifies
that
 $P^{-1}_q=ZP_{q^{-1}}Z$ after a brief computation.
\hfill $\Box$ \\

\begin{theorem}
\label{thm:rotP}
There exists a rotator for $V$
that is represented by
$P_q$ with
respect to each of the bases
$\lbrack x \rbrack_{row}$,
$\lbrack y\rbrack_{row}$,
$\lbrack z \rbrack_{row}$ for $V$.
Moreover, there exists a rotator
for $V^*$ that is represented by
$P_{q^{-1}}$ with
respect to each of the bases
$\lbrack x \rbrack_{row}$,
$\lbrack y\rbrack_{row}$,
$\lbrack z \rbrack_{row}$ for $V^*$.
\end{theorem}
\noindent {\it Proof:}
We first verify our assertion about $V$.
By Theorem
\ref{thm:omegaP},
each of the transition matrices
\begin{eqnarray*}
\lbrack x \rbrack_{row}\to
\lbrack y \rbrack_{row},
\qquad \quad 
\lbrack y \rbrack_{row}\to
\lbrack z \rbrack_{row},
\qquad \quad 
\lbrack z \rbrack_{row}\to
\lbrack x \rbrack_{row}
\end{eqnarray*}
is contained in  $\F P_q$.
By Lemma \ref{lem:om} there exists a rotator for
$V$. Denote this rotator by $R $.
For $\xi \in \lbrace x,y,z\rbrace$ let
$T_\xi$ denote the matrix that represents
$R$ with respect to  
the basis $\lbrack \xi \rbrack_{row}$.
By our initial comment and 
since $R$ is a rotator,
there exists $0 \not=\alpha_\xi \in \F$ such that
$T_\xi =\alpha_\xi P_q$.
By our initial comment and linear
algebra,
\begin{eqnarray*}
P_q^{-1}T_x P_q = T_y,
\qquad \qquad
P_q^{-1}T_y P_q = T_z,
\qquad \qquad
P_q^{-1}T_z P_q = T_x.
\end{eqnarray*}
Therefore 
$\alpha_x=\alpha_y=\alpha_z$. Let $\alpha$ denote this common
value and note that
 $R/\alpha$ is the desired rotator
for $V$.
We have verified our assertion about $V$.
The  assertion
about $V^*$ is similarly verified.
\hfill $\Box$ \\

\begin{definition}\rm
Let $\mathcal R$ denote the rotator
for $V$ or $V^*$ referred to in Theorem
\ref{thm:rotP}.
\end{definition}

\begin{theorem}
\label{thm:rot12}
In the table below we display the matrices that
represent $\mathcal R$ with respect to the twelve
bases for $V$ from
{\rm (\ref{eq:basisr1})--(\ref{eq:basisr4})}.

\medskip

\centerline{
\begin{tabular}[t]{c|c }
{\rm basis}
& 
   {\rm matrix rep. $\mathcal R$ }
   \\ \hline 
   $\lbrack \xi \rbrack_{row}$ &
    $P_q$ 
 \\ 
   $\lbrack \xi \rbrack_{col}$ &
    $P^t_q$ 
 \\ 
   $\lbrack \xi \rbrack^{inv}_{row}$ &
    $ZP_qZ$ 
 \\ 
   $\lbrack \xi \rbrack^{inv}_{col}$ &
    $ZP^t_qZ$ 
     \end{tabular}}

\medskip
\noindent In the above table $\xi \in \lbrace x,y,z\rbrace$.
\end{theorem}
\noindent {\it Proof:}
By Theorem
\ref{thm:rotP} the matrix $P_q$ represents $\mathcal R$
with respect to 
   $\lbrack \xi \rbrack_{row}$.
We now show that $P_q^t$ represents $\mathcal R$
with respect to
   $\lbrack \xi \rbrack_{col}$.
Let $D_\xi$ denote the transition matrix 
   $\lbrack \xi \rbrack_{row} \to
   \lbrack \xi \rbrack_{col}$.
By linear algebra
the matrix $D_\xi^{-1} P_q D_\xi$
represents $\mathcal R$ with respect to
   $\lbrack \xi \rbrack_{col}$.
The entries of $D_\xi $ are given in
Theorem
\ref{thm:diagtrans}. By this data 
and Definition
\ref{lem:pmat},
 $D_\xi^{-1} P_q D_\xi=P_q^t$.
Therefore $P_q^t$ represents $\mathcal R$
with respect to
   $\lbrack \xi \rbrack_{col}$.
We have verified the first two rows of the table.
The remaining rows are readily verified.
\hfill $\Box$ \\

\begin{theorem}
\label{thm:rot12dual}
In the table below we display the matrices that
represent $\mathcal R$ with respect to the twelve
bases for $V^*$ from
{\rm (\ref{eq:basisr1})--(\ref{eq:basisr4})}.

\medskip

\centerline{
\begin{tabular}[t]{c|c }
{\rm basis}
& 
   {\rm matrix rep. $\mathcal R$ }
   \\ \hline 
   $\lbrack \xi \rbrack_{row}$ &
    $P_{q^{-1}}$ 
 \\ 
   $\lbrack \xi \rbrack_{col}$ &
    $P^t_{q^{-1}}$ 
 \\ 
   $\lbrack \xi \rbrack^{inv}_{row}$ &
    $ZP_{q^{-1}}Z$ 
 \\ 
   $\lbrack \xi \rbrack^{inv}_{col}$ &
    $ZP^t_{q^{-1}}Z$ 
     \end{tabular}}

\medskip
\noindent In the above table $\xi \in \lbrace x,y,z\rbrace$.
\end{theorem}
\noindent {\it Proof:}
In Theorem \ref{thm:rot12} replace $q$ by $q^{-1}$ and invoke
Theorem
\ref{thm:rotP}.
\hfill $\Box$ \\

\begin{theorem} For the rotator $\mathcal R$ of $V$, the
 inverse of the adjoint 
is the rotator $\mathcal R$ for $V^*$.
\end{theorem}
\noindent {\it Proof:} 
By Theorem
\ref{thm:rot12}
the matrix $P_q$ represents 
$\mathcal R$ with respect to the basis
$\lbrack x \rbrack_{row} $ for $V$.
By Lemma
\ref{lem:whatdualwhat2}
the basis 
$\lbrack x \rbrack^{inv}_{col} $ for $V^*$
is dual to the basis
$\lbrack x \rbrack_{row} $ for $V$.
Therefore 
$P_q^t$ represents 
 ${\mathcal R}^{adj}$ with respect to the basis
$\lbrack x \rbrack^{inv}_{col} $ for $V^*$.
Therefore
$(P_q^t)^{-1}$ represents 
 $({\mathcal R}^{adj})^{-1}$ with respect to the basis
$\lbrack x \rbrack^{inv}_{col} $ for $V^*$.
By Theorem
\ref{thm:rot12dual}
the matrix $ZP^t_{q^{-1}}Z$ represents 
$\mathcal R$ with respect to the basis
$\lbrack x \rbrack^{inv}_{col} $ for $V^*$.
We have
$P^{-1}_q = Z P_{q^{-1}}Z$ by
Lemma
\ref{lem:inverseP} so
$(P_q^t)^{-1} =
ZP^t_{q^{-1}}Z$.
The result follows.
\hfill $\Box$ \\

\section
{A characterization of $y$ and $n_y$}

\noindent Recall the elements
$x,y,z$ and $n_x,n_y,n_z$ of
$U_q(\mathfrak{sl}_2)$.
In this section we 
 characterize 
 $y$ and $n_y$ using
the $U_q(\mathfrak{sl}_2)$-module $V$.
Similar characterizations apply to
$x, z$ and
$n_x, n_z$. We will be using
 Definition
\ref{def:lowering}.

\begin{theorem} Given $\phi \in 
{\rm End}(V)$. Then
$\phi \in \F n_y$ if and only if both
\begin{enumerate}
\item[\rm (i)]
$\phi$ is lowering for the decomposition
$\lbrack x \rbrack$ of $V$;
\item[\rm (ii)]
$\phi$ is 
 raising for the decomposition
$\lbrack z \rbrack$ of $V$.
\end{enumerate}
\end{theorem}
\noindent {\it Proof:}  First assume that 
$\phi \in \F n_y$.
Then $\phi$ satisfies the above conditions
(i), (ii) by Theorem
\ref{thm:nxnynz}.
Conversely, assume that $\phi$ satisfies
(i), (ii). We  show
$\phi \in \F n_y$. To avoid trivialities assume 
$\phi \not=0$.
Let $\lbrace V_i\rbrace_{i=0}^d$ denote
the decomposition
$\lbrack x \rbrack$ of $V$.
Let $\lbrace v_i\rbrace_{i=0}^d$ denote the
basis
$\lbrack x \rbrack_{row}$ for $V$.
So $V_i$ has basis $v_i$ for $0 \leq i \leq d$.
Recall
that $\eta_x = \sum_{i=0}^d v_i$ is a basis for
 component 0 of the decomposition
$\lbrack z \rbrack$ of $V$.
By assumption $\phi$ is raising for
$\lbrack z \rbrack$. We assume $\phi\not=0$ so $d\geq 1$.
Moreover 
$\phi \eta_x$ is contained in component 1 of 
$\lbrack z \rbrack$. 
By Theorem
\ref{thm:closedform},
$n_y \eta_x$ is a basis for
component 1 of 
$\lbrack z \rbrack$.
Therefore there exists $\alpha \in \F$ such that
$\phi \eta_x =\alpha n_y \eta_x$. By this and
$\eta_x = \sum_{i=0}^d v_i$, 
\begin{eqnarray}
0 =
\sum_{i=0}^d (\phi-\alpha n_y)v_i.
\label{eq:exp0}
\end{eqnarray}
Each of $\phi$, $n_y$ is lowering for
 $\lbrack x \rbrack$. Therefore
$\phi-\alpha n_y$ is lowering for
 $\lbrack x \rbrack$.
Therefore in (\ref{eq:exp0}) the $i$th summand is 
zero for $i=0$ and
contained
in $V_{i-1}$
for $1 \leq i \leq d$.
Now since the sum $\sum_{j=0}^d V_j$ is direct,
 in (\ref{eq:exp0}) the $i$th summand is zero for
$0 \leq i \leq d$.
Thus
$\phi-\alpha n_y$ vanishes on each vector in the basis
$\lbrace v_i\rbrace_{i=0}^d$ for $V$.
Therefore $\phi=\alpha n_y$, 
 so $\phi \in \F n_y$ as desired.
 The 
 result follows.
\hfill $\Box$ \\

\begin{theorem} Given $\phi \in 
{\rm End}(V)$. Then
$\phi \in \F y + \F 1$  if and only if both
\begin{enumerate}
\item[\rm (i)]
$\phi$ is quasi-raising for the decomposition
$\lbrack x \rbrack$ of $V$;
\item[\rm (ii)]
$\phi$ is 
 quasi-lowering for the decomposition
$\lbrack z \rbrack$ of $V$.
\end{enumerate}
\end{theorem}
\noindent {\it Proof:} 
First assume that
$\phi \in \F y + \F 1$. Then
$\phi$ satisfies the above conditions (i), (ii) by
Theorem
\ref{thm:xyzV}. Conversely, assume that
$\phi$ satisfies (i), (ii). We show
$\phi \in \F y + \F 1$. 
To avoid trivialities assume $\phi \not=0$.
Let $\lbrace V_i\rbrace_{i=0}^d$
denote the decomposition $\lbrack y \rbrack$
of $V$. We show that $\phi V_i \subseteq V_i$
for $0 \leq i \leq d$. Let $i$ be
given. By Lemma
\ref{lem:decompdesc}, $V_i$
is equal to the intersection of
  $n_z^{d-i}V$
  and
  $n_x^iV$.
 By Lemma
\ref{lem:nnnpre}(iii),
  $n_z^{d-i}V$ is the sum
  of components $d-i,d-i+1,\ldots,d$ for the decomposition
  $\lbrack x \rbrack $ of $V$. By assumption
  $\phi$ is quasi-raising for
  $\lbrack x \rbrack $.
  Therefore 
  $n_z^{d-i}V$ is $\phi$-invariant.
By Lemma \ref{lem:nnnpre}(i),
   $n_x^iV$ is the sum of components
  $0,1,\ldots, d-i$ for 
  the decomposition
  $\lbrack z \rbrack $ of $V$. By assumption
  $\phi$ is quasi-lowering for
  $\lbrack z \rbrack$.
  Therefore $n_x^iV$ is $\phi$-invariant.
 By these comments 
\begin{eqnarray*}
 \phi V_i =
\phi ( 
  n_z^{d-i}V \cap
  n_x^iV
)
\subseteq 
\phi ( 
  n_z^{d-i}V) \cap
  \phi(n_x^iV)
\subseteq 
  n_z^{d-i}V \cap
   n_x^iV
   = V_i.
\end{eqnarray*}
We have shown that $\phi V_i \subseteq V_i$
for $0 \leq i \leq d$.
The $\lbrace V_i\rbrace_{i=0}^d$ are the
eigenspaces for $y$ on $V$, so
$\phi$ commutes with $y$ 
  on $V$.
By this and since $y$ is
multiplicity-free on $V$, we see that
$\phi$ is contained in the subalgebra of
${\rm End}(V)$ generated by $y$.
This subalgebra has basis
$\lbrace y^i\rbrace_{i=0}^d$.
This subalgebra has another basis 
$\lbrace y_i\rbrace_{i=0}^d$ where
\begin{eqnarray*}
y_i = (y-q^{-d})(y-q^{2-d})\cdots (y-q^{2i-2-d})
\qquad \qquad (0 \leq i \leq d).
\end{eqnarray*}
By construction there exist
scalars $\lbrace \alpha_i\rbrace_{i=0}^d$ in 
$\F$ such that $\phi = \sum_{i=0}^d \alpha_i y_i$
on $V$.
Recall $\phi \not=0$ so 
$\lbrace \alpha_i\rbrace_{i=0}^d$ are not all zero.
Define $s={\rm max}\lbrace i|0 \leq i \leq d,
\alpha_i\not=0\rbrace $. We show $s\leq 1$. 
To this end we assume $s\geq 2$ and get a contradiction.
By construction
\begin{eqnarray}
\phi - \sum_{i=0}^{s-1}
\alpha_i y_i = \alpha_s y_s.
\label{eq:ys}
\end{eqnarray}
Let $\lbrace U_i\rbrace_{i=0}^d$ denote
the decomposition
$\lbrack x \rbrack$ of $V$.
Referring to equation 
(\ref{eq:ys}), we will apply each side to 
$U_0$.
By assumption $\phi$ is quasi-raising for
$\lbrack x \rbrack$.
Therefore $\phi U_0 \subseteq U_0+U_1$.
By Theorem
\ref{thm:xyzV} $y_i U_0 = U_i$
for $0 \leq i \leq d$.
Now for the equation 
(\ref{eq:ys}), apply each side to 
$U_0$ and consider the image. For the
left-hand side the image is contained in
$\sum_{i=0}^{s-1} U_i$. For the
right-hand side the image is $U_s$.
This is a contradiction, so $s\leq 1$. 
Therefore
$\phi \in \F y + \F 1$, as desired.
\hfill $\Box$ \\


\section
{A characterization of
$U_q(\mathfrak{sl}_2)$ 
}

\noindent In this section we give a
characterization of
$U_q(\mathfrak{sl}_2)$ in its equitable presentation.
This characterization extends some work of
Darren Funk-Neubauer
\cite{neubauer} concerning bidiagonal pairs of linear
transformations.
In order to motivate our result, we
 consider some implications of
Theorem \ref{thm:xyzV}.
Referring to the
$U_q(\mathfrak{sl}_2)$-module $V$ from that theorem,
let 
basis 1 (resp. basis 2) (resp. basis 3)
denote a basis for $V$ that induces the
decomposition 
$\lbrack x \rbrack$
(resp. $\lbrack y \rbrack$)
(resp. $\lbrack z \rbrack$) for $V$. 
On these bases 
 $x,y,z$ act as follows:

\medskip

\centerline{
\begin{tabular}[t]{c| c c c}
& 
   {\rm matrix rep. $x$ }
   &
   {\rm matrix rep. $y$ }
   &
   {\rm matrix rep. $z$ }
   \\ \hline 
   {\rm basis 1} &
    {\rm diagonal}
     & 
    {\rm lower bidiagonal}
     & 
    {\rm upper bidiagonal}
 \\ 
  {\rm basis 2} &
    {\rm upper bidiagonal}
    &
    {\rm diagonal}
     & 
    {\rm lower bidiagonal}
   \\ 
   {\rm basis 3} &
    {\rm lower bidiagonal}
     & 
    {\rm upper bidiagonal}
    &
    {\rm diagonal}
     \end{tabular}}

\medskip

\noindent The above  pattern appears not only for
irreducible $U_q(\mathfrak{sl}_2)$-modules. 
It also appears for
irreducible 
$\mathfrak{sl}_2$-modules \cite[Section~8]{bt},
as we now explain.

\begin{definition}
\label{def:sl2def}
\rm
\cite[Line~(2.2)]{bt}
Assume that $\F$ has characteristic 0.
For the Lie algebra
$\mathfrak{sl}_2$ over $\F$, the equitable basis
$x,y,z$ satisfies
\begin{eqnarray*}
\lbrack x,y\rbrack = 2x+2y,
\qquad \qquad
\lbrack y,z\rbrack = 2y+2z,
\qquad \qquad
\lbrack z,x\rbrack = 2z+2x.
\end{eqnarray*}
\end{definition}

\noindent Referring to 
Definition
\ref{def:sl2def},
let $W$ denote an irreducible 
$\mathfrak{sl}_2$-module  with dimension $d+1$.
By \cite[Section~8]{bt}
each of $x,y,z$ is multiplicity-free on $W$
with eigenvalues $\lbrace d-2i\rbrace_{i=0}^d$.
For $u \in \lbrace x,y,z\rbrace$
 define a decomposition $\lbrack u \rbrack$ of $W$
as follows. For $0 \leq i \leq d$ the $i$th component of
 $\lbrack u \rbrack$ is the eigenspace for $u$ 
on $W$ with eigenvalue $2i-d$.
Let 
basis 1 (resp. basis 2) (resp. basis 3)
denote a basis for $W$ that induces the
decomposition 
$\lbrack x \rbrack$
(resp. $\lbrack y \rbrack$)
(resp. $\lbrack z \rbrack$) for $W$. 
On these bases the
$\mathfrak{sl}_2$  elements
 $x,y,z$ act as in the above  table
\cite[Section~8]{bt}.

\begin{definition}\rm 
Let $0 \not=b\in \F$.
Let $\lbrace \alpha_i\rbrace_{i=0}^d$
denote a sequence of scalars taken from $\F$.  This seqeunce
is called {\it $b$-recurrent} whenever 
$\alpha_{i-1}\not=\alpha_i$ for 
$1 \leq i \leq d$ and 
\begin{eqnarray*}
\frac{\alpha_{i}-\alpha_{i+1}}{
\alpha_{i-1} - \alpha_i} = b
\qquad \qquad (1 \leq i \leq d-1).
\end{eqnarray*}
\end{definition}

\noindent The following theorem extends a result of 
Funk-Neubauer
\cite[Theorem~5.11]{neubauer}.

\begin{theorem}
Assume that the
field $\F$ is algebraically closed with characteristic
0.
Let $V$ denote a vector space over $\F$ with finite
positive dimension. Suppose we are given
$X,Y,Z$ in ${\rm End}(V)$.
Assume that there exist three bases for $V$ on which
 $X,Y,Z$ act as follows:

\medskip

\centerline{
\begin{tabular}[t]{c| c c c}
& 
   {\rm matrix rep. $X$ }
   &
   {\rm matrix rep. $Y$ }
   &
   {\rm matrix rep. $Z$ }
   \\ \hline 
   {\rm basis 1} &
    {\rm diagonal}
     & 
    {\rm lower bidiagonal}
     & 
    {\rm upper bidiagonal}
 \\ 
  {\rm basis 2} &
    {\rm upper bidiagonal}
    &
    {\rm diagonal}
     & 
    {\rm lower bidiagonal}
   \\ 
   {\rm basis 3} &
    {\rm lower bidiagonal}
     & 
    {\rm upper bidiagonal}
    &
    {\rm diagonal}
     \end{tabular}}
     
     \medskip
\noindent 
Then there exists $0 \not=b \in \F$ such that
for each diagonal matrix in the above table the sequence of
diagonal entries (top left to bottom right) is $b$-recurrent.
First assume $b\not=1$ and pick $q \in \F$ such that
$b=q^{-2}$.
Then 
there exists 
an irreducible  
 $U_q(\mathfrak{sl}_2)$-module  structure for $V$
such that on $V$,
\begin{eqnarray}
\label{eq:xyzaction}
x \in 
 \F X+\F I,
\qquad \qquad
y \in 
\F Y+ \F I,
\qquad \qquad
z \in 
\F Z+ \F I.
\end{eqnarray}
Next assume  $b =1$. Then
there exists  
 an irreducible  
 $\mathfrak{sl}_2$-module  structure for $V$
such that 
{\rm (\ref{eq:xyzaction})} holds 
on $V$.
\end{theorem}
\noindent {\it Proof:} 
For notational convenience, assume that the dimension of $V$ is $d+1$.
We now show that $X$
multiplicity-free.
With respect to basis 1 the matrix representing
$X$ is diagonal.
Therefore $X$ is diagonalizable on $V$.
With respect to basis 2 the matrix representing
$X$ is upper bidiagonal. 
Call this matrix $\mathbb X$.
Recall the definition of upper bidiagonal
from 
below Lemma
\ref{lem:ZBZ}.
By this definition
the matrices
${\lbrace {\mathbb X}^i \rbrace}_{i=0}^d$ are linearly
independent over $\F$. Therefore
${\lbrace X^i \rbrace}_{i=0}^d$ are linearly
independent over $\F$.
Consequently 
the minimal
polynomial of $X$ has degree $d+1$,
so $X$ has $d+1$ eigenspaces.
These eigenspaces must have
dimension 1,  so $X$ is
multiplicity-free. By a similar argument
$Y$ and $Z$ are
multiplicity-free. 
Now by the table in the theorem statement,
the maps $X,Y,Z$
act on each other's eigenspaces in
a bidiagonal fashion.
Consequently any two of $X, Y, Z$ form a bidiagonal
pair in the sense of Funk-Neubauer
\cite[Definition~2.2]{neubauer}.
Let
$\mathcal X$ (resp. $\mathcal Y$) (resp. $\mathcal Z$)
denote the  matrix  in 
${\rm Mat}_{d+1}(\F)$
that represents $X$ with respect to basis 1
(resp. 
 $Y$ with respect to basis 2)
(resp. 
 $Z$ with respect to basis 3).
Each of $\mathcal X, \mathcal Y, \mathcal Z$ is diagonal. 
By 
\cite[Theorem~5.1]{neubauer} there exists
$0 \not=b\in \F$ such that 
each sequence
of diagonal entries
$\lbrace \mathcal X_{ii}\rbrace_{i=0}^d$,
$\lbrace \mathcal Y_{ii}\rbrace_{i=0}^d$,
$\lbrace \mathcal Z_{ii}\rbrace_{i=0}^d$ is $b$-recurrent.
First assume $b\not=1$ and pick $q \in \F$
such that $b=q^{-2}$.
By the $b$-recurrence
there exist $a_1, a_2 \in \F$ such that $a_2\not=0$
and 
$\mathcal X_{ii} = a_1 + a_2 q^{d-2i}$ for $0 \leq i \leq d$.
After replacing $X$ by $(X-a_1 I)/a_2$
we obtain $\mathcal X_{ii}= q^{d-2i}$
for $0 \leq i \leq d$. Similarly adjusting $Y, Z$ we obtain
 $\mathcal Y_{ii}= q^{d-2i}$ and
 $\mathcal Z_{ii}= q^{d-2i}$
for $0 \leq i \leq d$. 
Now each of $X, Y, Z$ is multiplicity-free with
eigenvalues $\lbrace q^{d-2i} \rbrace_{i=0}^d$. 
These eigenvalues are nonzero
so $X, Y, Z$ are invertible.
Moreover by
\cite[Lemma~8.1]{neubauer},
\begin{eqnarray*}
\frac{qXY-q^{-1}YX}{q-q^{-1}}=I,
\qquad
\frac{qYZ-q^{-1}ZY}{q-q^{-1}}=I,
\qquad
\frac{qZX-q^{-1}XZ}{q-q^{-1}}=I.
\end{eqnarray*}
By these comments
 $V$
becomes a
$U_q(\mathfrak{sl}_2)$-module
on which $x,y,z$ act as $X,Y,Z$ respectively.
One checks that this 
$U_q(\mathfrak{sl}_2)$-module is irreducible.
Next assume
$b=1$.
By the $1$-recurrence 
there exist $a_1, a_2 \in \F$ such that $a_2\not=0$
and 
$\mathcal X_{ii} = a_1 + a_2 (2i-d)$ for $0 \leq i \leq d$.
After replacing $X$ by $(X-a_1 I)/a_2$
we obtain $\mathcal X_{ii}= 2i-d$
for $0 \leq i \leq d$. Similarly adjusting $Y$, $Z$ we obtain
 $\mathcal Y_{ii}= 2i-d$ and
 $\mathcal Z_{ii}= 2i-d$
for $0 \leq i \leq d$. 
Now each of $X, Y, Z$ is multiplicity-free with
eigenvalues $\lbrace 2i-d \rbrace_{i=0}^d$. 
By
\cite[Lemma~8.1]{neubauer},
\begin{eqnarray*}
XY-YX=2X+2Y,
\qquad
YZ-ZY=2Y+2Z,
\qquad
ZX-XZ=2Z+2X.
\end{eqnarray*}
Consequently 
$V$ becomes an
$\mathfrak{sl}_2$-module
on which $x,y,z$ act as $X,Y,Z$ respectively.
One  checks that this 
$\mathfrak{sl}_2$-module is irreducible.
The result follows.
\hfill $\Box$ \\


\section{Acknowledgments}
The author thanks Kazumasa Nomura
for giving this paper a close reading and offering many valuable
suggestions.

\noindent Paul Terwilliger \hfil\break
\noindent Department of Mathematics \hfil\break
\noindent University of Wisconsin \hfil\break
\noindent 480 Lincoln Drive \hfil\break
\noindent Madison, WI 53706-1388 USA \hfil\break
\noindent email: {\tt terwilli@math.wisc.edu }\hfil\break


\end{document}